\documentclass[11pt]{article}
\usepackage{amsmath, upgreek}
\usepackage{amssymb}
\usepackage{amscd}
\usepackage{xspace}
\usepackage{verbatim}
\usepackage{graphicx}	
\usepackage{color}

\newcommand{\mysection}[1]{
\section{#1}\setcounter{equation}{0}}
\title{\bf Singular solutions of some elliptic equations involving mixed absorption-reaction}
\author{{\bf Marie-Fran\c{c}oise Bidaut-V\'eron\footnote{\noindent Laboratoire de Math\'{e}matiques et Physique Th\'{e}orique,
Universit\'e de Tours, 37200 Tours, France. E-mail: veronmf@univ-tours.fr},} \\{\bf Marta Garcia-Huidobro \footnote{\noindent
Departamento de Matematicas, Pontifica Universidad Catolica de Chile
Casilla 307, Correo 2, Santiago de Chile. E-mail: mgarcia@mat.puc.cl}}\\
 {\bf Laurent V\'eron \footnote{\noindent
Laboratoire de Math\'{e}matiques et Physique Th\'{e}orique, Universit\'e de Tours, 37200 Tours, France. E-mail: veronl@univ-tours.fr}}\\[2mm]
}

\date{}
\begin{document}
 \maketitle


\newcommand{\txt}[1]{\;\text{ #1 }\;}
\newcommand{\tbf}{\textbf}
\newcommand{\tit}{\textit}
\newcommand{\tsc}{\textsc}
\newcommand{\trm}{\textrm}
\newcommand{\mbf}{\mathbf}
\newcommand{\mrm}{\mathrm}
\newcommand{\bsym}{\boldsymbol}
\newcommand{\scs}{\scriptstyle}
\newcommand{\sss}{\scriptscriptstyle}
\newcommand{\txts}{\textstyle}
\newcommand{\dsps}{\displaystyle}
\newcommand{\fnz}{\footnotesize}
\newcommand{\scz}{\scriptsize}
\newcommand{\be}{\begin{equation}}
\newcommand{\bel}[1]{\begin{equation}\label{#1}}
\newcommand{\ee}{\end{equation}}
\newcommand{\eqnl}[2]{\begin{equation}\label{#1}{#2}\end{equation}}
\newcommand{\barr}{\begin{eqnarray}}
\newcommand{\earr}{\end{eqnarray}}
\newcommand{\bars}{\begin{eqnarray*}}
\newcommand{\ears}{\end{eqnarray*}}
\newcommand{\nnu}{\nonumber \\}
\newtheorem{subn}{\name}
\renewcommand{\thesubn}{}
\newcommand{\bsn}[1]{\def\name{#1}\begin{subn}}
\newcommand{\esn}{\end{subn}}
\newtheorem{sub}{\name}[section]
\newcommand{\dn}[1]{\def\name{#1}}   
\newcommand{\bs}{\begin{sub}}
\newcommand{\es}{\end{sub}}
\newcommand{\bsl}[1]{\begin{sub}\label{#1}}
\newcommand{\bth}[1]{\def\name{Theorem}
\begin{sub}\label{t:#1}}
\newcommand{\blemma}[1]{\def\name{Lemma}
\begin{sub}\label{l:#1}}
\newcommand{\bcor}[1]{\def\name{Corollary}
\begin{sub}\label{c:#1}}
\newcommand{\bdef}[1]{\def\name{Definition}
\begin{sub}\label{d:#1}}
\newcommand{\bprop}[1]{\def\name{Proposition}
\begin{sub}\label{p:#1}}
\newcommand{\R}{\eqref}
\newcommand{\rth}[1]{Theorem~\ref{t:#1}}
\newcommand{\rlemma}[1]{Lemma~\ref{l:#1}}
\newcommand{\rcor}[1]{Corollary~\ref{c:#1}}
\newcommand{\rdef}[1]{Definition~\ref{d:#1}}
\newcommand{\rprop}[1]{Proposition~\ref{p:#1}}
\newcommand{\BA}{\begin{array}}
\newcommand{\EA}{\end{array}}
\newcommand{\BAN}{\renewcommand{\arraystretch}{1.2}
\setlength{\arraycolsep}{2pt}\begin{array}}
\newcommand{\BAV}[2]{\renewcommand{\arraystretch}{#1}
\setlength{\arraycolsep}{#2}\begin{array}}
\newcommand{\BSA}{\begin{subarray}}
\newcommand{\ESA}{\end{subarray}}
\newcommand{\BAL}{\begin{aligned}}
\newcommand{\EAL}{\end{aligned}}
\newcommand{\BALG}{\begin{alignat}}
\newcommand{\EALG}{\end{alignat}}
\newcommand{\BALGN}{\begin{alignat*}}
\newcommand{\EALGN}{\end{alignat*}}
\newcommand{\note}[1]{\textit{#1.}\hspace{2mm}}
\newcommand{\Proof}{\note{Proof}}
\newcommand{\qeda}{\hspace{10mm}\hfill $\square$}
\newcommand{\qed}{\\
${}$ \hfill $\square$}
\newcommand{\Remark}{\note{Remark}}
\newcommand{\modin}{$\,$\\[-4mm] \indent}
\newcommand{\forevery}{\quad \forall}
\newcommand{\set}[1]{\{#1\}}
\newcommand{\setdef}[2]{\{\,#1:\,#2\,\}}
\newcommand{\setm}[2]{\{\,#1\mid #2\,\}}
\newcommand{\mt}{\mapsto}
\newcommand{\lra}{\longrightarrow}
\newcommand{\lla}{\longleftarrow}
\newcommand{\llra}{\longleftrightarrow}
\newcommand{\Lra}{\Longrightarrow}
\newcommand{\Lla}{\Longleftarrow}
\newcommand{\Llra}{\Longleftrightarrow}
\newcommand{\warrow}{\rightharpoonup}
\newcommand{
\paran}[1]{\left (#1 \right )}
\newcommand{\sqbr}[1]{\left [#1 \right ]}
\newcommand{\curlybr}[1]{\left \{#1 \right \}}
\newcommand{\abs}[1]{\left |#1\right |}
\newcommand{\norm}[1]{\left \|#1\right \|}
\newcommand{
\paranb}[1]{\big (#1 \big )}
\newcommand{\lsqbrb}[1]{\big [#1 \big ]}
\newcommand{\lcurlybrb}[1]{\big \{#1 \big \}}
\newcommand{\absb}[1]{\big |#1\big |}
\newcommand{\normb}[1]{\big \|#1\big \|}
\newcommand{
\paranB}[1]{\Big (#1 \Big )}
\newcommand{\absB}[1]{\Big |#1\Big |}
\newcommand{\normB}[1]{\Big \|#1\Big \|}
\newcommand{\produal}[1]{\langle #1 \rangle}

\newcommand{\thkl}{\rule[-.5mm]{.3mm}{3mm}}
\newcommand{\thknorm}[1]{\thkl #1 \thkl\,}
\newcommand{\trinorm}[1]{|\!|\!| #1 |\!|\!|\,}
\newcommand{\bang}[1]{\langle #1 \rangle}
\def\angb<#1>{\langle #1 \rangle}
\newcommand{\vstrut}[1]{\rule{0mm}{#1}}
\newcommand{\rec}[1]{\frac{1}{#1}}
\newcommand{\opname}[1]{\mbox{\rm #1}\,}
\newcommand{\supp}{\opname{supp}}
\newcommand{\dist}{\opname{dist}}
\newcommand{\myfrac}[2]{{\displaystyle \frac{#1}{#2} }}
\newcommand{\myint}[2]{{\displaystyle \int_{#1}^{#2}}}
\newcommand{\mysum}[2]{{\displaystyle \sum_{#1}^{#2}}}
\newcommand {\dint}{{\displaystyle \myint\!\!\myint}}
\newcommand{\q}{\quad}
\newcommand{\qq}{\qquad}
\newcommand{\hsp}[1]{\hspace{#1mm}}
\newcommand{\vsp}[1]{\vspace{#1mm}}
\newcommand{\ity}{\infty}
\newcommand{\prt}{\partial}
\newcommand{\sms}{\setminus}
\newcommand{\ems}{\emptyset}
\newcommand{\ti}{\times}
\newcommand{\pr}{^\prime}
\newcommand{\ppr}{^{\prime\prime}}
\newcommand{\tl}{\tilde}
\newcommand{\sbs}{\subset}
\newcommand{\sbeq}{\subseteq}
\newcommand{\nind}{\noindent}
\newcommand{\ind}{\indent}
\newcommand{\ovl}{\overline}
\newcommand{\unl}{\underline}
\newcommand{\nin}{\not\in}
\newcommand{\pfrac}[2]{\genfrac{(}{)}{}{}{#1}{#2}}

\def\ga{\alpha}     \def\gb{\beta}       \def\gg{\gamma}
\def\gc{\chi}       \def\gd{\delta}      \def\ge{\epsilon}
\def\gth{\theta}                         \def\vge{\varepsilon}
\def\gf{\phi}       \def\vgf{\varphi}    \def\gh{\eta}
\def\gi{\iota}      \def\gk{\kappa}      \def\gl{\lambda}
\def\gm{\mu}        \def\gn{\nu}         \def\gp{\pi}
\def\vgp{\varpi}    \def\gr{\rho}        \def\vgr{\varrho}
\def\gs{\sigma}     \def\vgs{\varsigma}  \def\gt{\tau}
\def\gu{\upsilon}   \def\gv{\vartheta}   \def\gw{\omega}
\def\gx{\xi}        \def\gy{\psi}        \def\gz{\zeta}
\def\Gg{\Gamma}     \def\Gd{\Delta}      \def\Gf{\Phi}
\def\Gth{\Theta}
\def\Gl{\Lambda}    \def\Gs{\Sigma}      \def\Gp{\Pi}
\def\Gw{\Omega}     \def\Gx{\Xi}         \def\Gy{\Psi}

\def\CS{{\mathcal S}}   \def\CM{{\mathcal M}}   \def\CN{{\mathcal N}}
\def\CR{{\mathcal R}}   \def\CO{{\mathcal O}}   \def\CP{{\mathcal P}}
\def\CA{{\mathcal A}}   \def\CB{{\mathcal B}}   \def\CC{{\mathcal C}}
\def\CD{{\mathcal D}}   \def\CE{{\mathcal E}}   \def\CF{{\mathcal F}}
\def\CG{{\mathcal G}}   \def\CH{{\mathcal H}}   \def\CI{{\mathcal I}}
\def\CJ{{\mathcal J}}   \def\CK{{\mathcal K}}   \def\CL{{\mathcal L}}
\def\CT{{\mathcal T}}   \def\CU{{\mathcal U}}   \def\CV{{\mathcal V}}
\def\CZ{{\mathcal Z}}   \def\CX{{\mathcal X}}   \def\CY{{\mathcal Y}}
\def\CW{{\mathcal W}} \def\CQ{{\mathcal Q}}
\def\BBA {\mathbb A}   \def\BBb {\mathbb B}    \def\BBC {\mathbb C}
\def\BBD {\mathbb D}   \def\BBE {\mathbb E}    \def\BBF {\mathbb F}
\def\BBG {\mathbb G}   \def\BBH {\mathbb H}    \def\BBI {\mathbb I}
\def\BBJ {\mathbb J}   \def\BBK {\mathbb K}    \def\BBL {\mathbb L}
\def\BBM {\mathbb M}   \def\BBN {\mathbb N}    \def\BBO {\mathbb O}
\def\BBP {\mathbb P}   \def\BBR {\mathbb R}    \def\BBS {\mathbb S}
\def\BBT {\mathbb T}   \def\BBU {\mathbb U}    \def\BBV {\mathbb V}
\def\BBW {\mathbb W}   \def\BBX {\mathbb X}    \def\BBY {\mathbb Y}
\def\BBZ {\mathbb Z}   \def\BBQ {\mathbb Q}

\def\GTA {\mathfrak A}   \def\GTB {\mathfrak B}    \def\GTC {\mathfrak C}
\def\GTD {\mathfrak D}   \def\GTE {\mathfrak E}    \def\GTF {\mathfrak F}
\def\GTG {\mathfrak G}   \def\GTH {\mathfrak H}    \def\GTI {\mathfrak I}
\def\GTJ {\mathfrak J}   \def\GTK {\mathfrak K}    \def\GTL {\mathfrak L}
\def\GTM {\mathfrak M}   \def\GTN {\mathfrak N}    \def\GTO {\mathfrak O}
\def\GTP {\mathfrak P}   \def\GTR {\mathfrak R}    \def\GTS {\mathfrak S}
\def\GTT {\mathfrak T}   \def\GTU {\mathfrak U}    \def\GTV {\mathfrak V}
\def\GTW {\mathfrak W}   \def\GTX {\mathfrak X}    \def\GTY {\mathfrak Y}
\def\GTZ {\mathfrak Z}   \def\GTQ {\mathfrak Q}

\font\Sym= msam10 
\def\SYM#1{\hbox{\Sym #1}}
\newcommand{\bdw}{\prt\Gw\xspace}
\maketitle\medskip

{\abstract We study properties of nonnegative functions satisfying (E)$\;-\Gd u+u^p-M|\nabla u|^q=0$ in a domain of $\BBR^N$ when $p>1$, $M>0$ and $1<q<p$. We concentrate our analysis on the solutions of (E) with an isolated singularity, or in an exterior domain, or in the whole space. The existence of such solutions and their behaviours depend strongly on the values of the exponents $p$ and $q$ and in particular according to the sign of $q-\frac{2p}{p+1}$, and when $q=\frac{2p}{p+1}$, also on the value of the  parameter $M$ which becomes a key element. The description of the different behaviours is made possible by a sharp analysis of the radial solutions of (E). 
}\medskip

\nind {\it 2010 Mathematics Subject Classification:} 35J60-35J62-35J75-34A34-34C37\smallskip

\nind{\it  Keywords:} Elliptic equations, isolated singularities, dynamical systems, linearization, super and sub solutions, Lyapounov functions, central manifold.

\tableofcontents
\date{}



\mysection{Introduction}
The aim of this article is to study existence and properties of nonnegative singular solutions  of the following equation
\begin{equation}\label{Z1}
\BA{lll}
\CL^M_{p,q}u:=-\Gd u+u^p-M|\nabla u|^{q}=0
\EA
\end{equation}
in a domain $\Gw$ of $\BBR^{N}$ or in $\BBR^N$ ($N\geq 2$), where $M$ is a real number and $p>q>1$. 

In the case $M<0$ many results dealing with isolated singularities have been obtained in \cite{NPT}. Therefore we will mainly concentrate on the case $M>0$  where the two nonlinear terms act in a opposite direction: one is an absorption and
 the other is a source. Furthermore they are not of the same type, one involves the function and the other its gradient. 
 
 \smallskip
 First we consider the case $q=\frac{2p}{p+1}$. Then $(\ref{Z1})$ becomes
\begin{equation}\label{Z7}
\BA{lll}
\CL^M_{p,\frac{2p}{p+1}} u:=-\Gd u+u^p-M|\nabla u|^{\frac{2p}{p+1}}=0,
\EA
\end{equation}
and this equation is invariant under the scaling transformation $T_\ell$, $\ell>0$,  defined by
 \begin{equation}\label{Z2}
T_\ell [u](x)=\ell^{\ga}u(\ell x).
\end{equation}
In that case there may exist self-similar solutions, necessarily under the form $u(x)=u(r,s)=r^{-\ga}\gw(s)$, where $(r,s)\in\BBR_+\ti S^{N-1}$ are the spherical coordinates in $\BBR^N$. The function $\gw$ is a solution of the following equation on $S^{N-1}$
 \begin{equation}\label{Z3}
-\Gd'\gw+\ell_{N,p}\gw+\gw^p-M\left(\ga^2\gw^2+|\nabla '\gw|^2\right)^{\frac{p}{p+1}}=0,
\end{equation}
where $\Gd'$ and $\nabla '$ denote respectively the Laplace-Beltrami operator and the tangential gradient on $S^{N-1}$, identified with the covariant gradient on $S^{N-1}$ for the metric  induced by the standard one in $\BBR^N$, and where
 \begin{equation}\label{Z4}
\ga=\frac {2}{p-1}\,\text{ and }\,\ell_{N,p}=\ga K\,
\end{equation}
with 
 \begin{equation}\label{Z4+}
K=N-2-\ga=\myfrac{(N-2)p-N}{N-2}. 
\end{equation}
The nonzero constant solutions of $(\ref{Z3})$ are the positive zeros of the function 
 \begin{equation}\label{Z5}
 \CP_{\!_M}(x)=x^{p-1}-M\ga^{\frac{2p}{p+1}}x^{\frac{p-1}{p+1}}+\ell_{N,p}.
 \end{equation}
 The following value of the parameter $M$, which exists only if $N\geq 3$ and $p\geq\frac{N}{N-2}$, plays an important role in the study of $(\ref{Z3})$:
  \begin{equation}\label{Z6}
m^*:=(p+1)\left(\myfrac{(N-2)p-N}{2p}\right)^{\frac p{p+1}}.
 \end{equation}
 The separable solutions obtained in the next theorem are at the core of the process of describing the behaviour of positive solutions of $(\ref{Z1})$ near an isolated singularity or in an exterior domain of $\BBR^N$. 
\bth{T1} Let $p>1$, then \smallskip

\nind 1- If $M\leq 0$ equation $(\ref{Z3})$ admits a positive solution if and only if $N=2$ or $N\geq 3$ and $p<\frac{N}{N-2}$. Furthermore this solution is constant, unique and denoted by $x_{_{\!M}}$. \smallskip

\nind 2-   If $M>0$ and $p\leq \frac{N}{N-2}$ if $N\geq 3$, or any $p>1$ if $N=2$, equation $(\ref{Z3})$ admits a unique positive solution. It is constant and denoted by $x_{_{\!M}}$.\smallskip

\nind 3-  If $N\geq 3$, $p>\frac{N}{N-2}$ and $M=m^*$ there exists one positive solution to $(\ref{Z3})$.  It is constant and denoted by $x_{{m^*}}$.\smallskip

\nind 4-  If $N\geq 3$, $p>\frac{N}{N-2}$ and $0<M<m^*$ there exists no positive solution to $(\ref{Z3})$.
\smallskip

\nind 5-  If $N\geq 3$, $p>\frac{N}{N-2}$ and $M>m^*$ there exist two constant positive solutions $x_{_{1,M}}<x_{_{2,M}}$ to $(\ref{Z3})$ and any positive solution $\gw$ satisfies
  \begin{equation}\label{Z7-}\displaystyle
0<\min_{S^{N-1}}\gw\leq x_{_{1,M}}\leq \max_{S^{N-1}}\gw \leq x_{_{2,M}}.
 \end{equation}
Furthermore, if 
  \begin{equation}\label{Z7-1}\displaystyle
m^*<M<\tilde m:=\myfrac{(p+1)^2}{2}\left(\myfrac{(N-2)p^2-(N+2)}{4p^2}\right)^{\frac{p}{p+1}},
 \end{equation}
then $x_{_{1,M}}$ and $x_{_{2,M}}$ are the only positive solutions, and
  \begin{equation}\label{Z7-2}\frac{\tilde m}{m^*}>\left(\myfrac{p+1}{2p}\right)^{\frac{p}{p+1}}\myfrac{p+1}{2} >\myfrac{N-1}{N-2}
  \left(\myfrac{N-1}{N}\right)^{\frac{N}{2(N-1)}}>1.
  \end{equation}.
  \es

Not all the singular positive solutions of $(\ref{Z7})$
 are self-similar since there exist solutions with a {\it weak singularity}, which means
 \begin{equation}\label{Z8}
\BA{lll}\displaystyle
(i)\qquad\qquad&\displaystyle\lim_{x\to 0}|x|^{N-2}u(x)=k\qquad\qquad&\text{if }N\geq 3,\qquad\qquad\\[2mm]\displaystyle
(ii)\qquad\qquad &\displaystyle\lim_{x\to 0}\left|\ln|x|^{\phantom{m^m}}\!\!\!\!\!\!\!\!\right|^{-1}u(x)=k\qquad\qquad&\text{if }N=2.\qquad\qquad
\EA
\end{equation}

Thanks to the existence of positive radial sigular solutions in  $\BBR^N\setminus\{0\}$ we are able to prove the existence of non-radial positive solution in 
a punctured bounded domain with prescribed boundary value. This is a very general tool which is developed in Section 4 for obtaining singular solutions, and as an example we prove the following result.

 \bth{T2} Let $\Gw$ be a bounded smooth domain of $\BBR^N$ ($N\geq 3$) containing $0$ and $\gf\in W^{1,\infty}(\prt\Gw)$. If $1<p<\frac{N}{N-2}$ then for any real $M>0$ and $k>0$ there exists a minimal positive solution $u_k$ of $(\ref{Z7})$ in $\Gw\setminus\{0\}$ satisfying $(\ref{Z8})$ and such that $u=\phi$ on $\prt\Gw$.  Furthermore, $k\mapsto u_k$ is increasing and $u_k\uparrow u_\infty$ where $u_\infty$ is the minimal solution of $(\ref{Z7})$ in $\Gw\setminus\{0\}$ satisfying $(\ref{Z8})$, such that $u=\phi$ on $\prt\Gw$  and satisfying 
\begin{equation}\label{Z9}
\BA{lll}\displaystyle
\lim_{x\to 0}|x|^{\frac 2{p-1}}u_\infty(x)=x_{_{M}}.
\EA
\end{equation}
\es

If $p\geq \frac{N}{N-2}$  it is proved in \cite{BVGHV3}  that there exists no positive solution of $(\ref{Z7})$ with weak singularity at $0$ and that any positive solution  in $\Gw\setminus\{0\}$ can be extended as a weak solution in whole $\Gw$. However weak solutions may be unbounded. The different kinds of singular solutions play a key role for describing the behaviour near $0$ of any positive solution of $(\ref{Z7})$ in $\Gw\setminus\{0\}$. If $\Gw$ is replaced by $\BBR^N$ there holds:

 \bth{T3} Let $N\geq3$ and $1<p<\frac{N}{N-2}$. Then for any real $M>0$ and $k>0$ there exists a unique positive solution $u_k$ of $(\ref{Z7})$ in $\BBR^N\setminus\{0\}$ satisfying $(\ref{Z8})$ and 
 \begin{equation}\label{Z10}
\BA{lll}\displaystyle
\lim_{|x|\to \infty}|x|^{\frac 2{p-1}}u_k(x)=x_{_{M}}.
\EA
\end{equation}
Furthermore $u_k$ is radial and $u_k\uparrow u_{x_{_{M}}}$ as $k\to\infty$, where $u_{x_{_{M}}}(x)=x_{_{M}}|x|^{-\ga}$.\es
 
 When $p\geq \frac{N}{N-2}$ new phenomena appear.
  \bth{T4} Let $N\geq 3$, $p=\frac{N}{N-2}$ and $M>0$. Then the function  
  $$u_{x_{_{M}}}(x)=x_{_{M}}|x|^{2-N}=\left((N-2)M^{\frac{N-1}{N}}\right)^{N-2}|x|^{2-N}$$ 
  is the unique radial positive solution of $(\ref{Z7})$ in $\BBR^N\setminus\{0\}$ satisfying $\displaystyle\lim _{|x|\to 0}|x|^{N-2}u(x)=x_{_{M}}$.
  Moreover
  there exists a positive solution $u_{_S}$ of $(\ref{Z7})$ in $\BBR^N\setminus\{0\}$ satisfying 
\begin{equation}\label{Z11}
\BA{lll}\displaystyle
(i)\qquad&\displaystyle\lim_{|x|\to 0}|x|^{N-2}|\ln |x||^{N-1} u_{_S}(x)=\myfrac{1}{(N-2)\left((N-1)M\right)^{N-1}}\qquad\quad\\[4mm]
(ii)\qquad&\displaystyle\lim_{|x|\to \infty}|x|^{N-2} u_{_S}(x)=\left((N-2)M^{\frac{N-1}{N}}\right)^{N-2}.
\EA
\end{equation}
Furthermore $u_{_s}$ is the unique positive solution (not only radial) satisfying $(\ref{Z11})$. 
\es

The proof of existence is based upon a dynamical system formulation of the equation, see $(\ref{W4})$. Such a formulation, as well as similar ones, will be much used in the sequel. 
 \bth{T5} Let $N\geq 3$ and $p>\frac{N}{N-2}$.\\
\nind 1-  If $M>m^*$, besides the two self-similar solutions $u_{{x_{_{j,M}}}}$ ($j=1,2$), there exists a radial positive solution $u_{s}$ of $(\ref{Z7})$ in $\BBR^N\setminus\{0\}$, unique among the radial ones up to  the scaling transformation $T_\ell$,  satisfying 
 \begin{equation}\label{Z12}
\BA{lll}\displaystyle \lim_{|x|\to 0}|x|^{\ga}u_{s}(x)=x_{_{1,M}}\quad\text{and }\;\lim_{|x|\to \infty}|x|^{\ga}u_{s}(x)=x_{_{2,M}}.
\EA
\end{equation}
For any $k>0$ there exists also a radial positive solution $u$ of $(\ref{Z7})$ in $\BBR^N\setminus\{0\}$ satisfying
  \begin{equation}\label{Z13}
\BA{lll}\displaystyle \lim_{|x|\to 0}|x|^{\ga}u(x)=x_{_{1,M}}\quad\text{and }\;\lim_{|x|\to \infty}|x|^{N-2}u(x)=k>0.
\EA
\end{equation}
It is unique among the radial positive solutions satisfying $(\ref{Z13})$. Furthermore $T_\ell[u_{c}]=u_{c{\ell^{\ga+2-N}}}$. \\
\nind 2-  If $M=m^*$, the self-similar solution $u_{{m^*}}(x)=x_{m^*}|x|^{-\ga}$ is the unique among the radial positive solutions of $(\ref{Z7})$ in $\BBR^N\setminus\{0\}$  
satisfying 
\begin{equation}\label{Z14}
\BA{lll}\displaystyle \lim_{|x|\to 0}|x|^{\ga}u(x)=x_{{m^*}}\quad\text{and }\;\lim_{|x|\to \infty}|x|^{\ga}u(x)=x_{_{m^*}}.
\EA
\end{equation}
Furthermore,  for any $k>0$ there exists also a radial positive solution $u_{k}$ of $(\ref{Z7})$ in $\BBR^N\setminus\{0\}$ satisfying 
$(\ref{Z13})$ with $x_{_{1,M}}$ replaced by $x_{{m^*}}$. It is unique among the positive radial solutions satisfying $(\ref{Z13})$ and it satisfies the same scaling invariance as (i).
\es

The previous results allow to describe the behaviour at infinity of radial positive solutions of $(\ref{Z7})$ in the complement of a ball. The next result will be partially extended to non-radial solutions in Section 5.
\bprop{hardc} Let $N\geq 1$, $p>1$, $M>0$ and $u$ be a positive radial solution of $(\ref{Z7})$ in $\BBR^N\setminus B_R$ for some $R>0$. \smallskip

\nind 1- If $N=2$, or $N\geq 3$ and $1<p<\frac{N}{N-2}$, then 
$\displaystyle \lim_{|x|\to\infty}|x|^{\ga}u(x)=x_{_M}$.\smallskip

\nind 2- If $N\geq 3$ and $p=\frac{N}{N-2}$, then 
$\displaystyle \lim_{|x|\to\infty}|x|^{N-2}(\ln |x|)^{\frac{N-2}{2}}u(x)=\left(\tfrac{N-2}{\sqrt 2}\right)^{N-2}$ or 
 $\lim_{|x|\to\infty}|x|^{N-2}u(x)=x_{_M}$.\smallskip

\nind 3- If $N\geq 3$ and $p>\frac{N}{N-2}$, \smallskip

\nind  3-a- if $0<M<m^*$, then $\displaystyle \lim_{|x|\to\infty}|x|^{N-2}u(x)=k$ for some $k>0$.\smallskip

\nind  3-b- if $M=m^*$, then either  $\displaystyle \lim_{|x|\to\infty}|x|^{\ga}u(x)=x_{m^*}$, or $\displaystyle \lim_{r\to\infty}r^{N-2}u(r)=k$ for some $k>0$. \smallskip

\nind  2-c- if $M>m^*$, then either  $\displaystyle \lim_{|x|\to\infty}|x|^{\ga}u(x)=x_{_{1,M}}$, or $\displaystyle \lim_{|x|\to\infty}|x|^{\ga}u(x)=x_{_{2,M}}$ or $\displaystyle \lim_{|x|\to\infty}|x|^{N-2}u(x)=k$ for some $k>0$.

\es

 Next, we consider equation  $(\ref{Z1})$ when $q\neq\frac{2p}{p+1}$. In that case, the asymptotics of the solutions are governed either by the Emden-Fowler operator
\begin{equation}\label{Z15}
u\mapsto \CL_pu:=-\Gd u+u^p,
\end{equation}
or by the Riccati operator
\begin{equation}\label{Z16}
u\mapsto \CR^M_qu:=-\Gd u-M|\nabla u|^q,
\end{equation}
or by the eikonal operator
\begin{equation}\label{Z16-1}
u\mapsto \CE^M_{p,q}u:=u^p-M|\nabla u|^q.
\end{equation}
When $1<q<\frac{2p}{p+1}$ the governing equation is the Emden-Fowler equation $\CL_pu=0$ near a singularity and the Riccati equation $\CR^M_qu=0$ at infinity. When $\frac{2p}{p+1}<q<p$, the situation is reversed. The following exponents play a crucial role
\begin{equation}\label{Z17}
\ga=\myfrac{2}{p-1}\text{ },\; \gb=\myfrac{2-q}{q-1}\text{ }\text{ and }\;\gg=\myfrac{q}{p-q}\,\text{ if }\;q\neq p,
\end{equation}
and
\bel{Usigma}\BA {lll}\displaystyle  
\gs=(p+1)q-2p.
\EA\ee
We also define 
\begin{equation}\label{Z18}
\gk=\myfrac{(N-1)q-N}{q-1}\,\text{ if }\;q>1,
\end{equation}
and
\begin{equation}\label{Z18'}
\gth=\myfrac{(N-1)q-(N-2)p}{p-q}=\gg+2-N\,\text{ if }\;q\neq p.
\end{equation}

\bth{T6} Let $N\geq 1$, $M>0$ and $\frac{2p}{p+1}<q<p$. If there exists a radial positive solution $u$ of $(\ref{Z1})$ in $B_R\setminus\{0\}$ which is unbounded near $0$, then \smallskip

\nind 1- either 
\begin{equation}\label{Z19}\displaystyle
\lim_{x\to 0}|x|^\gg u(x)=X_{_M}\quad\text{where }\;X_{_M}=(M\gg^q)^{\frac 1{p-q}},
\end{equation}

\nind 2- or $(\ref{Z19})$ does not hold. In that case $q\leq 2$, $N\geq 2$ and the following situation occurs:\smallskip

\nind 2-a- if $\frac {N}{N-1}<q<2$, then
\begin{equation}\label{Z20}\displaystyle
\lim_{x\to 0}|x|^\gb u(x)=\xi_{_M}\quad\text{where }\;\xi_{_M}=\myfrac{1}{\gb}\left(\myfrac{\gk}{M}\right)(M\gg^q)^{\frac 1{q-1}},
\end{equation}

\nind 2-b- if $q=2$, then
\begin{equation}\label{Z21}\displaystyle
\lim_{x\to 0}|\ln |x||^{-1} u(x)=\myfrac{N-2}{M}\;\text{ if }\;N\geq 3,\;\text{or }\;\lim_{x\to 0}\left(\ln|\ln |x|| \right)u(x)=\myfrac{1}{M}\;\text{ if }\;N=2,
\end{equation}

\nind 2-c- if $q<\frac N{N-1}$, then there exists $k>0$ such that
\begin{equation}\label{Z22}\BA{lll}\displaystyle
\lim_{x\to 0}|x|^{N-2} u(x)=k\,\text{ if }\,N\geq 3\,
\text{ and }\,
\lim_{x\to 0}|\ln |x||^{-1}u(x)=k\quad\,\text{ if }\,N=2,
\EA\end{equation}

\nind 2-d- if $q=\frac N{N-1}$, then 
\begin{equation}\label{Z23}\BA{lll}\displaystyle
(i)\quad\quad&\displaystyle\lim_{x\to 0}||x|\ln |x||^{N-1} u(x)=\myfrac 1{N-1}\left(\myfrac{N-1}{M}\right)^{N-1}\quad&\text{if }N\geq 3\\[2mm]
(ii)\quad\quad&\displaystyle\lim_{x\to 0}|\ln |x||^{-1} u(x)=k>0&\text{if }N= 2
\EA\end{equation}
\es

In the case $1<q<\frac{2p}{p+1}$ the description of isolated singularities is simpler and it is similar to the one of the positive solutions of $(\ref{Z15})$.
\bth{T7} Let $M>0$, $1<p<\frac{N}{N-2}$ if $N\geq 3$, any $p>1$ if $N=1,2$, and $1<q<\frac{2p}{p+1}$. Assume that there exists a radial positive solution $u$ of $(\ref{Z1})$ in $B_R\setminus\{0\}$ 
which is unbounded near $0$. Then the following alternative holds:\smallskip

\nind 1- either 
\begin{equation}\label{Z24}\displaystyle
\lim_{x\to 0}|x|^\ga u(x)=x_{_0}:=(\ga |K|)^{\frac 1{p-1}},
\end{equation}
\smallskip

\nind 2- or $N\geq 2$ and
\begin{equation}\label{Z25}\BA{lll}\displaystyle
(i)\quad\quad&\displaystyle\lim_{x\to 0}|x|^{N-2}u(x)=k>0\quad&\text{if }N\geq 3\qquad\qquad\quad\quad\quad\quad\\[2mm]
(ii)\quad\quad&\displaystyle\lim_{x\to 0}|\ln |x||^{-1} u(x)=k>0\quad&\text{if }N=2.\quad\quad
\EA\end{equation}
\es

It is noticeable that all the behaviours described in the previous  two theorems occur. The behaviour at infinity of positive solutions of  
$(\ref{Z1})$ in $B^c_R$ inherits this complexity due to the value of $q$ with respect to $\frac{2p}{p+1}$, and the situation is less
 intricated in the case $\frac{2p}{p+1}<q<p$ than in the case $1<q<\frac{2p}{p+1}$.

\bth{T8} Let $N\geq 1$, $M>0$ and $\frac{2p}{p+1}<q<p$. Assume that there exists a radial positive solution $u$ of $(\ref{Z1})$ in $B^c_R$. Then
 
\nind 1- If $1<p<\frac{N}{N-2}$ (any $p>1$ if $N=1,2$), there holds
\begin{equation}\label{Z26}\displaystyle
\lim_{|x|\to \infty}|x|^{\ga}u(x)=x_{_0}.
\end{equation}
\nind 2- If $N\geq 3$ and $p>\frac{N}{N-2}$, there holds
\begin{equation}\label{Z27}\displaystyle
\lim_{|x|\to \infty}|x|^{N-2}u(x)=k>0.
\end{equation}
\nind 3- If $N\geq 3$ and  $p=\frac{N}{N-2}$, there holds
\begin{equation}\label{Z28}\displaystyle
\lim_{|x|\to \infty}(\ln |x|)^{\frac{N-2}{2}}|x|^{N-2}u(x)=\left(\myfrac{N-2}{\sqrt 2}\right)^{N-2}.
\end{equation}
\es
\bth{T9} Let $N\geq 2$, $M>0$ and $1<q<\frac{2p}{p+1}$. If  $u$ is a radial positive solution of $(\ref{Z1})$ in $B_R^c$, there holds.
 
\nind 1- If $q>\frac{N}{N-1}$, one of the three following situations occurs:\smallskip

\nind 1-a- either
\begin{equation}\label{Z29}\displaystyle
\lim_{|x|\to \infty}|x|^{\gg}u(x)=X_{_M},
\end{equation}
\nind 1-b- or
\begin{equation}\label{Z30}\displaystyle
\lim_{|x|\to \infty}|x|^{\gb}u(x)=\xi_{_M},
\end{equation}
\nind 1-c- or
\begin{equation}\label{Z31}\displaystyle
\lim_{|x|\to \infty}|x|^{N-2}u(x)=k>0.
\end{equation}

\nind 2- If $N=1$ or $1<q\leq\frac{N}{N-1}$, then only $(\ref{Z29})$ can occur.
\es

The existence of local or global singular solutions or asymptotic solutions with behaviour like $|x|^{-\gg}$ (eikonal type) or like $|x|^{-\gb}$ (Riccati type) near $0$ or $\infty$ will be proved in Section 3.6. 
For example we prove the following result by the method of sub and super solutions.
\bth{T10} Let $N\geq 1$, $p>1$ and $M>0$.  \smallskip

\nind 1- If $\frac{2p}{p+1}<q<p$, then there exists a unique global positive solution $u$ of $(\ref{Z1})$ such that $\displaystyle\lim_{x\to 0}|x|^\gg u(x)=X_{_M}$, and its behaviour at infinity is given by \rth{T8}. Moreover this solution is radial, and it is explicit if $N\geq 3$ and  $q=\frac{(N-2)p}{N-1}$. Furthermore, for any bounded smooth domain $\Gw$ containing $0$ there exists a positive solution of $(\ref{Z1})$ in  $\Gw\setminus\{0\}$ vanishing on $\prt\Gw$. \smallskip

\nind 2- If $\max\left\{1,\frac{(N-2)p}{N-1}\right\}<q<\frac{2p}{p+1}$, then for any $R>0$ there exists a positive solution $u$ in  $B^c_R$ satisfying $(\ref{Z29})$.
\es
Introducing a new powerful autonomous system of order 3, we can construct local solutions behaving like  $|x|^{-\gb}$ near $0$.
\bth{T11} Let $N\geq 2$, $p>1$ and $M>0$.  \smallskip

\nind 1- If $\max\left\{\frac{2p}{p+1},\frac{N}{N-1}\right\}<q<\min\{2,p\}$. Then there exists at least one radial positive solution $u$ of $(\ref{Z1})$ in a neighborhood of $0$ such that $\displaystyle\lim_{x\to 0}|x|^\gb u(x)=\xi_{_M}$. \smallskip

\nind 2- If $\frac{N}{N-1}<q<\frac{2p}{p+1}$ there exists a unique positive radial solution defined in a neighborhood of infinity satisfying such that $\displaystyle\lim_{|x|\to \infty}|x|^\gb u(x)=\xi_{_M}$. There exists no radial positive solution in $\BBR^N\setminus\{0\}$ with such a behaviour at infinity. 
\es

By a delicate method of super and sub solutions, we also prove the existence of radial positive singular solutions $u$ of $(\ref{Z1})$ in $\BBR^N\setminus\{0\}$ satisfying $(\ref{Z30})$ under more restrictive assumptions on the exponents $p$ and $q$. \smallskip

When $p<\frac{N}{N-2}$ we show the existence of the solutions of 
$(\ref{Z1})$ in $\BBR^N\setminus\{0\}$, or in a neighborhood of $0$, or at infinity having the behaviour described in \rth{T7} and \rth{T8}. Such solutions are associated to the Emden-Fowler operator.
\bth{T12} Let $M>0$, $N\geq 3$ and $1<p<\frac{N}{N-2}$, or $N=1,2$ and $p>1$.\\
\nind 1- If $1<q<\frac{2p}{p+1}$ there exists a unique positive solution of $(\ref{Z1})$ in $\BBR^N\setminus\{0\}$ satisfying 
\begin{equation}\BA{lll}\label{Z32}\displaystyle
(i) &\qquad \qquad \displaystyle\lim_{x\to 0}|x|^\ga u(x)=x_0\qquad \qquad \qquad\qquad \qquad \qquad\qquad\qquad  \\[2mm]
(ii) &\qquad \qquad \displaystyle \lim_{|x|\to \infty}|x|^\gg u(x)=X_{_M}.
\EA\end{equation}
Furthermore this solution is radial and $|x|^\ga u(x)\geq x_0$ in $\BBR^N\setminus\{0\}$. If $\Gw$ is a bounded domain containing $0$ there exists a positive solution $u$ of $(\ref{Z1})$ in $\Gw\setminus\{0\}$ satisfying $(\ref{Z32})$-(i) and vanishing on $\prt\Gw$.\\
\nind 2- If $\frac{2p}{p+1}<q<p$ there exists a positive radial solution of $(\ref{Z1})$ in $\BBR^N\setminus\{0\}$ satisfying 
\begin{equation}\BA{lll}\label{Z33}\displaystyle
(i) &\qquad \qquad \displaystyle\lim_{x\to 0}|x|^\gg u(x)=X_{_M}\qquad \qquad \qquad\qquad \qquad \qquad\qquad\qquad  \\[1mm]
(ii) &\qquad \qquad \displaystyle \lim_{|x|\to \infty}|x|^\ga u(x)=x_0.
\EA\end{equation}
Moreover this solution is unique among all the positive solutions satisfying $(\ref{Z33})$.
\es

We also give conditions on $p$ and $q$ for the existence of a positive radial solution of $(\ref{Z1})$ in $\BBR^N\setminus\{0\}$, necessarily singular at $0$, with a behaviour at $0$ given by $(\ref{Z19})$, $(\ref{Z24})$ or $(\ref{Z20})$, and an asymptotic behaviour at infinity given by $(\ref{Z31})$.\medskip

The last section of the article is devoted to non radial results. We first give a general existence statement which allows to construct positive singular solutions of $(\ref{Z1})$ in a punctured bounded domain with prescribed boundary value, provided there exists a radial singular solution  in $\BBR^N\setminus \{0\}$. This singular solution 
has been obtained by the phase plane analysis of Section 2 in the case $q=\frac{2p}{p+1}$, and by the radial analysis of section 3 in the other cases.

\bth{Exist1} Let $\Gw\subset B_R\subset\BBR^N$ be a bounded smooth domain containing $0$, $M$ a real number, $p>1$ and $1\leq q\leq 2$. If there exists a radial positive and decreasing function $v$ defined in $B_R\setminus \{0\}$ 
and satisfying $\CL_{p,q}^Mv=0$ in $\Gw \setminus \{0\}$ and 
$$\displaystyle \lim_{x\to 0}v(x)=\infty,
$$
then for any nonnegative function $\gf\in W^{1,\infty}(\Gw)$, there exists a solution $u$ of $\CL_{p,q}^Mu=0$ in $\Gw\setminus \{0\}$ satisfying 
$u=\gf$ on $\prt\Gw$ and 
$$\displaystyle \lim_{x\to 0}u(x)=\infty.
$$
Furthermore there holds
\begin{equation}\label{X3}\displaystyle 
(v(x)-\max_{z\in\prt\Gw}\phi(z))_+\leq u(x)\leq v(x)+\max_{z\in\prt\Gw}(\phi(z)-v(z))_+\quad\text{for all }x\in \Gw\setminus \{0\}.
\end{equation}
\es

A second key result deals with the uniqueness of positive solutions in $\BBR^N\setminus \{0\}$ or in a punctured bounded domain $\Gw\setminus\{0\}$ starshaped with respect to $0$. Using a general scaling method we prove the following
\bth{uni1} Assume $N\geq 2$, $p>1$, $1<q<2$ and $M>0$. Let $a$ such that 
\bel{UNX}\BA{lll}
(i)\phantom{-----}\displaystyle 0\leq a<\gb\quad&\text{if }q\leq \tfrac{2p}{p+1}\phantom{------------}\\[2mm]
(ii)\phantom{-----}\displaystyle \gb<a\quad&\text{if }q> \tfrac{2p}{p+1}.\phantom{------------}
\EA\ee
There exists at most one positive solution of $(\ref{Z1})$ in $\BBR^N\setminus\{0\}$ satisfying
\bel{UNX1}
\displaystyle \lim_{x\to 0}|x|^a|\ln |x||^{\tilde a}u(x)=\Gl
\ee
where $\Gl$ is some positive constant and $\tilde a$ is a real number. If $\Gw$ is a bounded domain containing $0$ and starshaped with respect to $0$ and $\gf\in C(\prt\Gw)$ is nonnegative, there is at most one positive solution $u$ of $(\ref{Z1})$ in $\Gw\setminus\{0\}$ satisfying $(\ref{UNX1})$ with value $\gf$ on $\prt\Gw$.
\es

This result admits various extensions valid when the exponent $a$ above is equal to $\gb$. \\
With the help of these results we characterize {\it all the local positive solutions of $(\ref{Z1})$, not necessarily radial}, either 
 near $0$ or near $\infty$. An important tool is the intensive use of the tangency property of graphs of global solutions which has been introduced in \cite{KiVe} for the studying of isolated singularities of $p$-harmonic functions. \medskip

\noindent{\bf Acknowledgements}. This article has been prepared with the support of the  FONDECYT grants 1210241 and 1190102 for the three authors. 
\mysection{The case $q=\frac{2p}{p+1}$}

\subsection{The equation on the sphere}
The existence of particular solutions of $(\ref{Z3})$, and eventually their uniqueness,  plays a key role in the description of the behaviour of all the solutions. Due to the invariance of the equation under the transformations $T_\ell$ these natural particular solutions are the ones which are self-similar, i.e. invariant by these transformations. In spherical coordinates $(r,s)\in \BBR_+\ti S^{N-1}$, they endow the form $(r,s)\mapsto u(r,s)=r^{-\ga}\gw(s)$, and $\gw$ is a solution of $(\ref{Z3})$. Since we are dealing with nonnegative solutions,   by the strong maximum principle they are either positive or identically zero. This fact does not depend on the sign of $M$.
\subsubsection{ Proof of \rth{T1}: constant positive solutions}
Assume $M\leq 0$ and $\gw$ is a nonnegative solution of $(\ref{Z3})$. Multiplying the equation by $\gw$ and integrating over $S^{N-1}$ yields
$$\myint{S^{N-1}}{}\left(|\nabla' \gw|^2+\ell_{N,p}\gw^2+\gw^{p+1}-M\left(\ga^2\gw^2+|\nabla' \gw|^2\right)^{\frac{p}{p+1}}\gw\right) dS=0.
$$
Since $\ell_{N,p}\geq 0$ if and only if $p\geq \frac{N}{N-2}$, we infer the non-existence statement 1.  \smallskip

For any $M$, constant positive solutions are the positive roots of $ \CP_{\!_M}(x)=0$. If we set $z=x^{\frac{p-1}{p+1}}$, $ \CP_{\!_M}(x)=0$ is equivalent to
$\tilde  \CP_{\!_M}(z)=0$ where
$$\tilde  \CP_{\!_M}(z)=z^{p+1}-M\ga^{\frac{2p}{p+1}}z+\ell_{N,p}.
$$
Since $\tilde  \CP'_{\!_M}(z)=(p+1)z^{p}-M\ga^{\frac{2p}{p+1}}$ the minimum of $\tilde  \CP_{\!_M}$ on $\BBR_+$ is achieved at $z=0$ if $M\leq 0$, or at $z_0=\left(\frac{M}{p+1}\right)^{\frac{1}{p}}\ga^{\frac{2}{p+1}}$ if $M> 0$. In the first case the function $\tilde  \CP_{\!_M}$ is increasing on $(0,\infty)$. It vanishes therein if and only if 
$\tilde  \CP_{\!_M}(0)=\ell_{N,p}<0$, or equivalently $1<p<\frac{N}{N-2}$. In the second case, $\tilde  \CP_{\!_M}$ is decreasing on $(0,z_0)$ and increasing on $(z_0,\infty)$. Its minimal value is
\begin{equation}\label{Y4}\tilde  \CP_{\!_M}( z_0)=\ell_{N,p}-p\left(\myfrac{M}{p+1}\right)^{\frac{p+1}{p}}\!\!\!\ga^2=\myfrac{4p}{(p-1)^2}\left(\left(\myfrac{m^*}{p+1}\right)^{\frac{p+1}{p}}\!-\left(\myfrac{M}{p+1}\right)^{\frac{p+1}{p}}\right).
\end{equation}
If $p\leq \frac{N}{N-2}$, then $\tilde  \CP_{\!_M}(z_0)<\tilde  \CP_{\!_M}(0)=\ell_{N,p}\leq 0$, hence $\tilde  \CP_{\!_M}$ admits a unique positive zero and so does $ \CP_{\!_M}$. This gives the existence of $x_{_M}$ in case 2. \\
\nind If $p> \frac{N}{N-2}$, then $\tilde  \CP_{\!_M}(0)=\ell_{N,p}> 0$. We obtain the existence  of constant solutions in 3, 4 and 5 according $M>m^*$, $M=m^*$ and $0<M<m^*$.\qeda\medskip

\subsubsection{ Proof of \rth{T1}: positive solutions} Let $\gw$ be a nonnegative solution of ($\ref{Z3}$). By regularity it is $C^2$ and either positive or identically $0$. If it is not constant, we denote by $\overline\gw$ and $\underline\gw$ respectively the maximum and the minimum of $\gw$ on $S^{N-1}$. There holds $ \CP_{\!_M}(\overline\gw)\leq 0$ and $ \CP_{\!_M}(\underline\gw)\geq 0$,
and if we set $\overline\gw^{\frac{p-1}{p+1}}=\overline z$ and $\underline\gw^{\frac{p-1}{p+1}}=\underline z$, we have that $\tilde  \CP_{\!_M}(\overline z)\leq 0$ and $\tilde  \CP_{\!_M}(\underline z)\geq 0$.

\nind 1- First we consider the case where $M<0$  and $1<p<\frac{N}{N-2}$. Since $\tilde  \CP_{\!_M}$ is increasing on $\BBR_+$ we deduce $\overline\gw=\underline\gw=\gw=x_{_M}$.\\
\nind  2- Next we assume $M>0$ and $p\leq \frac{N}{N-2}$, then $\tilde  \CP_{\!_M}$ is increasing on $(z_0,\infty)$. Hence it is negative on  
$[0,x_{_M}^{\frac{p-1}{p+1}})$ and positive on $(x_{_M}^{\frac{p-1}{p+1}},\infty)$. This implies 
$\overline\gw\leq x_{_M}\leq \underline\gw$ and finally $\gw=x_{_M}$.\\ 
\nind 3- If $p>\frac{N}{N-2}$ and $M=m^*$, $\CP_{m^*}$ is positive on $\BBR^+\setminus\{x_{m^*}\}$. This implies $\overline\gw\leq x_{m^*}\leq \underline\gw$ and finally $\gw=x_{m^*}$. \\ 
\nind 4- If $p>\frac{N}{N-2}$ and $M<m^*$, $\CP_{m^*}$ is positive on $\BBR^+$, hence there exists no positive solution.\\
\nind 5- Finally, if $p>\frac{N}{N-2}$ and $M>m^*$, $ \CP_{\!_M}$ is positive on 
$(0,x_{_{1,M}})\cup (x_{_{2,M}},\infty)$ and negative on $(x_{_{1,M}},x_{_{2,M}})$. This implies $(\ref{Z7-})$. The proof of the second assertion is more involved. Set $z=|\nabla'\gw|^2$ and $y=\ga^2\gw^2+z^2$. Then
$$\Gd'\gw=\ga K\gw+\gw^p-My^{\frac{p}{p+1}}.
$$
By Weitzenb\"ock's formula
\bel {Y5}
\myfrac{1}{2}\Gd'z=|Hess (\gw)|^2+\langle\nabla'\Gd'\gw,\nabla'\gw\rangle+Ricc_g(\nabla'\gw,\nabla'\gw),
\ee
where $Hess (\gw)$ is the Hessian and $Ricc_g$ is the curvature 2-tensor on $S^{N-1}$. In that case we have that $Ricc_g=(N-2)g$.  By Schwarz inequality 
$$|Hess (\gw)|^2\geq\frac{1}{N-1}(\Gd'\gw)^2,$$
 therefore, replacing $\Gd'\gw$ by its value, we obtain the inequality
$$\BA {lll}-\myfrac{1}{2}\Gd'z +(N-2)z+\myfrac{1}{N-1}\left(\Gd'\gw\right)^2+(\ga K+p\gw^{p-1})z
-\myfrac{Mp}{p+1}y^{-\frac{1}{p+1}}\langle \nabla'y,\nabla'\gw\rangle\leq 0.
\EA$$
Since $\nabla'z=\nabla'y-2\ga^2\gw\nabla'\gw$ we infer
$$\BA {lll}-\myfrac{1}{2}\Gd'z+\left(N-2+\ga K+p\gw^{p-1}-\myfrac{2Mp}{p+1}y^{-\frac{1}{p+1}}\gw\right)z+\myfrac{1}{N-1}(\Gd'\gw)^2\\[3mm]
\phantom{-------------------}
-\myfrac{Mp}{p+1}y^{-\frac{1}{p+1}}\langle \nabla'z,\nabla'\gw\rangle\leq 0.
\EA$$
Let $s_0\in S^{N-1}$ where $z$ is maximal. Then $\nabla'z(s_0)=0$ and $\Gd'z(s_0)\leq 0$. Hence at $s=s_0$ there holds
$$\BA {lll}\left(\ga K+N-2+p\gw^{p-1}-\myfrac{2Mp\ga^2\gw}{(p+1)(\ga^2\gw^2+z)^{\frac{1}{p+1}}}\right)z\\[3mm]
\phantom{-------\left(\myfrac{2Mp\ga^2\gw}{(p+1)(\ga^2\gw^2+z)^{\frac{1}{p+1}}}\right)z}
+\myfrac{1}{N-1}\left(\ga K\gw+\gw^p-My^{\frac{p}{p+1}}\right)^2\leq 0.
\EA$$
Therefore
\bel {Y6}\BA {lll}
\left(\ga K+N-2+p\gw^{p-1}-\myfrac{2Mp\ga^\frac{2p}{p+1}\gw^\frac{p-1}{p+1}}{p+1}\right)z\\[2mm]
\phantom{----------}
+\myfrac{1}{N-1}\left(\ga K\gw+\gw^p-M(\ga^2\gw^2+z)^{\frac{p}{p+1}}\right)^2\leq 0.
\EA\ee
Set  
$$F(t)=pt^{p+1}-\myfrac{2Mp\ga^\frac{2p}{p+1}}{p+1}t+N-2+\ga K$$
and $t_0=\gw^{\frac{p-1}{p+1}}(s_0)$. If $\gw$ is non-constant, $z(s_0)>0$, hence 
$F(t_0)\leq 0.$ . If $t_i=x_{_{i,M}}^{\frac{p-1}{p+1}}$, for $i=1,2$, there holds 
$$t_i^{p+1}-M\ga^\frac{2p}{p+1}t_i+\ga K=0,
$$
hence
$$\BA{lll}F(t_i)=pt_i^{p+1}-\frac{2p}{p+1}\left(\ga K+t_i^{p+1}\right)+N-2+\ga K=\myfrac{p(p-1)}{p+1}t_i^{p+1}+N-2-\myfrac{2K}{p+1}.
\EA$$
Since 
$$F'(t)=p(p+1)t^{p}-\myfrac{2Mp\ga^\frac{2p}{p+1}}{p+1},
$$
$F$ is minimal for $t=t_*=\left(\frac{2M}{(p+1)^2}\right)^{\frac 1p}\ga^{\frac{2}{p+1}}$ and
\bel {Y7}\BA {lll}
F(t_*)=N-2+\ga K-p^2\left(\myfrac{2M}{(p+1)^2}\right)^{\frac{p+1}{p}}\ga^2\leq F(t_0)\leq 0.
\EA\ee
This implies  
\bel {Y8}\BA {lll}
\myfrac{2M}{(p+1)^2}\ga^{\frac{2p}{p+1}}\geq \left(\myfrac{N-2+\ga K}{p^2}\right)^{\frac{p}{p+1}},
\EA\ee
and equivalently $M\geq \tilde m$ where $\tilde m$ is defined in $(\ref{Z7-1})$. Therefore, if $M<\tilde m$ there cannot exist non-constant positive solution. In order to prove $(\ref{Z7-2})$
we first notice that if $p> \frac{N}{N-2}$, then 
$$\myfrac{(N-2)p^2-(N+2)}{(N-2)p-N}>p+1.
$$
Therefore 
$$\myfrac{\tilde m}{m^*}>\myfrac{p+1}{2}\left(\myfrac{p+1}{2p}\right)^{\frac p{p+1}}.
$$
By taking the logarithm it is easy to check that the function $p\mapsto \frac{p+1}{2}\left(\frac{p+1}{2p}\right)^{\frac p{p+1}}$ is increasing, hence
the right-hand side of the previous inequality is minorized by 
\bel {Y9}\gth_N:=\myfrac{N-1}{N-2}\left(\myfrac {N-1}{N}\right)^{\frac{N}{2(N-1)}},
\ee 
which is the desired estimate. Notice that $\gth_3\sim1.47$ and $\gth_N$ decreases to $1$ when $N\to\infty$.
\qeda\medskip

\nind\Remark The following monotonicity properties of the points  $P_{_M}$ are straightforward: in cases (i) and (ii) $x_{_M}$ is increasing with $M$. In case (iii)  
$M\mapsto x_{_{1,M}}$ is decreasing while $M\mapsto x_{_{2,M}}$ is increasing. Furthermore, if $M'>M>m^*$, 
\bel {Y10}
x_{_{1,M'}}<x_{_{1,M}}<x_{m^*}<x_{_{2,M}}<x_{_{2,M'}}.
\ee
The value of $x_{m^*}$ is explicit
\bel {Y11}
x_{m^*}=\left(\myfrac{2K}{p(p-1)}\right)^{\frac{1}{p-1}}=\left(\myfrac{\ga K}{p}\right)^\frac{1}{p-1}.
\ee

We end this section by proving a result dealing with bifurcation from constant solutions.
\bth{bif}When $M>m^*$ the solution $x_{_{j,M}}$, $j=1,2$, is never a bifurcation point in the sense that the linearized equation at this point is singular.
\es
\Proof If we look for 
solutions of $(\ref{Z3})$ under the form $\gw=x_{_{j,M}}+\ge\phi_k$ where $\gf_k$ is an eigenfunction of 
$-\Gd'$ in $H^1(S^{N-1})$ associated to the eigenvalue $\gl_k=k(N-2+k)$, we obtain that 
\bel {Y11+}
\gl_k+\ga K+px^{p-1}_{_{j,M}}-\myfrac{2p}{p+1}\ga^{\frac{2p}{p+1}}Mx^{\frac{p-1}{p+1}}_{_{j,M}}=0
\ee
We recall that $\CP_{_M}$ is defined in $(\ref{Z5})$. If $\CQ_{_M}(x)=x\CP_{_M}(x)$, then $(\ref{Y11+})$ is equivalent to 
\bel {Y12}
\gl_k+x_{_{j,M}}\CP'_{_M}(x_{_{j,M}})=0.
\ee
We know that $\CP'_{_M}(x_{_{2,M}})>0$, then for any $k\in\BBN^*$ identity $(\ref{Y11+})$ is impossible with $j=2$.  Concerning the case $j=1$, $(\ref{Y11+})$ combined with $\CP_{_M}(x_{_{1,M}})=0$ and the value of $\ga$ yields
\bel {Y13}\gl_k-2K+\myfrac{p(p-1)}{p+1}\ga^{\frac{2p}{p+1}}Mx^{\frac{p-1}{p+1}}_{_{1,M}}=0,
\ee
which never occurs when $k\geq 2$ because of the values of $\gl_k$ and $K$. When $k=1$, since $\CP_{_M}(x_{_{1,M}})=0$ there holds
$$M\ga^{\frac{2p}{p+1}}x^{\frac{p-1}{p+1}}_{_{1,M}}>\ga K
$$ 
thus $(\ref{Y13})$ yields
$$\gl_1-2K+\myfrac{p(p-1)}{p+1}\ga K<0.
$$
Because $\gl_1=N-1$ and $\ga(p-1)=2$ we obtain $N-1-2K+\frac{2p}{p+1}K< 0$, which is equivalent to 
$(p-1)(N-1)+2+2\ga< 0$, a contradiction.\qeda
\subsection{Radial solutions}
In this section we study in detail the nonnegative solutions of the ordinary differential equation
\bel{W1}
-u''-\myfrac{N-1}{r}u'+|u|^{p-1}u-M|u'|^{\frac{2p}{p+1}}=0,
\ee
when  $p>1$. Because of the scaling invariance $(\ref{Z2})$ the equation can be transformed into an autonomous equation by setting
\bel{W1+}u(r)=r^{-\ga}x(t),\quad t=\ln r.
\ee
Then $x(t)$ satisfies
\bel{W2}
x_{tt}+Lx_t-\ga Kx-|x|^{p-1}x+M\left|\ga x-x_t\right|^{\frac{2p}{p+1}}=0,
\ee
where we recall that $K=\frac{(N-2)p-N}{p-1}=N-2-\ga$ and where we set
\bel{W3}
L=\myfrac{(N-2)p-(N+2)}{p-1}=K-\ga.
\ee
If we set $u'(r)=-r^{-(\ga+1)}y(t)$, then $(\ref{W1})$ is equivalent to
\bel{W4}\BA {lll}
x_t=\ga x-y\\
y_t=-Ky-|x|^{p-1}x+M|y|^{\frac{2p}{p+1}}.
\EA\ee
Since we are interested in positive $u$ we restrict to solutions of  $(\ref{W4})$ in the half-space 
$\BBR^2_+=\{(x,y)\in\BBR^2: x>0\}=Q_1\cup Q_4$ where 
$$Q_1=\{(x,y)\in \BBR^2_+:y>0\}$$ is the first quadrant and $$Q_4=\{(x,y)\in \BBR^2_+:y<0\}$$
is the fourth quadrant. The regular solutions of $(\ref{W1})$ (with $u(0)=u_0>0$ and $u'(0)=0$) are increasing near 
$r=0$, so their trajectory $\CT_1:=\{(x(t),y(t))\}$ lies in $Q_4$ as $t\to-\infty$. The solutions defined in a neighborhood of 
$r=0$ and unbounded near $0$ are decreasing, so their trajectory lie in $Q_1$ as $t\to-\infty$. The solutions defined near $r=\infty$ are decreasing, so their trajectory remain in $Q_1$ as $t\to\infty$. \smallskip

\rth{T1} can be reformulated in the following way:\medskip

{\it
\nind 1-  If $M\leq 0$ and $1<p<\frac N{N-2}$, the only non-trivial equilibrium in $\overline\BBR^2_+$, is $P_{_M}=(x_{_M},\ga x_{_M})$. If 
$p\geq\frac N{N-2}$ there exists no non-trivial equilibrium in this region.\smallskip

\nind 2-  If $1<p\leq \frac N{N-2}$ and $M>0$, the only non-trivial equilibrium in $\overline\BBR^2_+$ is $P_{_M}=(x_{_M},\ga x_{_M})$.\smallskip

\nind 3-  If $p\geq \frac N{N-2}$ and $M> m^*$, there exist two non-trivial equilibria in  $\overline\BBR^2_+$, $P_{_{1,M}}=(x_{_{1,M}},\ga x_{_{1,M}})$ and 
$P_{_{2,M}}=(x_{_{2,M}},\ga x_{_{2,M}})$.\smallskip

\nind 4-  If $p\geq \frac N{N-2}$ and $M= m^*$, there exists one non-trivial equilibria in  $\overline\BBR^2_+$, $P_{m^*}=(x_{m*},\ga x_{m^*})$.\smallskip

\nind 5-  If $p\geq \frac N{N-2}$ and $0<M< m^*$ there exists no non-trivial equilibrium in $\overline\BBR^2_+$.}\medskip

We also recall the classical result concerning regular solutions, not only in the case $q=\frac{2p}{p+1}$. 
\bprop{reg} Let $N\geq 1$, $1<q<p$ and $M\geq 0$. Then for any $a>0$ there exists a unique radial maximal positive solution $u$ of $(\ref{Z1})$ satisfying $u(0)=a$, $u'(0)=0$. This solution denoted by $u_{[a]}$ is defined in $B_{R}$, where $R=R_a>0$, and it satisfies $\displaystyle\lim_{|x|\uparrow R}u_{[a]}=\infty$.
\es
\Proof In the case $M=0$ the result is a standard combination of Cauchy-Lipschitz theorem with the Keller-Osserman estimate. In the case $M>0$ the proof can be obtained in a slightly similar way using also \rprop{Oss}. See also 
\cite{BM}, \cite{Ve92} and \cite {BG} for many extensions concerning these regular (or large) solutions.\qeda

\subsubsection{Linearisation at $(0,0)$}
The linearization at $(0,0)$ is given by the system
\bel{W5}\BA {lll}
x_t=\ga x-y\\
y_t=-Ky
\EA\ee
with eigenvalues $\gl_1=-K$, $\gl_2=\ga$ and corresponding eigenvectors $\xi_1=(1,N-2)$ and $\xi_2=(1,0)$ if $K\neq -\ga$ or equivalently $N\neq 2$. If $N=2$, $\gl_1=\gl_2=\ga$, the only eigenspace is $span(\xi_2)$. In any cases $\gl_2-\gl_1=N-2$. There exists one trajectory  located in $Q_4$ of the linearized system
converging to $0$ when $t\to -\infty$. To this trajectory is associated a trajectory $\CT_r$ of $(\ref{W4})$ such that $(x(t), y(t))\approx c \left(e^{\ga t},-\tfrac{1}{N}e^{(\ga+1)t}\right)$ when $t\to-\infty$. These solutions are associated to the one parameter family of regular solutions mentioned above with $u(0)=u_0$ and $u'(0)=0$.  \
\smallskip

\nind (i) Assume first that $N\geq 3$. \\
\nind If $p<\frac{N}{N-2}$ then $K<0$ and $(0,0)$ is a source. Then all trajectories of $(\ref{W5})$ defined in a neighborhood of 
$(0,0)$ converge to this point when $t\to-\infty$. Besides the trajectory $\CT_r$, all the other trajectories converging to zero when $t\to-\infty$
start in $Q_1$ with initial slope $N-2$. They satisfy $x(t)\sim ae^{-Kt}$ for some $a>0$ by \rlemma {eigen}. This means that $r^{N-2}u(r)\to a$ when $r\to 0$.
\\
\nind If $p>\frac{N}{N-2}$, then $K>0$ and $(0,0)$ is a saddle point. The trajectory $\CT_r$  converges to $(0,0)$ at $-\infty$. There is also the unique trajectory $\CT_s$ which converges to $(0,0)$ when $t\to\infty$. Their slope at $(0,0)$ is $N-2$ and they correspond to solutions $u(r)\sim cr^{2-N}$ when $r\to\infty$. \\
\nind If $p=\frac{N}{N-2}$, then $K=0$. Besides the regular trajectory which always exists, there exists an invariant trajectory passing through $(0,0)$, with slope $N-2$, by the theorem of the central manifold. We will see later on that it converges to $(0,0)$ as $t\to-\infty$. \smallskip

\nind (ii) Assume now that $N=1$ or $2$ there still exists the regular trajectory $\CT_r$. \\
\nind If $N=1$, then $\gl_1=\frac{p+1}{p-1}>\gl_2$.  There exist infinitely many trajectories different from $\CT_r$, converging to $(0,0)$ at $-\infty$, in $Q_1$ or $Q_4$, corresponding to solutions such that $u(0)=u_0>0$ and 
$u'(0)=a\in\BBR\setminus\{0\}$. There exists one trajectory converging to $(0,0)$ at $-\infty$ with slope $-1$ and located in $Q_4$. It corresponds to solutions such that  $u(0)=0$ and $u'(0)>0$.\\
\nind If $N=2$, then $\gl_1=\gl_2=\ga$. The point $(0,0)$ is a degenerate node. All the trajectories in a neighborhood of $(0,0)$ tend to $(0,0)$ when $t\to-\infty$ and are tangent to $\xi_1$. However they behave like $c(-t)^{-1}$ for any $c>0$. They correspond to solutions $u$ such that $\displaystyle \lim_{r\to 0}(-\ln r)^{-1}u(r)=a>0$.
\subsubsection{Linearisation at the non-trivial equilibrium points}
\blemma {equili} 1- If $1<p\leq \frac{N}{N-2}$ and $M> 0$, or $1<p<\frac{N}{N-2}$ and $M=0$, $P_{_M}$ is a saddle point.\smallskip

\nind 2-  If $p> \frac{N}{N-2}$ and $M>m^*$, $P_{_{1,M}}$ is  a node point and a source and $P_{_{2,M}}$ is a saddle point.\smallskip

\nind 3-  If $p> \frac{N}{N-2}$ and $M=m^*$, $P_{_{m^*}}$ is not hyperbolic. One eigenvalue is $N-2$ and the other is $0$. 
\es
\Proof Set $y_{_M}=\ga x_{_M}$. In view of $(\ref{Z5})$, $y_{_M}$ satisfies 
\bel{W7}\BA {lll}
\ga^{-p}y^{p-1}_{_M}-My^{\frac{p-1}{p+1}}_{_M}+K=0.
\EA\ee
Setting $x=x_{_M}+\overline x$, $y=y_{_M}+\overline y$, the linearized equation at $(x_{_M},y_{_M})$ is 
\bel{W8}\BA {lll}
\overline x_t=\ga \overline x-\overline y\\
\overline y_t=-px^{p-1}_{_M}\overline x+\left(\frac {2Mp}{p+1}y_{_M}^{\frac{p-1}{p+1}}-K\right)\overline y.
\EA\ee
The characteristic polynomial of the corresponding matrix is 
\bel{W9}\BA {lll}
T_{{y_{_M}}}(X)=X^2-\left(\myfrac {2Mp}{p+1}y_{_M}^{\frac{p-1}{p+1}}-L\right)X+2K-\myfrac {2Mp}{p+1}y_{_M}^{\frac{p-1}{p+1}}.
\EA\ee
\nind 1- If $p\leq \frac{N}{N-2}$ and $M>0$, or  $p<\frac{N}{N-2}$ and $M=0$, $P_{_M}$ is unique. Since either $K\leq 0$ and $M>0$ or $K< 0$ and $M=0$, the product of the roots is negative. Hence $P_{_M}$ is a saddle point.\\
\nind 2-3- Next we assume $N\geq 3$ and $M\geq m^*$. 
The sum of the roots of $T_{{y_{_M}}}(X)$ is equal to 
$$\BA {lll}
\myfrac {2Mp}{p+1}y_{_M}^{\frac{p-1}{p+1}}-L=\myfrac {2p}{p+1}\left(\ga^{-p}y^{p-1}_{_M}+K\right)-L>\myfrac{p-1}{p+1}K+\ga>0.
\EA$$
Concerning the product $\Gp(y_{_M})$ of the roots, we deduce from $(\ref{Y6})$ that $T_{y_{m^*}}(0)=0$ hence $\Gp(y_{m^*})=0$. Since by $(\ref{W7})$
$$ \Gp(y_{_M})=2K-\myfrac {2Mp}{p+1}y_{_M}^{\frac{p-1}{p+1}}=\myfrac {2K}{p+1}-\myfrac {2p}{p+1}\ga^{-p}y_{_M}^{p-1},
$$
we infer that for $M>m^*$, 
\bel{W10}\BA {lll}\Gp (y_{_{2,M}})<\Gp(y_{m^*})=0<\Gp (y_{_{1,M}}).
\EA\ee
Hence $P_{_{2,M}}$ is a saddle point and  $P_{_{1,M}}$ is a source. In order to characterize the nature of this source we denote by $D(T_{{y_{_M}}})$ the discriminant of $T_{{y_{_M}}}$. Then 
$$\BA{lll}D(T_{{y_{_{1,M}}}})=\left(\myfrac {2Mp}{p+1}y_{_{1,M}}^{\frac{p-1}{p+1}}-L\right)^2+4\left(\myfrac {2Mp}{p+1}y_{_{1,M}}^{\frac{p-1}{p+1}}-2K\right)\\[4mm]
\phantom{D(T_{{y_{_{1,M}}}})}
=\left(\myfrac {2Mp}{p+1}y_{_{1,M}}^{\frac{p-1}{p+1}}-L+2+2\sqrt{N-1}\right)\left(\myfrac {2Mp}{p+1}y_{_{1,M}}^{\frac{p-1}{p+1}}-L+2-2\sqrt{N-1}\right).
\EA$$
By $(\ref{W7})$, $My_{_{1,M}}^{\frac{p-1}{p+1}}>K$, hence 
$$\BA {lll}\myfrac {2Mp}{p+1}y_{_{1,M}}^{\frac{p-1}{p+1}}-L+2-2\sqrt{N-1}>\myfrac {p-1}{p+1}K+\myfrac{2}{p-1}+2-2\sqrt{N-1}\\[4mm]
\phantom{\myfrac {2Mp}{p+1}y_{_{1,M}}^{\frac{p-1}{p+1}}-L+2-2\sqrt{N-1}}=\myfrac{N-2\sqrt{N-1}}{p^2-1}\left(p-\myfrac{N-2}{N-2\sqrt{N-1}}\right)^2.
\EA$$
Hence $D(T_{{y_{_{1,M}}}})>0$ which implies that the roots are real and $P_{_{1,M}}$ is a node. \\
\nind If $M=m^*$ the product of the roots is $0$, hence one root is $0$. Since their sum is equal to $\frac{2p}{p+1}m^*y_{_{m^*}}^{\frac{p-1}{p+1}}-L=N-2$, the nonzero root is equal to $N-2$. 

\qeda
\subsubsection{The vanishing curves of the vector field}
The vector field associated to $(\ref{W4})$ is defined by
\bel{W11}\BA {lll}
(x,y)\mapsto H(x,y)=(H_1(x,y), H_2(x,y)):=\left(\ga x-y,-Ky-x^p+M|y|^{\frac{2p}{p+1}}\right).
\EA\ee
We call vanishing curves of $H$ in $\BBR^2_+$ the set of points where $H_1$ or $H_2$ vanishes.
$$\CL=\left\{(x,y)\in \BBR^2_+:y=\ga x\right\},
$$
and 
$$\CC=\left\{(x,y)\in \BBR^2_+:H_2(x,y)=0\right\}=\CC_1\cup\CC_4,
$$
where
$$\CC_1=\left\{(x,y)\in Q_1:x^p={My^{\frac{2p}{p+1}}-Ky}:=\Phi(y)\right\},
$$
and
$$\CC_4=\left\{(x,y)\in Q_4:x^p=M|y|^{\frac{2p}{p+1}}-Ky:=\Psi(y)\right\}.
$$

Those vanishing curves are the boundary of some semi-invariant regions in $\BBR^2_+$. Their configuration depends on the intersection of these curves.\smallskip

\nind I- {\it If $K\leq 0$ and $M>0$ we denote by }\\
\nind (A) is the set of points $(x,y)\in Q_1$ such that $y>\max\left\{\ga x,\Phi^{-1} (x)\right\}$.\\
\nind (B) is the set of points $(x,y)\in Q_1$ such that $x\geq x_{_M}$ and $\ga x<y<\Phi^{-1} (x)$.\\
\nind (C)  is the union of the set of points $(x,y)\in Q_1$ such that $y<\min\left\{\ga x,\Phi^{-1} (x)\right\}$ and  the set of points  $(x,y)\in Q_4$ such that $x>\Psi(y)$.\\
\nind D) is the set of points $(x,y)\in Q_1$ such that $\Phi^{-1} (x)<y<\ga x$.\\
\nind (E) is the set of points $(x,y)\in Q_4$ such that $x<\Psi(y)$.\medskip

\nind II- {\it If $K> 0$ and $M>m^*$ we denote by} \\
\nind (A) is the set of points $(x,y)\in Q_1$ such that $y>\max\{\ga x,\Phi^{-1} (x)\}$.\\
\nind (B) is the set of points $(x,y)\in Q_1$ such that $x\geq x_{_{2,M}}$ and $\ga x<y<\Phi^{-1} (x)$.\\
\nind (C)  is the union of the set of points $(x,y)\in Q_1$ such that $x>\Phi(y)$ and  the set of points  $(x,y)\in Q_4$ such that $x>\Psi(y)$.\\
\nind (D) is the set of points $(x,y)\in Q_1$ such that $\Phi^{-1} (x)<y<\ga x$.\\
\nind (E) is the set of points $(x,y)\in Q_4$ such that $x<\Psi(y)$.\\
\nind (F) is the set of points $(x,y)\in Q_1$ such that $x\leq x_{_{1,M}}$ and $\ga x<y<\Phi^{-1} (x)$.\medskip

\nind III- {\it If $K> 0$ and $M=m^*$}, (D) is empty.\medskip

\nind IV- {\it If $K> 0$ and $0<M<m^*$}, (D) is empty and (B) and (F) are replaced by the set $(\tilde B)=\left\{(x,y)\in Q_1\right.$ such that $\ga x<y<\Phi^{-1} (x)\}$ (note that $\tilde B$ is connected).
 \medskip
 
 \mbox{}
   
   \nind   We present below some graphics of the vector field $H$ associated to system $(\ref{W4})$.
   We show the vanishing curves of the vector field $H$ as well as the direction of the vector field along these curves.\medskip

\begin{figure}[!h]
\begin{center}
 \includegraphics[keepaspectratio, width=9cm]{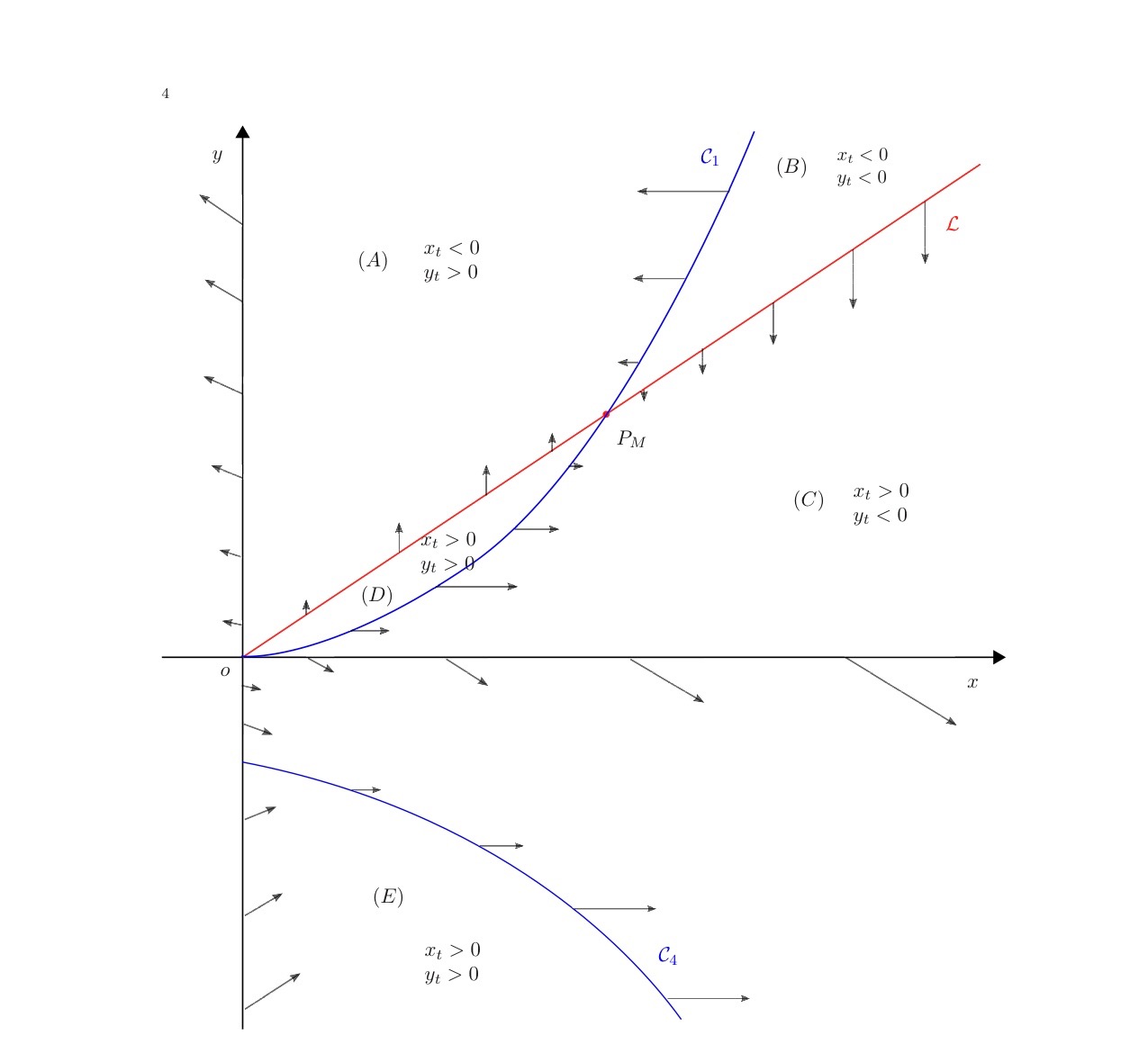}
 \end{center}
 \caption{$M>0$, $K<0\Longleftrightarrow p<\frac{N}{N-2}$.}
 \end{figure}

\mbox{   }
\begin{figure}[!h]
\begin{center}
 \includegraphics[keepaspectratio, width=9cm]{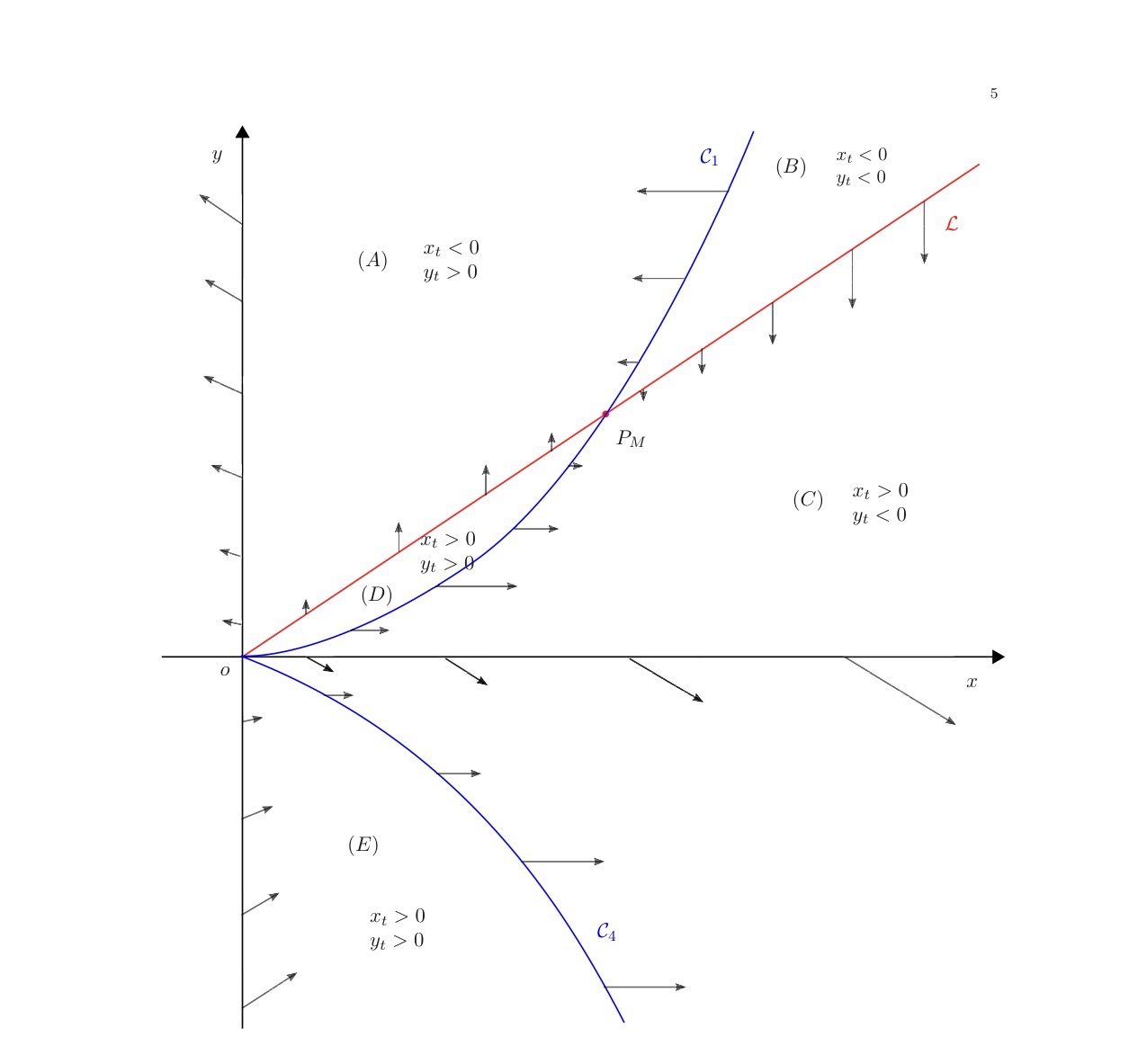}
 \end{center}
 \caption{$M>0$, $K=0\Longleftrightarrow p=\frac{N}{N-2}$.}
 \end{figure}

\mbox{   }\mbox{   }
\begin{figure}[!h]
\begin{center}
 \includegraphics[keepaspectratio, width=9cm]{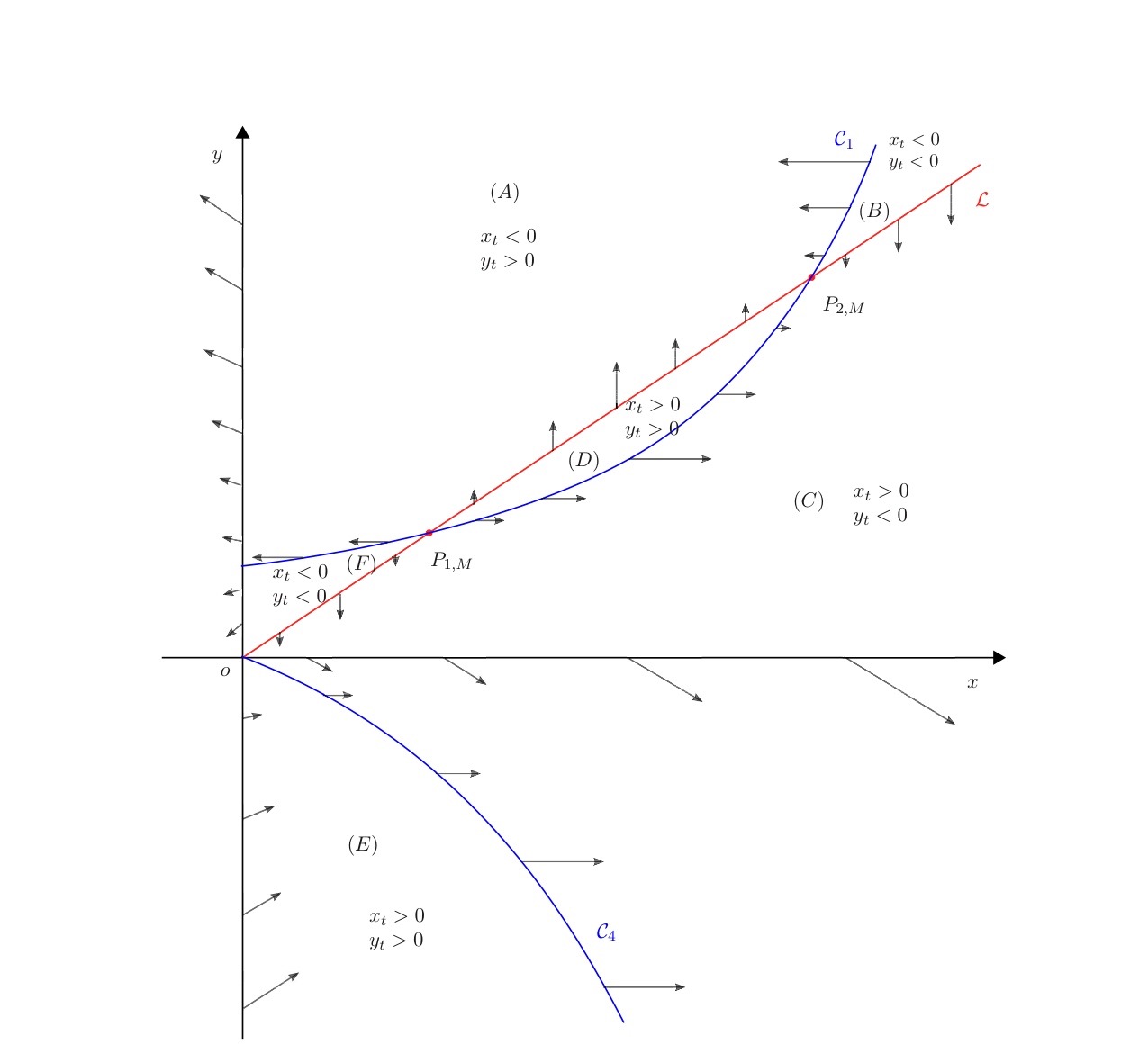}
 \end{center}
 \caption{$M>m^*$, $K>0\Longleftrightarrow p>\frac{N}{N-2}$.}
 \end{figure}

\mbox{   }
\mbox{   }\mbox{   }
\begin{figure}[!h]
\begin{center}
 \includegraphics[keepaspectratio, width=9cm]{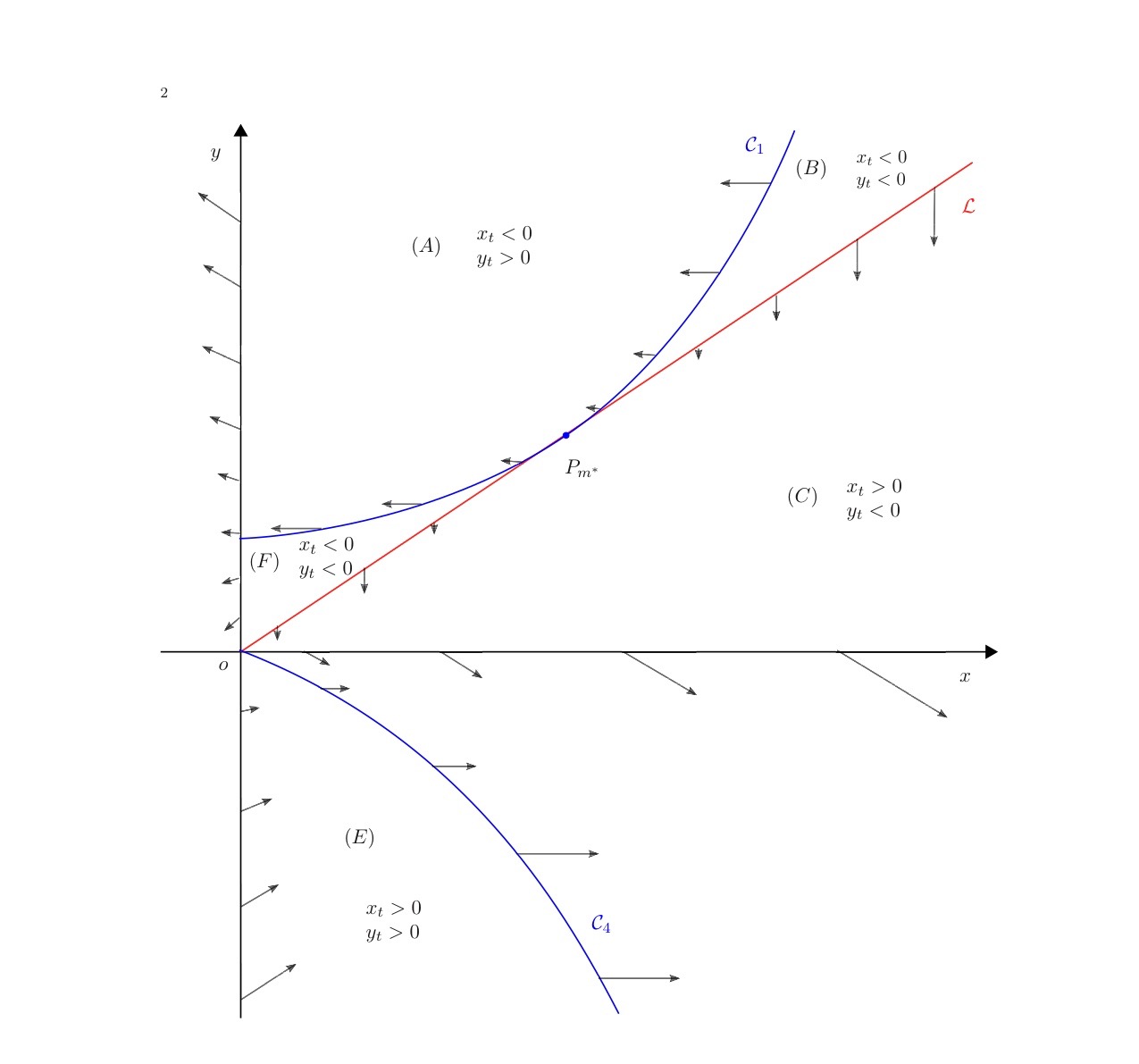}
 \end{center}
 \caption{$M=m^*$, $K>0\Longleftrightarrow p>\frac{N}{N-2}$.}
 \end{figure}
\mbox{   }
\mbox{   }
\mbox{   }\mbox{   }
\begin{figure}[!h]
\begin{center}
 \includegraphics[keepaspectratio, width=9cm]{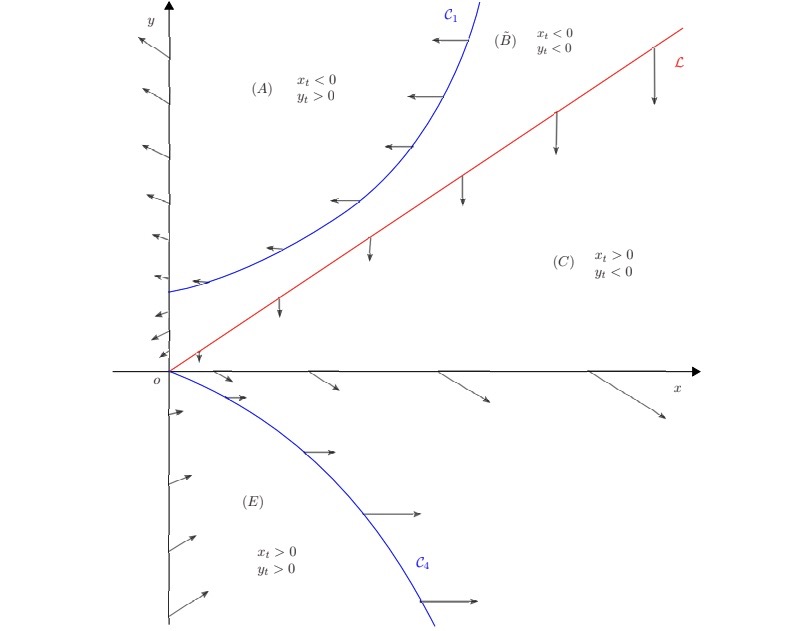}
 \end{center}
 \caption{$0<M<m^*$, $K>0\Longleftrightarrow p>\frac{N}{N-2}$.}
 \end{figure}

\newpage


\subsection{Description of the radial solutions defined near 0}
In this section we use the dynamical system $(\ref{W8})$ to describe all the positive solutions of $(\ref{W1})$ defined on a maximal interval $(0,R)$, $R\leq\infty$. The case $M=0$ which is well-known 
will be used as a comparison model. 

\subsubsection{The case $1<p<\frac{N}{N-2}$ and $M\geq 0$}
In this range of exponents the fixed point $P_{_M}$ is unique, the problem is more rigid and some of our existence and uniqueness results hold without the assumption of radiality as shown in \rth{uni}. 

\bth{subcrit-th} Let $N=1,2$ and $p>1$ or $N\geq 3$ and $p<\frac N{N-2}$, and  $M>0$.\\
\nind 1-  The function $u_{x_{_M}}$ is the unique positive solution of $(\ref{Z7})$ in $\BBR^N\setminus\{0\}$ satisfying
\bel{V1}\displaystyle 
\lim_{x\to 0}|x|^{\ga}u(x)=x_{_M}.
\ee
\nind 2-  For any $k>0$ there exists a unique  positive solution $u=u_k$ of $(\ref{Z7})$ in $\BBR^N\setminus\{0\}$ satisfying $(\ref{Z8})$. Furthermore $u_k$ is radial and
 \bel{V2}\displaystyle 
\lim_{|x|\to \infty}|x|^{\ga}u_k(x)=x_{_M}.
\ee
To this set of solutions $u_k$ is associated a unique heteroclinic orbit $\CT_1$ of the system $(\ref{W4})$ connecting the origin when $t\to -\infty$ to $P_{_M}$ when $t\to\infty$.\smallskip

\nind 3- For any $R>0$ there exists a unique positive solution of $(\ref{Z7})$ in $\BBR^N\setminus \overline B_R$ satisfying $(\ref{V2})$ and $\displaystyle \lim_{|x|\to R}u(x)=\infty$. This solution is radial. 

\nind 4-  For any $R>0$ there exists a unique positive solution of $(\ref{Z7})$ in $B_R\setminus\{0\}$ satisfying 
\bel{V3}\displaystyle 
\lim_{x\to 0}|x|^{\ga}u(x)=x_{_M}\quad\text{and }\,\lim_{|x|\to R}u(x)=\infty,
\ee
and a unique positive solution satisfying
\bel{V4}\displaystyle 
\lim_{x\to 0}|x|^{\ga}u(x)=x_{_M}\quad\text{and }\,\lim_{|x|\to R}u(x)=0.
\ee
Moreover these solutions are radial.
\smallskip

\nind 5-  Assume $N>2$. For any $k>0$ there exists $R_k>0$ and a unique radial positive solution of $(\ref{Z7})$ in $B_{R_k}\setminus\{0\}$ satisfying $(\ref{Z8})$ and vanishing on $\prt B_{B_{R_k}}$ or such that
\bel{V3'} \displaystyle
\lim_{|x|\to R_k}u(x)=\infty.
\ee
Furthermore  the mapping $k\to R_k$ is decreasing from $(0,\infty)$ onto $(0,\infty)$.
\es
\Proof 1- All the uniqueness results, which are valid not only for radial solutions, are proved in Section 4, in particular  in \rth{uni}.\smallskip

\nind 2- In the phase plane $(x,y)$, recall that $(0,0)$ is a source equilibrium and there exist infinitely many trajectories 
different from the regular one $\CT_r$, converging to $(0,0)$ when $t\to-\infty$, with the initial slope $N-2<\ga$ since $p<\frac{N}{N-2}$, so they start from $(0,0)$ in Region (D) of Figure 1. The point $P_{_M}$ is a saddle point with eigenvalues $\gl<0<\tilde\gl$ and associated eigenvectors
$$\xi_1=(1,\ga+|\gl|)\,\text{ and }\;\xi_2=(1,\ga-\tilde\gl). 
$$
We denote by $\CT_1$ the trajectory such that $x(t)$ increases and converges to $x_{_M}$  when $t\to\infty$, and by 
$\CT_2$ the trajectory such that $x(t)$ decreases and converges to $x_{_M}$  when $t\to\infty$. Their common slope is larger than $\ga$, then $\CT_1$ lies in the region $(D)$ and $\CT_2$ lies in the region $(B)$ when $t\to\infty$. Because 
$(D)$ is negatively invariant and bounded, $\CT_1$ is contained in $(D)$, a region in which $x(t)$ and $y(t)$ are monotone. Hence 
$\CT_1$ converges to a fixed point in $(\overline D)$ which is necessarily $(0,0y)$ as $t\to-\infty$. Therefore $\CT_1$ is an heteroclinic orbit joining $(0,0)$ to $P_{_M}$, and it is necessarily unique since $P_{_M}$ is a saddle point. Its slope at $(0,0)$ is $N-2$. It corresponds to a solution $u_k$ of $(\ref{W1})$ satisfying  $(\ref{Z8})$-(i). This solution $u_k$ is unique by 
\rth{uni}. Furthermore if $k<k'$ then $u_k\leq u_{k'}$. We also notice that for any $\ell>0$ and $x\in\BBR^N\setminus\{0\}$, 
$$T_\ell[u_k](r):=\ell^\ga u_k(\ell r)=u_{k\ell^{\ga+2-N}}(r).
$$
If we denote by $u_\infty$ the limit of the increasing sequence $\{u_k\}$, then 
$$T_\ell[u_\infty](r)=u_{\infty}(r)
$$
This implies that $u_\infty$ is a self-similar solution of $(\ref{W1})$, hence $u_\infty=u_{x_M}$ and $(\ref{V2})$ holds.
\smallskip

\nind 3- The trajectory $\CT_2$ converges to $P_{_M}$ when $t\to \infty$ and remains in the region $(B)$ which is negatively invariant. If $\CT_2$ were defined on whole $\BBR$, it would imply that it remains bounded because of the a priori estimate  \rprop{Oss}. But in the region $(B)$ the two functions $x(t)$ and $y(t)$ are decreasing. Hence the trajectory would converge to an equilibrium in the closure of $(B)$ different from $P_{_M}$, which is impossible. Therefore the two functions $x(t)$ and $y(t)$ with image $\CT_2$ are defined on some maximal interval $(T,\infty)$ and if we set $R=e^T$, the corresponding solution $u$ of
 $(\ref{W1})$ satisfies $\displaystyle \lim_{r\downarrow R}u(r)=\infty$. Uniqueness of a solution defined on 
 $(R,\infty)$ and blowing-up at $r=R$ follows from \rth{uni3}. By the scaling $T_\ell$, the function $u$ is transformed in a solution of $(\ref{W1})$ which blows-up at $r=\ell^{-1}R$ and which is associated with the same trajectory $\CT_2$. Hence $R$ can be any positive real number.\smallskip

\nind 4- There exist two unstable trajectories $\CT_3$ and $\CT_4$ converging to $P_{_M}$ when $t\to-\infty$. They are associated to the eigenvalue $\tilde\gl$, and their slope at $P_{_M}$ is $\ga-\tilde\gl$. We denote by $\CT_4$ the trajectory which enters in the region $(C)$. Since this region is positively invariant, $\CT_4$ remains in it. Then either 
its components are defined on some maximal interval $(-\infty,T)$ and the corresponding solution $u$ of $(\ref{W1})$ tends to $\infty$ when 
$r\uparrow R:=e^T$, or they are defined on wole $\BBR$. In that case $u$ would coincide with $u_{x_M}$ by 1- which is contradictory. Hence $u$ is defined on the maximal interval $(0,R)$. Notice also that $u$ is decreasing on some interval $(0,r_0)$ and increasing on $(r_0,R)$ by the phase plane analysis. Thanks to the scaling $T_\ell$,  $R$ can be taken arbitrarily. Since 
this solution is uniquely determined by  $\CT_4$, it is unique. This corresponds to a uniqueness result for solutions of 
$(\ref{Z1})$ in the class of radial solutions. This proves  $(\ref{V3})$\\
Consider now the trajectory $\CT_3$. It belongs to region $(A)$ when $t\to-\infty$, and in this region $x_t<0$ and $y_t>0$. Since $P_{_M}$ is the only equilibrium in the quadrant $Q_1$ the trajectory intersects the straight line 
$x=0, y>0$ at some $y_0>y_{_M}$ for some $t=T$. Hence the corresponding solution $u$ vanishes for $r=R=e^T$. Furthermore $R$ can be taken arbitrarily. This proves  $(\ref{V4})$. \smallskip

\nind 5- Since $(0,0)$ is a source, there exists $\ge>0$ such that any backward trajectory issued from 
$(x_0,y_0)\in B_\ge(0)$ converges to $(0,0)$ when $t\to-\infty$. All these trajectories in the first quadrant $Q_1$ have initial slope $N-2$. If $(x_0,y_0)\in B_\ge(0)\cap (D)$ is above the heteroclinic orbit $\CT_1$, it cannot converge to $P_{_M}$, then it crosses $\CL$, enters in $(A)$ and crosses the axis $\{x=0,y>0\}$ for some $T$. By Appendix A2  the associated solution $u$ of $(\ref{W1})$ satisfies $(\ref{Z8})$ for some $k>0$ and $u(e^T)=0$. If $(x_0,y_0)\in B_\ge(0)\cap (D)$ is below $\CT_1$, it enters the region $(C)$ which is positively invariant and for the same reasons as in Step 3 it blows-up for some $t=T$. The corresponding solution $u$ satisfies $(\ref{Z8})$ for some $k>0$ and blows up for $r=e^T:=R$. Uniqueness of this type of solutions in $B_R$ follows either from the general result \rth{uni} or from 
the uniqueness of the trajectories of the system $(W8)$. The correspondance $k\mapsto R$ is decreasing and onto from $(0,\infty)$ to $(0,\infty)$ by uniqueness and using the scaling transformation $T_\ell$. 
\qeda\medskip

\nind\Remark It follows from the analysis of the phase plane that all the positive radial solutions of $(\ref{Z7})$ defined in a neighborhood of $x=0$ or in the complement of a ball have their behaviour described by 1 or 2.  \smallskip

\nind{\it Proof of \rth{T3}}. It is a direct consequence of \rth{subcrit-th} and \rth{uni}.

\subsubsection{The case $p=\frac{N}{N-2}$, $q=\frac{2p}{p+1}=\frac{N}{N-1}$ and $M> 0$}
When $p=\frac N{N-2}$ and $M=0$ the isolated singularities of solutions of  $(\ref{Z1})$ are removable and the behaviour at infinity of these solutions is described in \cite{Veasym}. When $M>0$ it is no longer the case and the interaction of the two reaction terms yields new phenomena. The next result covers \rth{T4}, up to uniqueness which will follow from \rth{uni}.
 \bth{crit-th} Let $N\geq 3$, $p=\frac{N}{N-2}$ and $M>0$. \smallskip
 
\nind 1-  If $M=0$ any isolated singularity of a solution, not necessarily radial neither nonnegative, of $-\Gd u+|u|^{p-1}u=0$ is removable. If $u$ is any solution of this equation in $B_R^c$, there exists $\gl$ such that
 \bel{V10}\displaystyle
\lim_{|x|\to\infty}|x|^{N-2}\left(\ln|x|\right)^{\frac {N-2}{2}}u(x)=\gl,
\ee
 and $\gl$ can only take the three values $\left(\frac{N-2}{\sqrt 2}\right)^{N-2}$, $-\left(\frac{N-2}{\sqrt 2}\right)^{N-2}$ and $0$. \smallskip
 
\nind 2-  If $M>0$, the function $u_s(x)=\left((N-2)M^{\frac {N}{N-1}}\right)^{N-2}|x|^{2-N}$ is the unique positive separable solution of $(\ref{Z1})$ in $\BBR^N\setminus\{0\}$. There exists a positive radial solution $u$,  unique up to the scaling transformations $T_\ell$, satisfying 
  \bel{V11}\BA {lll}
(i) \qquad \qquad&\displaystyle\lim_{r\to 0}r^{N-2}\left|\ln r\right|^{N-1}u(r)=\myfrac{((N-1)M)^{1-N}}{N-2} \qquad \qquad\\[2mm]
  
(ii)&\displaystyle\lim_{r\to\infty}r^{N-2}u(r)=\left((N-2)M^{\frac {N}{N-1}}\right)^{N-2}. \qquad \qquad
\EA\ee
Furthermore, for any $R>0$ there exists a positive and radial solution $u:=u_{_R}$ of $(\ref{Z1})$ in $B_R\setminus\{0\}$ satisfying 
   \bel{V14}\displaystyle\lim_{r\to0}r^{N-2}u(r)=\left((N-2)M^{\frac N{N-1}}\right)^{N-2},
\ee
and    
\bel{V15}\displaystyle\lim_{|x|\to R}u(x)=\infty.\ee
Finally, there exists also a unique positive solution $u:=\tilde u_{_R}$ of $(\ref{Z1})$ in $B_R\setminus\{0\}$ satisfying 
$(\ref{V11})$-(i) and $(\ref{V15})$ or $(\ref{V11})$-(i) and $u=0$ on $\prt B_R$. In both cases the solution is radial.
 \es
 \Proof The results of assertion 1 is proved in \cite{Veasym}, \cite{BM}.\\ 
Assertion 2- Since $p=\frac {N}{N-2}$, $K=0$. As a consequence, the vanishing curve $\CC_4$ goes through $(0,0)$. \rlemma{equili} is still valid. The point $P_{_M}=(x_{_M},y_{_M})$ is a saddle point and with $x_{_M}=(N-2)^{N-2}M^{\frac {N(N-2)}{N-1}}$ and $y_{_M}=(N-2)^{N-1}M^{\frac {N(N-2)}{N-1}}$. The stable curve $\CT_1$ is an heteroclinic orbit staying in the region (D) and connecting $(0,0)$ to $P_M$. The point $(0,0)$ is no longer hyperbolic since the charcteristic values are $\gl_1=0$ and $\gl_2=N-2$, and the behaviour of the solutions in its neighbourhood is more delicate. The vector $\xi_2=(1,0)$ is the eigenvector associated to the nonzero eigenvalue $N-2$ and the unstable curve corresponds to the regular solutions $u_{[a]}$. By the central manifold theorem, the curve $\CT_1$ is the central manifold of $(0,0)$ and is tangent at this point to the eigenvector $\xi_1=(1,N-2)$. Therefore
 $\displaystyle\lim_{t\to-\infty}\tfrac{y(t)}{x(t)}=N-2$ on $\CT_1$. As a consequence
 $$\myfrac{(u(r))^{\frac{N}{N-2}}}{|u'(r)|^{\frac{N}{N-1}}}=\myfrac{(x(t))^{\frac{N}{N-2}}}{(y(t))^{\frac{N}{N-1}}}\leq c(y(t))^{\frac{N}{(N-1)(N-2)}}\quad\text{for }r\leq \ge_0.
 $$
 Consequently $u^{\frac{N}{N-2}}=o\left(|u'|^{\frac{N}{N-1}}\right)$ in a neighborhood of $r=0$. Therefore, for any $\ge>0$ there exists $r_\ge>0$ such that 
  \bel{V12}(1-\ge)Mr^{N-1}|u'|^{\frac{N}{N-1}}\leq -(r^{N-1}u')'\leq Mr^{N-1}|u'|^{\frac{N}{N-1}}\quad\text{on }\,(0,r_\ge].
\ee
 Putting $W=r^{N-1}|u'|$ these inequalities become completely integrable and we derive
  \bel{V13}(N-1)(W(r))^{-\frac{1}{N-1}}=M|\ln r|^{-1}(1+o(1))\quad\text{as }r\to 0.
\ee
 Integrating $(\ref{V13})$ implies $(\ref{V11})$-(i). Uniqueness among the radial solutions follows from the uniqueness of the stable heteroclinic orbit. As in 
 \rth{subcrit-th} the unstable trajectory of $P_{_M}$ entering  (C) intersect the axis $y=0$ at some $P_0=(x_0,0)$ and any corresponding solution $u$ is defined in some $B_R$, R>0, and it blows-up when $|x|\uparrow R$. Using the transformation $T_\ell$, $T_\ell[u]$ is a solution defined in the ball $B_{\frac R\ell}$ which blows-up for $|x|=\frac R\ell$ and still satisfies $(\ref{V14})$. The backward trajectory $\CT^-_{[P]}$ of any point $P=(P,0)$ on the seqment $(0,P_0)$ converges to 
 $(0,0)$ when $t\to-\infty$. For the same reason as for $\CT_1$ the corresponding solution $u$ satisfies $(\ref{V11})$-(i) and it blows-up ifor $r=R$ for some 
 $R>0$. Since the scaling transformation $T_\ell$ leaves $(\ref{V11})$-(i) unchanged, $\frac{R}{\ell}$ can take any value, this ends the proof.\qeda\medskip
 
 \nind\Remark

\subsubsection{The case $p>\frac{N}{N-2}$ and $M> 0$}

In the range of exponent $p>\frac N{N-2}$, the positive parameter $m^*$ defined by  $(\ref{Z6})$ plays a fundamental role. The following result covers \rth{T5} and describes all the positive solutions of $(\ref{W1})$ either defined near $\infty$ or near $0$.
\bth{supcrit} Let $N>2$ and $p>\frac N{N-2}$.\smallskip

\nind 1-  If $M>m^*$, then $u_{x_{_{1,M}}}$ and $u_{x_{_{2,M}}}$ are the two self-similar solutions. Moreover \\
\nind (i) there exists a unique, up to the transformation $T_\ell$, positive radial solution $u=u_{1,2}$ defined in $\BBR^N\setminus\{0\}$ satisfying 
  \bel{V16}\displaystyle
  \lim_{r\to 0}r^{\ga}u(r)=x_{_{1,M}}\quad\text{and }\,   \lim_{r\to \infty}r^{\ga}u(r)=x_{_{2,M}}.
  \ee
  \nind (ii) For any $k>0$ there exists a unique positive radial $u=u_{1,k}$ defined in $\BBR^N\setminus\{0\}$ satisfying 
    \bel{V17}\displaystyle
  \lim_{r\to 0}r^{\ga}u(r)=x_{_{1,M}}\quad\text{and }\,   \lim_{r\to \infty}r^{N-2}u(r)=k.
  \ee
  \nind (iii) For any $R>0$ there exists a positive radial solution $u=u_{j,R}$ in $B_R\setminus\{0\}$ with j=1,2, satisfying 
      \bel{V18}\displaystyle
  \lim_{r\to 0}r^{\ga}u(r)=x_{j,M}\quad\text{and }\,   \lim_{r\uparrow R}u(r)=\infty.
  \ee
This solution is  unique if $j=2$.
  There exists also a unique radial positive solution $\tilde u=\tilde u_{2,R}$ in $\BBR^N\setminus B_R$ satisfying
        \bel{V19}\displaystyle
  \lim_{r\to \infty}r^{\ga}\tilde u(r)=x_{_{2,M}}\quad\text{and }\,   \lim_{r\downarrow R}\tilde u(r)=\infty.
  \ee
  or 
        \bel{V18'}\displaystyle
  \lim_{r\to \infty}r^{\ga}\tilde u(r)=x_{_{2,M}}\quad\text{and }\,   \lim_{r\downarrow R}\tilde u(R)=0,
  \ee

\nind 2-  If $M=m^*$, $u_{x_{m^*}}$  is the unique self-similar solution, and statement 1-(ii)  still holds with $x_{_{1,M}}$ replaced by $x_{m^*}$ in $(\ref{V17})$. There exist infinitely many radial positive solutions in $B_R\setminus\{0\}$ satisfying $(\ref{V18})$ with $x_{_{j,M}}$ replaced by $x_{m^*}$, and at least one in $\BBR^N\setminus B_R$ satisfying $(\ref{V19})$ or $(\ref{V18'})$ with $x_{_{2,M}}$ replaced by $x_{m^*}$.
\smallskip

\nind 3-  If $0\leq M<m^*$, there exists no singular solution. For any $R>0$ there exist $k>0$ and a unique positive radial solution in $\BBR^N\setminus B_R$ satisfying
\bel{V20}\BA {lll}\displaystyle
 (i)\qquad\qquad &\displaystyle\lim_{r\to \infty}r^{N-2}u(r)=k\qquad\qquad\qquad\qquad\qquad\qquad\qquad\\[2mm]
 (ii)\qquad\qquad &\displaystyle\lim_{r\downarrow R}u(r)=\infty.
  \EA\ee
  Any positive radial solution defined in $\BBR^N\setminus B_R$ has the same asymptotic behaviour as in $(\ref{V20})$-(i).
\es
\Proof{\it Case  1: $M>m^*$}.  (i) From \rlemma{equili}-2, $P_{_{2,M}}$ is a saddle point with stable trajectories $\CT_1$, $\CT_2$, and unstable ones $\CT_3$ and $\CT_4$ defined as in the proof of \rth{subcrit-th}. The trajectory $\CT_1$ lies in $(D)$
as $t\to\infty$ and remains in $(D)$ for all $t$ because $(D)$ is negatively invariant. Hence it converges to $P_{_{1,M}}$ when $t\to-\infty$. Therefore $\CT_1$ is an heteroclinic orbit connecting $P_{_{1,M}}$ to $P_{_{2,M}}$. It is unique and it corresponds to a solution $u$ satisfying $(\ref{V16})$, thus $u$ is unique up to the scaling transformations $T_\ell$ for $\ell>0$. \smallskip

\nind (ii) The point $(0,0)$ is a saddle point with unstable trajectory $\CT_r$ and stable trajectory $\CT_s$ which converges to $(0,0)$ as $t\to\infty$ with initial slope $N-2$. To $\CT_s$ are associated the solutions $u$ of 
$(\ref{W1})$ satisfying $\displaystyle\lim_{r\to \infty}r^{N-2}u(r)=k$, this solution is unique for fixed $k$ and denoted by $u_k$. Since $N-2>\ga$,  this stable trajectory lies in the region $(F)$ at infinity. Since $(F)$ is negatively invariant, the two functions $x(t)$ and $y(t)$ are decreasing and thus $\CT_s$ converges to $P_{_{1,M}}$  when $t\to-\infty$. Hence 
$\CT_s$ is a heteroclinic orbit connecting $P_{_{1,M}}$ to $(0,0)$ and it is unique. To this trajectory is associated a solution $u$ of $(\ref{W1})$ satisfying $(\ref{V17})$ and unique up to the transformations $T_\ell$. \smallskip

\nind  (iii) The unstable trajectory $\CT_4$ of $P_{_{2,M}}$ enters the region (C), crosses the axis $0x$ and blows-up in finite time as in \rth{subcrit-th}. Since $P_{_{1,M}}$ is a source and a node, there exist trajectories different from $\CT_1$ converging to $P_{_{1,M}}$ when $t\to-\infty$ and with a slope at this point smaller than 
 $\ga$. Consider one of them below $\CT_1$ near $P_{1,M}$; either it enters the region (C), then intersects the axis $0x$ and finally blows-up, or it enters the region (D), but since it cannot converge to $P_{_{2,M}}$, it leaves (D) and finally blows up as in the first case. In any case such a trajectory corresponds 
 to a solution which satisfies $(\ref{V18})$ with $j=1$. Because of the scaling invariance of the condition, $R$ can take any positive value. Notice that since there may exist several trajectories converging to $P_{1,M}$ at $-\infty$ with the same slope at this point, the corresponding solution $u_{1,R}$ is not unique for $R$ fixed. \\
 As in the proof of \rth{subcrit-th} $\CT_3$ corresponds to a solution satisfying $(\ref{V18'})$. The stable trajectory 
 $\CT_2$ of $P_{_{2,M}}$ lies in the region $(B)$ near $\infty$. Since this region is negatively invariant the trajectory remains in it, hence $x(t)$ and $y(t)$ are decreasing. If they were defined on $\BBR$, they would remain bounded by \rprop {Oss} and the trajectory would converge to a fixed point in $(B)$, different from $P_{_{2,M}}$. Since such a point does not exist the functions  $x(t)$ and $y(t)$ are defined on a maximal interval $(T,\infty)$ and they blow-up when $t\downarrow T$. To this traectory is associated a solution $\tilde u$ of $(\ref{W1})$ satisfying $(\ref{V19})$ with $R=e^T$.
The trajectory $\CT_2$ is unique thus $R$ can be fixed arbitrarily by using the scaling transformation $T_\ell$.\smallskip

\nind {\it Case  2: $M=m^*$}. There exists a unique nontrivial equilibrium $P_{_{m^*}}$. To the eigenvalue $0$ is associated the eigenvector $\left(1,\ga\right)$, while to the eigenvalue $N-2$ is associated the eigenvector  is $\left(1,-K\right)$. There exist two trajectories $\CT_3$ and $\CT_4$ converging to 
 $P_{_{m^*}}$ when $t\to -\infty$. The trajectory $\CT_4$ with slope $\ga+2-N$ at $P_{_{m^*}}$ enters the region (C), crosses the axis $0x$ and blows-up in finite time. It corresponds to a solution $u$ with  $\displaystyle \lim_{r\to0}r^{\ga}u(r)=x_{_{m^*}}$ when $r\to 0$ and which blows-up at $r=R$. \\
The point $(0,0)$ is a saddle point. The stable manifold $\CT_1$ has initial slope $N-2$. It corresponds to a trajectory which converges to 
 $(0,0)$ when $t\to\infty$. Since the region (F) is negatively invariant this trajectory converges to $P_{_{m^*}}$ when $t\to-\infty$, and its slope at this point is 
 $\ga$. Hence $\CT_1$ is the central manifold at $P_{_{m^*}}$. As in case (1) this trajectory corresponds to a positive solution $u$ in $\BBR^N\setminus\{0\}$ which satisfies 
 \bel{V21}\displaystyle
  \lim_{r\to 0}r^{\ga}u(r)=x_{_{m^*}}\quad\text{and }\,    \lim_{r\to \infty}r^{N-2}u(r)=k,
  \ee
  for some $k>0$.\\
  Moreover, any trajectory which has one point in the bounded negatively invariant region delimited by $\CT_s$, $\CT_4$ and the axis $0x$, converges to $P_{m^*}$ when $t\to-\infty$, tangentially to the line $\CL$. Since it cannot converge to $(0,0)$, it crosses the axis $0x$ in finite time and it blows up for $t=T=\ln R$. This corresponds to a positive solution $u$ of $(\ref{W1})$ in $B_R\setminus\{0\}$ which satisfies  $(\ref{V18})$ with $x_{j,M}$ replaced by 
  $x_{m^*}$. \\
 {\it We claim that there exists at least one trajectory belonging to the central manifold at $P_{m^*}$ which converges to 
  $P_{m^*}$ when $t\to\infty$ and blows up in finite time}: the backward trajectory $\CT_P$ of any $P\in\CC_1\cap\{(x,y):x>x_{m^*}\}$, belongs locally to $(B)$ for $t<0$ since this region is negatively invariant. Furthermore its coordinates satisfy $x(t)>x_{m^*}$ and $y(t)>y_{m^*}$ for $t<0$. Next, the backward trajectory $\CT_P$ of any $P\in\CL\cap\{(x,y):x>x_{m^*}\}$ belongs to $(B)$ and its coordinates satisfy also $x(t)>x_{m^*}$ and $y(t)>y_{m^*}$ for $t<0$ and $x(t)>x_{m^*}$ for $t>0$. Let  $\CU$ be the set of points $P\in (B)$ such that $\CT_P$ crosses $\CC_1\cap\{(x,y):x>x_{m^*}\}$ for some $t>0$ and $\CV$ the set of points $P\in (B)$ such that $\CT_P$ crosses $\CL\cap\{(x,y):x>x_{m^*}\}$ for some $t>0$. By standard transversality arguments $\CU$ and $\CV$ are open and disjoint. 
Since $(B)$ is connected, it cannot be the union of the two sets   $\CU$ and $\CV$. Hence there exists $P_0$ in 
$(B)\setminus\{\CU\cup \CV\}$. By monotonicity, $\CT_{P_0}$ converges to $p_{m^*}$ when $t\to\infty$. Clearly this trajectory cannot be defined on whole $\BBR$ by \rprop{Oss}, hence it blows-up for $t=T=e^R$. This proves the existence of a solution $u$ which satisfies 
           \bel{V21+}\displaystyle
  \lim_{r\to \infty}r^{\ga}u(r)=x_{m^*}\quad\text{and }\,   \lim_{r\downarrow R}u(R)=0.
  \ee
\smallskip

\nind {\it Case  3: $0<M<m^*$}. There exists no equilibrium besides $(0,0)$ which is a saddle point with unstable trajectory $\CT_r$ and stable trajectory $\CT_s$ with initial slope $N-2>\ga$. The region $(\tilde B)$ between $\CC_1$ and $\CL$ is negatively invariant, hence $\CT_s$ remains in it and its two coordinate functions are decreasing and necessarily unbounded. The corresponding solution $u$ of $(\ref{W1})$ cannot be defined for all $r>0$ because of \rprop{Oss}, hence it blows-up for $r\downarrow R$. This proves $(\ref{V20})$.\qeda\medskip

\nind\Remark It is noticeable that in the case $m=m^*$, the equilibrium $P_{m^*}$ is not hyperbolic and the central manifold there  consist in curves with the same slope at $P_{m^*}$ but one is converging to this point when $t\to-\infty$ while the other (may be there are many) converges when $t\to \infty$.

\mysection{The radial case for $q\neq\frac{2p}{p+1}$}
In this section we study the nonnegative solutions of 
\bel{U0}
-u''-\myfrac{N-1}{r}u'+|u|^{p-1}u-M|u'|^{q}=0,
\ee
when $q\neq\frac{2p}{p+1}$. 

\subsection{Non-autonomous systems associated to the equation}
Since $q\neq \frac {2p}{p+1}$ there exists no autonomous 2-dimensional system in which equation $({\ref{W1}})$ can be transformed. The systems that we introduce below are suitable for specific range of singular phenomena characteristic of one of the following equations $\CL_pu=0$, $\CR^M_qu=0$ and $\CE^M_{p,q}u=0$.
\subsubsection{System describing the behaviour of Emden-Fowler equation}
We set 
\bel{U1}\displaystyle  
u(r)=r^{-\frac 2{p-1}}x(t)=r^{-\ga}x(t)\,,\;u'(r)=-r^{-\frac {p+1}{p-1}}y(t)=-r^{-\ga-1}y(t)\,,\; t=\ln r.
\ee
If $u$ is a positive solution of $(\ref{W1})$ there holds
\bel{U2}\BA {lll}\displaystyle  
x_t=\ga x-y\\
y_t=-Ky-x^p+Me^{-\frac{\gs t}{p-1}}|y|^{q},
\EA\ee
where, we recall it, $\sigma$ is defined in $(\ref{Usigma})$.
Equivalently
\bel{U3}\BA {lll}\displaystyle  
x_{tt}+Lx_t-\ga K x-x^p+Me^{-\frac{\gs t}{p-1}}|\ga x-x_t|^q=0,
\EA\ee
where $K=N-2-\ga$ and $L=K-\ga$. If $M=0$ this is the system which describes the radial solutions of  $\CL_pu=0$.
\subsubsection{System describing the behaviour of the Riccati equation}
We set 
\bel{U4}\displaystyle  
u(r)=r^{-\frac {2-q}{q-1}}\gx(t)=r^{-\gb}\gx(t)\,,\;u'(r)=-r^{-\frac {1}{q-1}}\eta(t)=-r^{-\gb-1}\eta(t)\,,\; t=\ln r.
\ee
If $u$ is a positive solution of $(\ref{W1})$, $(\gx,\eta)$ satisfies the system
\bel{U5}\BA {lll}\displaystyle  
\xi_t=\gb \xi-\eta\\
\eta_t=-\gk \eta-e^{\frac{\gs t}{q-1}}\xi^p+M|\eta|^{q},
\EA\ee
where $\gk=N-\gb$ is defined at $(\ref{Z18})$. The system admits a unique nontrivial equilibrium with $\xi\geq 0$ if and only if $\frac{N}{N-1}<q<2$: it is 
\bel{U5-1}\BA {lll}\displaystyle  
(\xi_{_M},\eta_{_M})=(\xi_{_M},\gb\xi_{_M})\quad\text{with }\;\xi_{_M}=\myfrac{1}{\gb}\left(\myfrac{\gk}{M}\right)^{\gb+1}.
\EA\ee
The system $(\ref{U5})$ is equivalent to
\bel{U6}\BA {lll}\displaystyle  
\gx_{tt}+(N-2-2\gb)\gx_t-\gb \gk \gx-e^{\frac{\gs t}{q-1}}\xi^p+M|\gb \xi-\xi_t|^q=0.
\EA\ee
According to the sign of $\gs$ this system is a perturbation at $-\infty$ or at $\infty$ of 
\bel{U7}\BA {lll}\displaystyle  
\xi_t=\gb \xi-\eta\\
\eta_t=-\gk \eta+M|\eta|^{q},
\EA\ee
which describes the radial positive solutions of $\CR^M_qu=0$.
\subsubsection{System describing the behaviour of the eikonal equation}
Assuming $p\neq q$, we set
\bel{U8}\displaystyle  
u(r)=r^{-\frac {q}{p-q}}X(t)=r^{-\gg}X(t)\,,\;u'(r)=-r^{-\frac {p}{p-q}}Y(t)=-r^{-\gg-1}Y(t)\,,\; t=\ln r.
\ee
Then if $u$ is a positive solution of $(\ref{W1})$, there holds
\bel{U9}\BA {lll}\displaystyle  
X_t=\gg X-Y\\
Y_t=\gth Y+e^{-\frac{\gs t}{p-q}}(M|Y|^{q}-X^p),
\EA\ee
where $\gth=\gg+2-N$ is defined at $(\ref{Z18})$. Equivalently
\bel{U10}\BA {lll}\displaystyle  
X_{tt}+(N-2-2\gg)X_t+\gth\gg X+e^{-\frac{\gs t}{p-q}}\left(M|\gg X-X_t|^{q}-X^p\right)=0.
\EA\ee
According to the sign of $\gs$ this equation is a perturbation at $-\infty$ or at $\infty$ of 
$$M|\gg X-X_t|^{q}-X^p=0
$$
which corresponds to the eikonal equation $\CE^M_{p,q}u=0$. We note that in the case $q=\frac{N-2}{N-1}p$, then $\gg=N-2$, there exists an explicit radial solution of $(\ref{Z1})$ which is 
\bel{U11}\BA {lll}\displaystyle  
u^*_{_{M,p}}(r)=C^*_pr^{2-N}\quad\text{with }\; C^*_p=\left(M(N-2)^{\frac{N-2}{(N-1)p}}\right)^{\frac{N-1}{p}}.
\EA\ee
The function $u^*_{_{M,p}}$ is harmonic and satisfies $\CE^M_{p,q}u^*_{_{M,p}}=0$. This solution which has already been noticed in \cite{RiVe}
 will be useful in the sequel. 
 \medskip
 
 \nind\Remark The following relations between the solutions of the systems $(\ref{U2})$, $(\ref{U5})$ and $(\ref{U9})$ hold,
\bel{U13}\BA {lll}\displaystyle  
(i)\qquad\qquad &u(r)=r^{-\ga}x(t)=r^{-\gb}\xi(t)=r^{-\gg}X(t)\qquad\qquad\qquad\qquad\\[2mm]
(ii)\qquad\qquad &u'(r)=-r^{\frac{p+1}{p-1}}y(t)=-r^{\frac{1}{q-1}}\eta(t)=-r^{\frac{p}{p-q}}Y(t),
\EA\ee
which implies
\bel{U14}\BA {lll}\displaystyle  
(i)\qquad\qquad &\xi(t)=e^{-\frac{\gs}{(q-1)(p-q)}t}X(t)=e^{-\frac{\gs}{(q-1)(p-1)}t}x(t)\qquad\qquad\qquad\\[2mm]
(ii)\qquad\qquad &\eta(t)=e^{-\frac{\gs}{(q-1)(p-q)}t}Y(t)=e^{-\frac{\gs}{(q-1)(p-1)}t}y(t).
\EA\ee
This yields the following relations
\bel{U14r}\BA {lll}\displaystyle  
x^{p-1}=X^{p-q}\xi^{q-1}\quad\text{and }\;y^{p-1}=Y^{p-q}\eta^{q-1}.
\EA\ee
\subsubsection{ Lyapounov and slope functions}
There are several functions the variation of which along trajectories will be analyzed in the sequel. They are specific to the change of variable we use. The most surprising one is the function $E$ described below.
\blemma {Lyap1}Let $N\geq 1$, $p,q>1$, $p\neq q$. We define $E$ on $\BBR_+\ti\BBR\ti\BBR$ by
\bel{U15}\BA {lll}\displaystyle  
E(X,Y,t)=\myfrac{X^{p+1}}{p+1}-M\gg^q\myfrac{X^{q+1}}{q+1}-e^{\frac{\gs t}{p-q}}\left(\myfrac{(\gg X-Y)^2}{2}+\myfrac{\gg\gth X^2}{2}\right).
\EA\ee
If $u$ is a positive solution of $(\ref{W1})$ and $X$ and $Y$ are defined from $(\ref{U8})$, set
\bel{U16}\BA {lll}\displaystyle  
\CE(t)=E(X(t),Y(t),t).
\EA\ee
Then
\bel{U17}\BA {lll}\displaystyle  
\CE_t(t)=-M\left(\gg X-Y\right)\left(\gg^q X^q-|Y|^q\right)\\[4mm]
\phantom{----------}-e^{\frac{\gs t}{p-q}}\left(\left(\myfrac \ga 2+\gg +\gth\right)\left(\gg X-Y\right)^2+\myfrac{\gs\gg\gth}{2(p-q)}X^2\right),
\EA\ee
where $(X,Y)=(X(t),Y(t))$.
\es
\Proof There holds $Y=\gg X-X_t$ and 
$$X_{tt}=-\theta\gg X-(N-2-2\gg)X_t-e^{-\frac{\gs t}{p-q}}(M|Y|^q-X^p).
$$
Multiplying by $e^{\frac{\gs t}{p-q}}X_t$ we get
$$e^{\frac{\gs t}{p-q}}\myfrac {d}{dt}\left(\myfrac{X_t^2}{2}+\gth\gg \myfrac {X^2}{2}\right)-\myfrac {d}{dt}\left(\myfrac{X^{p+1}}{p+1}.
\right)=(\gg+\gth)e^{\frac{\gs t}{p-q}}X_t^2-M|Y|^qX_t.
$$
Putting 
$$\CF(t)=\myfrac{X^{p+1}}{p+1}-e^{\frac{\gs t}{p-q}}\left(\myfrac{X_t^2}{2}+\gth\gg \myfrac {X^2}{2}\right),
$$
we obtain
$$\CF_t(t)=-e^{\frac{\gs t}{p-q}}\left(\myfrac{\gs}{p-q}\left(\myfrac{X_t^2}{2}+\gth\gg \myfrac {X^2}{2}\right)+\left(\gg+\gth\right)X_t^2\right)
+M|Y|^qX_t.
$$
Since
$$\CE(t)=\CF(t)-M\gg^q\myfrac{X^{q+1}}{q+1},
$$
and $X_t=\gg X-Y$, we obtain
$$\CE_t(t)=-M(\gg X-Y)(\gg^qX^q-|Y|^q)-e^{\frac{\gs t}{p-q}}\left(\left(\myfrac{\gs}{2(p-q)}+\gg+\gth\right)X_t^2+\myfrac{\gs\gg\gth}{2(p-q)}X^2\right).
$$
and $(\ref{U17})$ follows.\qeda\medskip

The slope of a  trajectory has shown its importance in the previous section when studying solutions of $(\ref{W1})$ near an equilibrium. We introduce it as a Lyapounov type function the variations of which will be of particular interest for studying solutions of eikonal type. 

\bdef{slope} The slope of a solution $u$ is 
\bel{U18}\BA {lll}\displaystyle  
S(t)=-\myfrac{ru'(r)}{u(r)}=\myfrac{y(t)}{x(t)}=\myfrac{\eta(t)}{\gx(t)}=\myfrac{Y(t)}{X(t)}.
\EA\ee
Since $\frac {S_t}{S}=\frac {\eta_t}{\eta}-\frac {\xi_t}{\xi}$, there holds
\bel{U19}\BA {lll}\displaystyle  
S_t=S(S+2-N)+\xi^{q-1}M|S|^q-x^{p-1}\\
\phantom{S_t}=S(S+2-N)+\xi^{q-1}(M|S|^q-X^{p-q})\quad\text{if }q\neq p.
\EA\ee
Note that $S>0$ if $Y>0$.
\es
\subsection{Asymptotic estimates for the Riccati equation}
The next lemma deals with estimates near $0$ (resp. $\infty$) of radial subsolutions (resp. radial supersolutions) of the equation $\CR_q^Mu=0$, which reduces to
\bel{T0}\BA {lll}\displaystyle  
\CR_q^Mu=-u''-\myfrac{N-1}{r}u'-M|u'|^q=0
\EA\ee
in the radial case.
\blemma{subsup} Assume $N\geq 1$, $q>1$ and $M>0$. \smallskip

\nind 1- Let $u$ be any $C^2$ radial decreasing function satisfying $\CR_q^Mu\leq 0$ near $0$. \\
\nind (i) If $q>\frac{N}{N-1}$, then
\bel{T1}\BA {lll}\displaystyle  
\liminf_{r\to 0}r^{\frac{1}{q-1}}|u'(r)|\geq \left(\myfrac{\gk}{M}\right)^\frac{1}{q-1}.
\EA\ee
Therefore 
\bel{T2}\BA {lll}\displaystyle  
\liminf_{r\to 0}r^{\gb}u(r)\geq \myfrac 1\gb\left(\myfrac{\gk}{M}\right)^{\gb+1}\quad\text{if }\,q<2,
\EA\ee
\bel{T3}\BA {lll}\displaystyle  
\liminf_{r\to 0}|\ln r|^{-1}u(r)\geq \myfrac{N-2}{M}\quad\text{if }\,q=2\text{ and }\;N>2.
\EA\ee
\nind (ii) If $q=\frac{N}{N-1}$, then
\bel{T4}\BA {lll}\displaystyle  
\liminf_{r\to 0}r^{N-1}|\ln r|^{N-1}|u'(r)|\geq ((N-1)M)^{1-N}.
\EA\ee
Therefore 
\bel{T5}\BA {lll}\displaystyle  
\liminf_{r\to 0}r^{N-2}|\ln r|^{N-1}u(r)\geq \myfrac {((N-1)M)^{1-N}}{N-2}\quad\text{if }\,N\geq 3,
\EA\ee
\bel{T6}\BA {lll}\displaystyle  
\liminf_{r\to 0}(\ln|\ln r|)^{-1}u(r)\geq \myfrac{1}{M}\quad\text{if }\,N=2.
\EA\ee
(iii) If $N=1$, or if $N\geq 2$ and $1<q<\frac{N}{N-1}$, then $r^{N-1}|u'(r)|$ admits a limit belonging to $(0,\infty]$. Therefore if $N\geq 3$, $r^{N-2}u(r)$ admits a limit $c$ belonging to $(0,\infty]$. If $N=2$, $r^{N-2}u(r)$ has to be replaced by $|\ln r|^{-1}u(r)$ and if $N=1$ by $u(r)$ in the previous expression.
\smallskip

\nind 2- Let $u$ be any $C^2$ radial decreasing function satisfying $\CR_q^Mu\geq 0$  near $0$. Then all the previous statements $(\ref{T1})$--$(\ref{T6})$ are valid,   provided $\geq$ is replaced by $\leq$,  $\liminf$ by $\limsup$ and, in case (iii), the limit $c$ belongs to $[0,\infty)$. Furthermore if $q>2$ the function $u$ is bounded. \smallskip

\nind 3- Let $u$ be any $C^2$ radial decreasing function satisfying $\CR_q^Mu\leq 0$ in $B_R^c$ and tending to $0$ at infinity, and assume $q>\frac N{N-1}$. Then\\
\nind (iv)  either $q<2$ and
\bel{T7}\BA {lll}\displaystyle  
r^{\frac{1}{q-1}}|u'(r)|\geq \left(\myfrac{\gk}{M}\right)^{\frac{1}{q-1}}\quad\text{and }\; r^{\gb}u(r)\geq \myfrac{1}{\gb}\left(\myfrac{\gk}{M}\right)^{\frac{1}{q-1}},
\EA\ee
for $r$ large enough, \\

\nind (v)  or $N>2$ and
\bel{T8}\BA {lll}\displaystyle  
\lim_{r\to \infty}r^{N-1}|u'(r)|=c\geq 0\quad\text{and }\;\lim_{r\to \infty}r^{N-2}u(r)=C\geq 0.
\EA\ee
\smallskip

\nind 4- Let $u$ be any $C^2$ radial decreasing function satisfying $\CR_q^Mu\geq 0$ in $B_R^c$ and tending to $0$ at infinity. Then $q>\frac{N}{N-1}$, and either 
$N\geq 3$ and $\displaystyle \lim_{r\to \infty}r^{N-1}|u'(r)|=c> 0$, or 
\bel{T9}\BA {lll}\displaystyle  
r^{\frac{1}{q-1}}|u'(r)|\leq \left(\myfrac{\gk}{M}\right)^{\frac{1}{q-1}};\quad\text{and if $q<2$, then}\; r^{\gb}u(r)\leq \myfrac{1}{\gb}\left(\myfrac{\gk}{M}\right)^{\gb+1},
\EA\ee
for $r$ large enough.
\es
\Proof If $u$ is a radial decreasing subsolution (resp. supersolution), there holds
$$(r^{N-1}u')'+Mr^{N-1}|u'|^q\geq 0\qquad\text{(resp. $\leq 0$)}.
$$
Set $W(r)=-r^{N-1}u'(r)=r^{N-1}|u'(r)|$, then
$$Mr^{-(N-1)(q-1)}-W^{-q}W'\geq 0\qquad\text{(resp. $\leq 0$)}.
$$
Hence the function 
\bel{T10}r\mapsto\phi(r)=\left\{\BA {lll}W^{1-q}(r)-\frac{M}{\gk}r^{N-(N-1)q}\quad&\text{if }\gk=\frac{(N-1)q-N}{q-1}\neq 0\\[2mm]
W^{1-q}(r)+\frac{M}{N-1}\ln r\quad&\text{if }\gk=0,
\EA\right.
\ee
is nondecreasing (resp. nonincreasing). Notice in particular that if $u$ is a decreasing radial solution, there holds
\bel{T11}
|u'(r)|=\left\{\BA {lll}r^{1-N}\left(C+\frac{M}{\gk}r^{N-(N-1)q}\right)^{-\frac{1}{q-1}}\quad&\text{if }\gk\neq 0\\[2mm]
r^{1-N}\left(C-\frac{M}{N-1}\ln r\right)^{-\frac{1}{q-1}}\quad&\text{if }\gk=0,\EA\right.
\ee
and the estimate on $u$ follows by integration since $\gb=\frac{2-q}{q-1}$.\smallskip

\nind 1- If $u$ is a decreasing subsolution, $\phi$ is nondecreasing. \\
(i)- If $\gk>0$, then  
$$W^{1-q}(r)\leq \frac{M}{\gk}r^{N-(N-1)q}+c_0\quad\text{for }0<r\leq r_0,
$$
where $c_0=W^{1-q}(r_0)-\frac{M}{\gk}r_0^{N-(N-1)q}$. Since $N-(N-1)q<0$, $(\ref{T1})$ follows.\\
(ii)- If $\gk=0$, then for any $\ge>0$, there exists $r_\ge>0$ such that 
$$W^{1-q}(r)\leq \left(\frac{M}{N-1}+\ge\right)|\ln r|\quad\text{for }0<r\leq r_\ge,
$$
and $(\ref{T4})$ follows.\\
(iii)- If $\gk<0$, then $r^{N-(N-1)q}\to 0$ as $r\to 0$. Therefore $W(r)$ admits a  limit belonging to $(0,\infty]$ when $r\to 0$. We derive the estimates on $u$ by integration. \smallskip

\nind 2- If $u$ is a decreasing supersolution the results follow in the same way. If $q>2$, the estimate
$$\displaystyle\limsup_{r\to 0}r^{\frac{1}{q-1}}|u'(r)|\leq \left(\frac{M}{\gk}\right)^{\frac{1}{q-1}},
$$
implies that $u$ is bounded near $0$.\smallskip

\nind 3- If $u$ is a decreasing subsolution in an exterior domain, the function $\gf$ defined in $(\ref{T10})$ is nondecreasing, hence it admits a limit $\gn$ in $(-\infty,\infty]$.\\
(iv)- If 
$\gk>0$, then $r^{N-(N-1)q}\to 0$ as $r\to\infty$, hence $W^{1-q}(r)\to\gn\in [0,\infty]$, therefore $r^{N-1}|u'(r)|\to c\in [0,\infty]$. If $\gn\in (0,\infty]$, then 
$c\in [0,\infty)$. Since $u$ tends to $0$ at infinity, we obtain $r^{N-2}u(r)\to \frac{c}{N-2}$ when $r\to\infty$. \\
If $\gn=0$, then $W^{1-q}(r)\leq \frac M\gk r^{N-(N-1)q}$. This yields the estimate from below $(\ref{T7})$ of $|u'(r)|$, and therefore for $u(r)$ if $q<2$. \\
If $q\geq 2$, we obtain $|u'(r)|\geq cr^{-\frac{1}{q-1}}$, and we derive a contradiction since  $r^{-\frac{1}{q-1}}$ is not integrable at infinity. \smallskip

\nind 4- If $u$ is a decreasing supersolution in an exterior domain, then $r^{N-1}|u'(r)|$ is nondecreasing. Hence there exists $c>a$ such that 
$r^{N-1}|u'(r)|\geq a$, which implies $u(r)\geq cr^{2-N}$ for some $c>0$. Since the function $\gf$ is nonincreasing, it admits a limit $\gn$ belonging to 
$[-\infty,\infty)$. If $\gk>0$ and because $r^{N-(N-1)q}\to 0$, it follows that $\gn\in [0,\infty)$. If $\gn>0$, then $r^{N-1}|u'(r)|$ has a limit in $(0,\infty)$, and this implies that $r^{N-2}u(r)$ admits a  limit in $(0,\infty)$ at infinity. If $\gn=0$, then $W^{1-q}(r)-\frac M\gk r^{N-(N-1)q}\geq 0$ and we obtain $(\ref{T9})$. If $\gk\leq 0$, then $\gf(r)\to \infty$ as $r\to\infty$, contradiction.\qeda 

\subsection{Estimates near $0$}
In this paragraph we prove \rth{T6} and \rth{T7}.
\subsubsection{The case $\frac {2p}{p+1}<q<p.$}
Here we prove  \rth{T6}. If $u$ is a positive solution of $(\ref{W1})$ unbounded near $0$, then $u'<0$, hence the variable $X$ and $Y$ defined in $(\ref{U8})$ satisfy
\bel{S1}\BA {lll}\displaystyle  
X_t=\gg X-Y\\
Y_t=\gth Y+e^{-\frac{\gs t}{p-q}}(MY^{q}-X^p),
\EA\ee
where, we recall it, $\gs=(p+1)q-2p>0$ and $\gth=\frac{(N-1)q-(N-2)p}{p-q}$. Since $q<p$, $X(t)$ remains bounded when $t\to-\infty$. The difficulty 
comes from the fact that the term $e^{-\frac{\gs t}{p-q}}$ tends to infinity when $t\to-\infty$. 
\blemma{conv} Assume $\frac {2p}{p+1}<q<p$. If $u$ is a positive solution of $(\ref{W1})$ in $B_R\setminus\{0\}$ such that $u'<0$, then $r^\gg u(r)$ admits a limit when $r\to 0$ which can take only  the values $X_{_M}$ or $0$.
\es
\Proof We use the function $\CE$ introduced in $(\ref{U16})$. Because of 
\rprop{Oss} and \rprop{Osgrad}, $X$ and $Y$ are bounded. By assumption $\gs$ is nonnegative, hence $\CE(t)$ is bounded when $t\to -\infty$.  Using   $(\ref{U17})$ we have that 
$$\displaystyle  
\CE_t(t)+M(\gg X-Y)(\gg^q X^q-Y^q)=-e^{\frac{\gs t}{p-q}}\left(\left(\myfrac{\ga}{2}+\gg+\gth\right)(\gg X-Y)^2+\myfrac{\gs\gg\gth}{2}X^2\right),
$$
which implies that
\bel{S2}\BA {lll}\displaystyle  
-\frac{C_2\gs }{p-q}e^{\frac{\gs t}{p-q}}\leq \CE_t(t)+M(\gg X-Y)(\gg^q X^q-Y^q)\leq \frac{C_2\gs }{p-q}e^{\frac{\gs t}{p-q}},
\EA\ee
for some $C_2>0$. Because $(\gg X-Y)(\gg^q X^q-Y^q)\geq 0$, we deduce that the function $t\mapsto \CE(t)-C_2e^{\frac{\gs t}{p-q}}$ is decreasing, therefore it admits a finite limit $\Gl$ when $t\to-\infty$, and $\Gl$ is also the limit of $\CE(t)$.  Hence 
\bel{S2-1}\displaystyle 
\lim_{t\to-\infty}\left(\frac{X^{p+1}(t)}{p+1}-M\gg^q\frac{X^{q+1}(t)}{q+1}\right)=\Gl.
\ee
Therefore $X(t)$ converges to some  $\gl$ satisfying $\frac{\gl^{p+1}}{p+1}-M\gg^q\frac{\gl^{q+1}}{q+1}=\Gl$. The omega-limit set at $-\infty$ of the trajectory 
$\{(X(t),Y(t))\}_{t\in\BBR_-}$ is the set $\Gg$ of couples $(X_0,Y_0)$ such that there exists a sequence $\{t_n\}$ decreasing to $-\infty$ such that 
$(X(t_n),Y(t_n))\to (X_0,Y_0)$. It is non-empty since the trajectory is bounded, connected and compact. By La Salle's theorem, the function 
$\CE(t)-C_2e^{\frac{\gs t}{p-q}}$ which is monotone decreasing is constant on $\Gg$. This implies $M(\gg X_0-Y_0)(\gg^q X_0^q-Y_0^q)=0$, hence 
$Y_0=\gg X_0$. Because $X(t)\to \gl$ then $X_0=\gl$, hence $Y_0=\gg\gl$ and $Y(t)\to\gg\gl$ when $t\to-\infty$. If $M\gg^q\gl^q\neq \gl^p$, it implies that 
$$Y_t(t)=\gth\gg\gl+e^{-\frac{\gs t}{p-q}}(M\gg^q\gl^q- \gl^p+\ge(t))\quad\text{where }\ge(t)\to 0\text{ as }t\to-\infty.
$$
Hence $Y_t(t)=ce^{-\frac{\gs t}{p-q}}(1+o(1)))$ where $c\neq 0$. Clearly this implies that $Y(t)$ cannot be bounded, contradiction. Therefore $M\gg^q\gl^q-\gl^p$. This implies that 
\bel{S3}\gl\in\left\{0,M^{\frac{1}{p-q}}\gg^{\gg}\right\},
\ee
which ends the proof.\qeda
\blemma{prem} Assume $N\geq 2$, and let $u$ be a positive solution of $(\ref{W1})$ in $B_R$ unbounded near $0$. \smallskip

\nind 1- If $q>\frac{N}{N-1}$ and $u(r)=o\left(r^{-\frac{q}{p(q-1)}}\right)$ near $r=0$, then necessarily $q\leq 2$ and 
\bel{S4}\displaystyle 
\lim_{r\to 0}r^{\frac 1{q-1}}|u'(r)|=\eta_{_M}:=\left(\myfrac{\gk}{M}\right)^{\frac 1{q-1}}.
\ee
Therefore
\bel{S5}\displaystyle 
\lim_{r\to 0}r^\gb u(r)=\xi_{_M}:=\frac 1\gb\left(\myfrac{\gk}{M}\right)^{\frac 1{q-1}}\quad\text{if }\,q<2,
\ee
\bel{S6}\displaystyle 
\lim_{r\to 0}|\ln r|^{-1} u(r)=\frac {N-2}{M}\quad\text{if }\,q=2.
\ee

\nind 2-  If $\frac{2p}{p+1}<q=\frac{N}{N-1}$ and $r^{\gb+\ge}u(r)=r^{N-2+\ge}u(r)$ is bounded for any $\ge>0$, then
\bel{S5-1}\displaystyle 
\lim_{r\to 0}r^{N-2}|\ln r|^{N-1} u(r)=\frac{1} {N-2}\left(\myfrac{N-1}{M}\right)^{N-1}\quad\text{if }\,N\geq 3.
\ee
\bel{S6-1}\displaystyle 
\lim_{r\to 0}r|\ln r||u'(r)|=\lim_{r\to 0}\ln(|\ln r|)u(r)=\myfrac{1}{M}\quad\text{if }\,q=N=2.
\ee

\nind 3-  If $\frac{2p}{p+1}\leq q<\frac{N}{N-1}$ and $r^{N-2}u(r)$ is bounded if $N\geq 3$ or $r^{\ge}u(r)$ is bounded for any $\ge>0$ if $N=2$, then there exists $k>0$ such that
\bel{S7}\BA{lll}\displaystyle 
\displaystyle \lim_{r\to 0}r^{N-2}u(r)=k\qquad\text{if }\,N\geq 3,
\EA
\ee
and 
\bel{S8}\BA{lll}\displaystyle 
\displaystyle \lim_{r\to 0}|\ln r|^{-1}u(r)=k\qquad\text{if }\, N=2.
\EA
\ee

\es
\Proof We first  notice  that if $u$ is unbounded near $0$, then $u'<0$ in a neighborhood of $0$ and we can apply the results of \rlemma{subsup} concerning subsolutions. Moreover, if $u^p(r)=o(|u'(r)|^q)$ when $r\to 0$, then for any $\gd>0$ there exists $r_\gd>0$ such that 
\bel{S9}\displaystyle 
M(1-\gd)|u'|^q\leq -\Gd u\leq M|u'|^q \quad\text{in }B_{r_\gd}\setminus\{0\},
\ee
and we can also use the results of \rlemma{subsup} dealing with supersolutions.\smallskip

\nind 1- Since $u$ is a decreasing subsolution of $\CR^M_qu=0$, $|u'(r)|\geq cr^{-\gb-1}$ near $0$, hence $|u'(r)|^q\geq cr^{-\frac{q}{q-1}}$. By assumption
$u^p(r)=o\left(r^{-\frac{q}{q-1}}\right)$. Then $u^p(r)=o(|u'(r)|^q)$ near $0$, hence $(\ref{S9})$ applies.
It follows by \rlemma{subsup}-(1)-(2) that
\bel{S10}\displaystyle 
 \left(\myfrac{\gk}{M}\right)^{\gb+1}\leq \liminf_{r\to 0}r^{\gb+1}|u'(r)|^q\leq \limsup_{r\to 0}r^{\gb+1}|u'(r)|^q\leq  \left(\myfrac{\gk}{M(1-\gd)}\right)^{\gb+1}.
\ee
Since $\gd>0$ is arbitrary, this implies $(\ref{S4})$. The other estimates $(\ref{S5})$ and $(\ref{S6})$ are obtained by integration.\smallskip

\nind 2- By $(\ref{T4})$, $|u'(r)|^q\geq cr^{q(1-N)}|\ln r|^{q(1-N)}=cr^{-N}|\ln r|^{-N}$.  From the assumptions, for any $\ge>0$,  $u^p(r)\leq c_\ge r^{p(2-N-\ge)}$, then 
$$\myfrac{u^p(r)}{|u'(r)|^q}\leq c_\ge'r^{N-(N-2+\ge)p}|\ln r|^{N}.
$$
Next $\frac{2p}{p+1}<\frac {N}{N-1}$ implies that $N>p(N-2)$. Therefore, we can take $\ge>0$ small enough such that $N-(N-2+\ge)p>0$. This implies that $(\ref{S9})$ holds in $B_{r_\gd}\setminus\{0\}$. Hence we get $(\ref{S5-1})$ and $(\ref{S6-1})$ by integration.\smallskip

\nind 3- Suppose $q<\frac{N}{N-1}$ then $p<\frac{N}{N-2}$ if $N\geq 3$. By \rlemma{subsup}-(1), $r^{N-1}|u'(r)|\geq c>0$ near $0$, hence 
$$\myfrac{u^p(r)}{|u'(r)|^q}\leq \tilde c^{-q}r^{(N-1)q-(N-2)p}.
$$
Since $(N-1)q-(N-2)p\geq (N-1)\frac{2p}{p+1}-(N-2)p=\frac{p}{p+1}\left(N-(N-2)p\right)>0$, we deduce that for any $\gd>0$, $(\ref{S9})$ holds in $B_{r_\gd}\setminus\{0\}$. Then we use \rlemma{subsup} and obtain $(\ref{S7})$ and $(\ref{S8})$ by integration. In the case $N=2$ there holds
$$\myfrac{u^p}{|u'|^q}\leq cr^{q-\ge p}
$$
for any $\ge>0$. Choosing $\ge<\frac qp$, we find again $(\ref{S10})$. 
\qeda\medskip

For obtaining the next result, the key is the introduction of the slope function $S$ which allows to make precise the behaviour of solutions such that $r^\gg u(r)\to 0$.
\blemma{mono} Assume $N\geq 2$, $\frac{2p}{p+1}<q<p$ and $M>0$. If $u$ is a positive solution of $(\ref{W1})$ unbounded near $0$ and such that $r^\gg u(r)\to 0$ when 
$r\to 0$, then $q\leq 2$ and the following trichotomy holds.\smallskip

\nind 1-  If $q>\frac{N}{N-1}$, then $(\ref{S5})$ or $(\ref{S6})$ is satisfied.\smallskip

\nind 2-   If $q=\frac{N}{N-1}$, then $(\ref{S5-1})$ or $(\ref{S6-1})$ is satisfied.\smallskip

\nind 3-   If $q<\frac{N}{N-1}$, then $(\ref{S7})$ or $(\ref{S8})$ is satisfied.
\es 
\Proof By assumption $X(t)\to 0$ as $t\to-\infty$. We recall that $S(t)=\frac{Y(t)}{X(t)}$ satisfies $(\ref{U19})$ hence
\bel{S11}\BA{lll}\displaystyle 
X_t=X(\gg-S)\\[1mm]
\,S_t=S(S+2-N)+\xi^{q-1}(MS^q-X^{p-q}).
\EA\ee
1- We first assume that $S(t)\to 0$ as $t\to-\infty$. Then for any $\ge>0$, there exists $r_\ge>0$ such that $0<-\frac{ru'(r)}{u(r)}\leq\ge$ on $(0,r_\ge]$. Hence 
$r\mapsto r^\ge u(r)$ is increasing. This implies that $r^\ge u(r)$ is bounded near $0$ and thus $q\leq 2$ by \rlemma{prem}. If $\frac{N}{N-1}<q<2$ it would follow from \rlemma {prem} that $(\ref{S5})$ holds, which is not possible. Hence $N\geq 3$, $q=2$ and $(\ref{S6})$ holds. If $q=\frac{N}{N-1}$ and $N\geq 3$, 
$(\ref{S5-1})$ cannot hold; hence $N=q=2$ and $(\ref{S6-1})$ holds. If $q<\frac {N}{N-1}$, $(\ref{S7})$ cannot be satisfied, hence $N=2$ and $(\ref{S8})$ holds.\smallskip

\nind 2- Now we assume that $\displaystyle \liminf_{t\to -\infty}S(t)=m> 0$. Then there exist $t_0>-\infty$ and $m_0\in (0,\infty)$ such that 
$S(t)\geq m_0$ for $t\leq t_0$. Hence $Y^q(t)\geq m_0X^q(t)$ therefore $X^p(t)=X^{p-q}(t)X^q(t)=o(Y^q(t))$ as $t\to-\infty$. This implies 
$u^p(r)=o(|u'(r)|^q$ as $r\to 0$. Then $(\ref{S9})$ holds. Using \rlemma{subsup}-(1)-(2), we have $(\ref{S5})$ or $(\ref{S6})$ if $q>\frac{N}{N-1}$, $(\ref{S5-1})$ or $(\ref{S6-1})$ if $q=\frac{N}{N-1}$ and $(\ref{S7})$ or $(\ref{S8})$ if $q<\frac N{N-1}$.\smallskip

\nind 3- Next we assume that $\displaystyle 0=\liminf_{t\to -\infty}S(t)<\limsup_{t\to -\infty}S(t)=\Gs\in (0,\infty]$. Then there exists a decreasing sequence $\{t_n\}$ converging to $-\infty$ such that $S_t(t_n)=0$ and $S_n:=S(t_n)$, which is a local maximum of $S(t)$, tends to $\Gs$. Put $X_n=X(t_n)$ and $\xi_n=\xi(t_t)$, then
\bel{S12}\xi^{q-1}_n=\myfrac{S_n(N-2-S_n)}{MS_n^q-X_n^{p-q}}=\myfrac{S_n(N-2-S_n)}{MS_n^q(1-\ge_n)}=\myfrac{N-2-S_n}{MS_n^{q-1}(1-\ge_n)},
\ee
with $\ge_n\to 0$. This implies in particular $N>2$ and $S_n<N-2$. Since it holds for all local maximum of $S_n$ we deduce $S<N-2$, which implies 
$u(r)\leq Cr^{2-N}$. If $q<\frac{N}{N-1}$ (resp. $q=\frac{N}{N-1}$) we obtain $(\ref{S7})$ from \rlemma{prem}-(3) (resp. $(\ref{S5-1})$ from \rlemma{prem}-(2)). If $q>\frac{N}{N-1}$ we write  $(\ref{S12})$ under the form
\bel{S12-1}\xi^{q-1}_nMS_n^{q-1}(1-\ge_n)=\eta^{q-1}_nM(1-\ge_n)=N-2-S_n.
\ee
From $(\ref{T1})$, $\eta^{q-1}_n\geq \frac{\gk}{M}(1-\ge'_n)$ for $n$ large enough and $\ge'_n\to 0$ when $n\to\infty$, hence 
\bel{S13}N-2-S_n\geq (1-\ge_n)(1-\ge'_n)\gk\Longrightarrow S_n\leq N-2-\gk+\ge''_n=\gb+\ge''_n.
\ee
This implies that for any $\ge>0$ there exists $n_\ge$ such that $S(t)\leq S_{n_\ge}\leq \gb+\frac{\ge}{2}$ for $t\leq t_{n_\ge}$. Hence $r^{\gb+\ge}u(r)\to 0 $ as $r\to 0$. Since $\frac{q}{p(q-1)}=\gb+\frac{q(p+1)-2p}{p(q-1)}$, it implies that $u(r)=o(r^{-\frac{q}{p(q-1)}})$ as $r\to 0$. Therefore $(\ref{S5})$ and $(\ref{S6})$ hold.
\qeda\medskip

\nind{\it Proof of \rth{T6}}. It follows from \rlemma {conv},  \rlemma {prem} and \rlemma {mono}. {\hspace{10mm}\hfill $\phantom{\square}$}\qeda

\subsubsection{The case $1<q<\frac{2p}{p+1}.$ }
{\it Proof of \rth{T7}}. If $1<q<\frac{2p}{p+1}$ and $p\geq \frac{N}{N-2}$ it is proved in \cite{BVGHV3} that positive solutions of $(\ref{Z1})$ in $B_R\setminus\{0\}$ can be extended as a $C^2$ solution in  $B_R$. Next we suppose that $p< \frac{N}{N-2}$, or $N=1,2$, hence $q<\frac N{N-1}$. We use the change of variable $(\ref{U1})$ and $(x, y)$ satisfies $(\ref{U2})$. It is important to notice that $\gs=(p+1)q-2p$ is negative, therefore the system satisfied by $(x, y)$ is a perturbation at $-\infty$ of the system 
\bel{R3}\BA {lll}
x_t=\ga x-y\\
y_t=-Ky-x^p
\EA\ee
where $K=N-2-\ga$, associated to the Emden-Fowler equation $\CL_pu=0$ by the same change of variable. Since $(x(t),y(t))$ is bounded, the omega-limit set at $-\infty$ of the trajectory $\{(x(t),x(t))\}_{t\in\BBR_-}$ is a non-empty compact connected subset of the set of stationary solutions of $(\ref{U3})$. Therefore
\bel{R4}\displaystyle
\lim_{t\to -\infty}(x(t),y(t))=(\ell,\ga\ell)\quad\text{where }\;\ell\in\{0,x_0\}.
\ee 
If $\ell=x_0$ the result is proved, thus let us assume that $\ell=0$. By \rlemma{subsup}-1-(iii) $r^{N-1}u'(r)$ admits a limit $c\in (0,\infty]$ when $r\to 0$. If $c<\infty$, $(\ref{Z25})$ follows by integration. Thus we are left with the case $c=\infty$. Hence $\displaystyle\lim_{r\to 0} r^{N-2}u(r)=\infty$ if $N\geq 3$, or $\displaystyle\lim_{r\to 0}|\ln r|^{-1}u(r)=\infty$ if $N=2$.
 Therefore, for any $k>0$, $u$ is bounded from below in $B_R\setminus\{0\}$ by the function $v_k$ which satisfies $\CL_pv_k=0$ in $B_R\setminus\{0\}$, $v_k=0$ on $\prt B_R$ and $\displaystyle\lim_{r\to 0} r^{N-2}v_k(r)=k$ if $N\geq 3$, or $\displaystyle\lim_{r\to 0}|\ln r|^{-1}v_k(r)=k$ if $N=2$. Letting $k\to\infty$, $v_k\uparrow v_\infty$, and $\displaystyle\lim_{r\to 0}r^\ga v_\infty(r)=x_0$ by \cite{Vesingsol}. This is a contradiction.{\hspace{10mm}\hfill $\phantom{\square}$}\qeda
  \subsection{Estimates at infinity}
  \subsubsection{ The case $q>\frac{2p}{p+1}.$}
 {\it  Proof of \rth{T8}}. We recall that by \rprop{Oss} and \rprop{Osgrad} in the Appendix all the positive solutions of $(\ref{Z1})$ in $B_R^c$ satisfy 
\bel{Q1}u(r)+r|u'(r)|\leq cr^{-\ga}\quad\text{in }B_{R+1}^c.
\ee
where $c=c(N,p,q)>0$, and by the maximum principle they are decreasing. Since $u$ is continuous in $B_R^c$, $\gu=\min\{u(r):r=R\}$ is well defined and positive. By the maximum principle, for any $n>R$, $u$ is bounded from below in $B_n\setminus B_R$ by the solution $\tilde v_n$ of 
\bel{Q2}\CL_p\tilde v_n=0\quad\text{in }B_n\setminus \overline B_R\,,\; \tilde v_n=\gu\quad\text{on }\prt B_R\,,\; \tilde v_n=0\quad\text{on }\prt B_n.
\ee
When $n\to\infty$, $\tilde v_n\uparrow\tilde v_\infty$ which satisfies $\CL_p\tilde v_\infty=0$ in $B_R^c$ and $\tilde v_\infty=\gu$ on $\prt B_R$. Then $u\geq \tilde v_\infty$ and by \cite{Veasym}, $\tilde v_\infty$ satisfies
\bel{Q3}\BA{lll}
(i)\displaystyle\;\;\lim_{r\to\infty}r^\ga \tilde v_\infty(r)=x_0\;&\text{if }1<p<\frac{N}{N-2}\\[2mm]
(ii)\displaystyle\;\;\lim_{r\to\infty}r^{N-2}\left(\ln r\right)^\frac{N-2}{2} \tilde v_\infty(r)=\left(\myfrac{N-2}{\sqrt 2}\right)^{N-2}\;&\text{if }N\geq 3\text{ and }p=\frac{N}{N-2}\\[3mm]
(ii)\displaystyle\;\;\lim_{r\to\infty}r^{N-2}\tilde v_\infty(r)=c>0\;&\text{if }N\geq 3\text{ and }p>\frac{N}{N-2}.
\EA
\ee
We make the change of variable $(\ref{U1})$ and obtain the system $(\ref{U2})$ satisfied by the functions $t\mapsto (x(t), y(t))$. Since $q>\frac{2p}{p+1}$, $\gs$ is positive. Hence the omega-limit set  of the trajectory of $\{(x(t), y(t))\}_{t\geq 0}$ as $t\to\infty$ is a non-empty compact connected set of the set of solutions of stationary solutions of 
$(\ref{R3})$, therefore 
\bel{Q4}\displaystyle
\lim_{t\to \infty}(x(t),y(t))=(\ell,\ga\ell)\quad\text{where }\;\ell\in\{0,x_0\}.
\ee 
Therefore if $1<p<\frac{N}{N-2}$ we obtain $(\ref{Z26})$, and if $p\geq \frac{N}{N-2}$ we have that $\ell=0$. \\
If $p>\!\frac{N}{N-2}$, then $q>\!\frac{N}{N-1}$. From \rlemma{subsup}-(3), we have that either $q<2$ and $(\ref{T7})$ holds, or $N>2$ and $(\ref{T8})$ holds. However, since $q>\frac{2p}{p+1}$, one has $r^{-\gb}=o(r^{-\ga})$ when $r\to\infty$, hence $(\ref{T7})$ does not hold and we deduce that $(\ref{T8})$ is verified. \\
Finally we consider the case $p=\frac{N}{N-2}$. Then $\ell=0$ and $x$ satisfies
\bel{Q5}\displaystyle
x_{tt}-(N-2)x_t-x^{\frac {N}{N-2}}+Me^{-((N-1)q-N)t}\left((N-2)x-x_t\right)^q=0 \quad\text{on }(\ln R,\infty),
\ee 
and $q>\frac{2p}{p+1}=\frac{N}{N-1}$. Since $u$ is bounded from below by $v_\infty$, we have that $x(t)\geq ct^{-\frac{N-2}{2}}$, with $c>0$, for $t$ large enough. Hence for any $\ge>0$ there exists $t_\ge>\ln R$ such that 
\bel{Q6}\displaystyle
x_{tt}-(N-2)x_t-x^{\frac {N}{N-2}}\leq 0\leq x_{tt}-(N-2)x_t-(1-\ge)x^{\frac {N}{N-2}} \quad\text{on }(\ln R,\infty).
\ee 
Therefore $\gth_{2,\ge}(t)\leq x(t)\leq \gth_{1,\ge}(t)$ where 
\bel{Q7}\BA{lll}\displaystyle
\frac{d^2}{dt^2}\gth_{j,\ge}-(N-2)\frac{d}{dt}\gth_{j,\ge}-(1+(-1)^j\ge)\gth_{j,\ge}^{\frac {N}{N-2}}=0 \quad\text{on }(t_\ge,\infty)\\
\gth_{j,\ge}(t_\ge)=x(t_\ge).
\EA\ee 
The asymptotic expansion of $\gth_{j,\ge}$ is obtained in \cite[Lemme 3.2 ]{Veasym} using an old result due to Hardy. We give below a simpler proof. 
\bel{Q8}\BA{lll}\displaystyle
\gth_{j,\ge}(t)=\left(\myfrac{N-2}{\sqrt{1+(-1)^j\ge}}\right)^{N-2}\left(\myfrac{1}{2t}\right)^{\frac{N-2}{2}}(1+o(1)).
\EA\ee
This implies that for any $\ge>0$ there holds
\bel{Q9}\BA{lll}\displaystyle
\left(\!\myfrac{N-2}{\sqrt{2(1+\ge)}}\!\right)^{N-2}\!\!\!\!\!\leq\liminf_{t\to-\infty}t^{\frac{N-2}{2}}x(t)\leq\limsup_{t\to-\infty}t^{\frac{N-2}{2}}x(t)\leq \left(\!\myfrac{N-2}{\sqrt{2(1-\ge)}}\!\right)^{N-2},\!\!\!\!
\EA\ee
which implies $(\ref{Z28})$.\qeda\medskip

\nind\Remark The proof of  Hardy's theorem  quoted in \cite{Bel} is not easy to find. An alternative proof is to consider the following equation, to which  $(\ref{Q7})$ reduces by a suitable scaling transformation,
\bel{Q10}\BA{lll}\displaystyle
\gth''-\gth'-\gth^n=0\quad\text{on }[0,\infty),
\EA\ee
where $n>1$ and $\gth>0$. Since $\gth(t)\to 0$ as $t\to\infty$, it is easy to see that for any $t>1$, $\gth(t)\leq Ct^{-\frac{1}{n-1}}$ by considering supersolutions under the form
$$\psi(t)=at^{-\frac{1}{n-1}}+bt^{-\frac{2}{n-1}}.
$$
Since $\phi(t)=\left(\frac{1}{(n-1)(t+t_0)}\right)^{\frac{1}{n-1}}$ is a subsolution for some $t_0>0$, it is smaller than $\gth(t)$. Furthermore, for any $\ge>0$, there exists $t_\ge>0$ such that
$t\mapsto\left(\frac{(1+\ge)}{(n-1)t}\right)^{\frac{1}{n-1}}$ is a supersolution on $[t_\ge,\infty)$ and is larger than $\gth$. From that we infer
\bel{Q11}\BA{lll}\displaystyle
\lim_{t\to\infty}t^{\frac{1}{n-1}}\gth(t)=\left(\frac{1}{n-1}\right)^{\frac{1}{n-1}}.
\EA\ee
An alternative proof of the convergence is to set $x(t)=t^{\frac{2-N}{2}}\CZ(t)$. We get
$$\CZ_{tt}-(N-2)\left(1+\frac 1t\right)\CZ_{t}+\frac 1t\left(\left(\myfrac{(N-2)^2}{2}+\myfrac{N(N-2)}{4t}\right)\CZ-\CZ^{\frac{N}{N-2}}\right)+\Gf=0
$$
where $\Gf(t)=O\left(t^{\frac{N}{N-2}}e^{(N-(N-1)q)t}\right)$. Applying \cite [Corollary 4.2]{BiRa} we deduce that $\CZ(t)$ converges to a limit $\ell$ which satisfies $\ell\left(\frac{(N-2)^2}{2}-\ell^{\frac{2}{N-2}}\right)=0$. From the lower bound 
$u\geq\tilde v_\infty$ and $(\ref{Q3})$ we infer that $\ell=0$ is impossible.

  \subsubsection{ The case $1<q<\frac{2p}{p+1}.$ }
 {\it Proof of \rth{T9}}.  If $1<q<\frac{2p}{p+1}$, then $q<p$. Therefore $\gg<\ga<\gb$. Hence if $u$ a nonnegative solution of  $(\ref{W1})$ on $[R,\infty)$, $r^\gg u(r)$ is bounded for $r\geq R>0$. Therefore the natural system for describing the solution is the system $(\ref{U9})$ with bounded $X(t)$ and $Y(t)$ and we use an argument similar to the poof of \rlemma{conv}.
\blemma{convbis} Assume $p>1$, $1<q<\frac{2p}{p+1}$ and $M>0$. If $u$ is a positive radial solution of  $(\ref{Z1})$ in $B_R^c$, there holds
  \bel{P1}\BA{lll}\displaystyle
\lim_{r\to\infty}r^\gg u(r)=\ell\in \{0,X_{_M}\}.
\EA\ee
\es
\Proof By \rprop {Osgrad} $Y$ is bounded, hence $(\ref{S2})$ still holds on $[T,\infty)$. We consider now the function $\CE$ defined in $(\ref{U16})$, then $t\mapsto \CE(t)-\frac{C_2|\gs|}{p-q}e^{\frac{\gs t}{p-1}}$ is decreasing and bounded at infinity since $\gs<0$. Therefore $\CE(t)$ converges to some real number $\Gl$ when $t\to\infty$. This implies that identity $(\ref{S2-1})$ is still valid provided 
$\displaystyle\lim_{t\to-\infty}$ is replaced by $\displaystyle\lim_{t\to\infty}$. Mutatis mutandis, the remaining of the proof of \rlemma{conv} still holds and we get 
$(\ref{P1})$.\qeda
\blemma{monobis}Let the structural assumptions of \rlemma{convbis} be satisfied. If $u$ is a positive radial solution of  $(\ref{Z1})$ in $B_R^c$, such that 
$r^\gg u(r)\to 0$ when $r\to\infty$, then necessarily $q>\frac{N}{N-1}$ and the following alternative holds:\\
\nind 1-  either $q<2$ and 
  \bel{P2}\BA{lll}\displaystyle
\lim_{r\to\infty}r^\gb u(r)=\xi_{_M},
\EA\ee
where we recall that $\xi_{_M}$ is defined in $(\ref{Z20})$,\smallskip

\nind 2-  or $N>2$ and
  \bel{P3}\BA{lll}\displaystyle
\lim_{r\to\infty}r^{N-2} u(r)=k>0.
\EA\ee
\es
\Proof Since $\CR^{^M}_qu\leq 0$, we can apply \rlemma{subsup}-(3) provided $q>\frac{N}{N-1}$. If this holds the following estimate from below of $u$ holds: \\
\nind either $q<2$ and $\displaystyle\liminf_{r\to\infty}r^\gb u(r)\geq \xi_{_M}$,\\
\nind  or $N>2$ and $\displaystyle\lim_{r\to\infty}r^{N-2} u(r)=k\geq 0$.\smallskip

\nind 1- We first prove that $r^\gb u(r)$ is bounded and we recall that $S(t)=\frac {Y(t)}{X(t)}$ denotes the slope function. \\
1-(i) If $S(t)\to 0$ as $t\to\infty$, then for any $\ge>0$, $r\mapsto r^\ge u(r)$ is nondecreasing. Hence $u(r)\geq c_\ge r^{-\ge}$ for $r\geq R$, for some 
$c_\ge>0$. This contradicts \rprop{Oss}\\
\nind 1-(ii) If $\displaystyle\liminf_{t\to\infty} S(t)=m>0$. Then there exists $t_0>\ln R$ and $m_0>0$ such that $S(t)\geq m_0$ on $[t_0,\infty)$. Hence 
$Y^q(t)\geq m_0^qX^q(t)$ for $t\geq t_0$ and $u^p=o(|u'|^q)$ as $r\to\infty$. Using \rlemma{subsup}-(3)-(4) we infer that $q>\frac{N}{N-1}$ and $(\ref{P2})$ or $(\ref{P3})$ holds, and in both cases $r^\gb u(r)$ is bounded.\\
\nind 1-(iii) If $S$ satisfies $\displaystyle 0=\liminf_{t\to\infty}S(t)<\limsup_{t\to\infty}S(t)=\Gs\in (0,\infty]$. There exists an increasing sequence $\{t_n\}$ tending
to infinity of local maximum of $S(t)$. As in the proof of \rlemma{mono}-(3) we obtain that $N>2$ and $u(r)\leq Cr^{2-N}$. \\
If $q\geq \frac N{N-1}$, equivalently $\gb\leq N-2$, then $u(r)\leq Cr^{2-N}\leq Cr^{-\gb}$ for $r\geq 1$.\\ 
If $q< \frac N{N-1}$, we have from $(\ref{U5})$ and $\eta=\xi S$ that 
$\xi_t=\xi(\gb-S)>\xi(\gb-N+2)\geq 0$. Hence $\xi(t)$ is increasing with limit $\xi_\infty\leq \infty$. Since at the points $t_n$ of local maximum of 
$S(t)$, we also have $(\ref{S12})$, we obtain the implication
  \bel{P4}\xi_n^{q-1} S_n^{q-1}=\myfrac{N-2-S_n}{M(1-\ge'_n)}\Longrightarrow \xi_\infty^{q-1} \Gs^{q-1}=\myfrac{N-2-\Gs}{M}.
\ee
Hence $\xi_\infty$ is finite, which implies again that $r^\gb u(r)$ is bounded. \smallskip

\nind 2- Convergence. Since $r^\gb u(r)$ is bounded, the trajectory $\{(\xi(t),\eta(t))_{t\geq\ln R}\}$ endows this property, and since $\gs<0$, its omega-limit set at infinity is non-empty, compact, connected and it is a subset of the nonnegative stationary solutions of $(\ref{U5})$.\\
 If $q\leq \frac {N}{N-1}$ the set is reduced to $(0,0)$. Since $\gk\leq 0$, we deduce from $(\ref{U6})$ that $\xi(t)$ is monotone decreasing. It follows from $(\ref{U5})$ that $S(t)>\gb$, hence $u^p=o(|u'|^q)$ as in 1-(ii) and by \rlemma{subsup}-(3)-(4) necessarily $q>\frac {N}{N-1}$, contradiction. \\
If $q> \frac {N}{N-1}$, then either $(\xi(t),\eta(t))$ converges to $(\xi_{_M},\eta_{_M})$ or it converges to $(0,0)$, in   which case $\displaystyle\lim_{r\to\infty}r^{N-2}u(r)=k\geq 0$ by \rlemma{subsup}. The function $u$ is bounded from below in $B_R^c$ by the solution $v$ of 
$$\BA {lll}-\Gd v+v^p=0\qquad&\text{in }\, B_R^c\\
\phantom{-\Gd +v^p}
\displaystyle v=\min_{|x|=R}u(x)&\text{on }\, \prt B_R
\EA$$ 
Since $\displaystyle\lim_{r\to\infty}r^{N-2}v(r)=k'> 0$ and $k\geq k'>0$, 
this ends the proof.
\qeda \medskip

\nind{\it Proof of \rth{T9}}. In all the cases, the basic convergence $(\ref{Z29})$ holds true from \rlemma{convbis}. If the limit of 
$r^\gg u(r)$ is zero, then necessarily $\frac{N}{N-1}<q<2$ and we have $(\ref{Z30})$ or $(\ref{Z31})$.\qeda
\subsection{Solutions of eikonal type}
In order to study the properties of solutions of eikonal type we first give some asymptotic expansion results.
\blemma{devlim} Let $M>0$,  $\frac{2p}{p+1}<q<p$ (resp. $1<q<\frac{2p}{p+1}$) and $\gth\neq 0$ (see $(\ref{Z18'})$ for the definition of $\gth$).
If $(X,Y)$ is a solution of $(\ref{U9})$ which converges to $(X_{_M},Y_{_M})$ when $t\to-\infty$ (resp. $t\to\infty$), then 
$t\mapsto X(t)-X_{_M}$ has a constant sign for $|t|$ large enough.
Furthermore
  \bel{O1}
  X(t)=X_{_M}+\myfrac{\gth\gg^{1-q}X^{2-q}_{_M}}{p(q-1)M}e^\frac{\gs t}{p-q}+O\left(e^\frac{2\gs t}{p-q}\right)
  \text{ as }t\to-\infty\;\; (\text{resp. }t\to\infty).
\ee
Equivalently, with $u(r)=r^{-\frac{q}{p-q}}X(t)$, 
  \bel{O2}
  u(r)=X_{_M}r^{-\frac{q}{p-q}}+\myfrac{\gth\gg^{1-q}X^{2-q}_{_M}}{p(q-1)M}r^{-\frac{p(2-q)}{p-q}}+
  O\left(r^{-\frac{(2p+1)q-2p}{p-q}}\right)
  \text{ as }r\to 0\;\; (\text{resp. }r\to\infty).
\ee
\es
\Proof  (i) {\it Expansion of $MY^q-X^p$}.  Set
$$\Psi(t)=Y_t-\gth Y=e^{-\frac{\gs t}{p-q}}(MY^q-X^p).
$$
Then
$$\Psi_t=e^{-\frac{\gs t}{p-q}}\left(-\myfrac{\gs}{p-q}(MY^q-X^p)+MqY^{q-1}(\gth Y+\Psi)-pX^{p-1}(\gg X-Y)\right).
$$
If $\Psi$ is not monotone, one has at the local extremum $t_n$ of $\Psi$, denoting $\Psi_n=\Psi(t_n)$, $X_n=X(t_n)$ and $Y_n=Y(t_n)$,
$$MqY_n^{q-1}\Psi_n=-Mq\gth Y_n^q+pX_n^{p-1}(\gg X_n-Y_n) +\myfrac{\gs}{p-q}(MY_n^q-X_n^p).
$$
But $\displaystyle\lim_{t\to-\infty}(X(t),Y(t))= (X_{_M},Y_{_M})$, then $pX_n^{p-1}(\gg X_n-Y_n) +\myfrac{\gs}{p-q}(MY_n^q-X_n^p)\to 0$ when $t_n\to-\infty$. Therefore $\Psi_n\to-\gth\gg X_{_M}$. Since the limit is valid for local minima or maxima it follows that $\displaystyle\lim_{t\to-\infty}\Psi(t)=-\gth\gg X_{_M}$.\\
If $\Psi$ is monotone, then $Y_t(t)-\gth Y(t)$ admits a limit when $t\to-\infty$. Since $Y(t)\to Y_{_M}$, it follows that $Y_t(t) $ has a also a limit at $-\infty$ and the only possible one is $0$. Hence $\Psi(t)\to -\gth \gg X_{_M}$. In both case it yields, since $\gth\neq 0$,
 \bel{O4}
MY(t)^q-X(t)^p=-\gth \gg X_{_M}e^{\frac{\gs t}{p-q}}(1+o(1))\quad\text{ as }t\to-\infty.
\ee

\nind (ii) {\it We claim that $X-X_{_M}$ has a constant sign}. If $X$ is nondecreasing (resp. nonincreasing) then $X(t)\geq X_{_M}$ (resp. $X(t)\leq X_{_M}$) for $t\leq 0$. Actually the inequality is strict, otherwise, if there is some $t_0$ such that $X(t_0)= X_{_M}$, we would have $X_t(t)=0$ for $t\leq t_0$ and  $X(t)=X_{_M}$ for $t\leq t_0$. If $\gth \neq 0$ this contradicts the fact that $M\gg^qX^q_{_M}-X^p_{_M}=0$. If $\gth> 0$ we deduce from $(\ref{U10})$ that if  $X(t_n)$ is a local minimum we have
$$e^{-\frac{\gs t_n}{p-q}}\left(M\gg ^qX^q(t_n)-X^p(t_n)\right)=-\gth X(t_n)-X_{tt}(t_n) < 0\Longrightarrow X(t_n)> X_{_M}.$$
This implies that $X(t)>X_{_M}$. Similarly, if $\gth<0$ we get $X(t)<X_{_M}$. \smallskip

\nind (iii) {\it Asymptotic expansion}. We write $X=X_{_M}(1+w)$ and $Y=\gg X_{_M}(1+z)$. Then 
 \bel{O5}\BA {lll}
w_t=\gg(w-z)\\
\,z_t=\gth (1+z)+M\gg^{q-1}X^{q-1}_{_M}e^{-\frac{\gs t}{p-q}}\left((1+z)^q-(1+w)^p\right).
\EA\ee
There holds
$$(1+z)^q-(1+w)^p=qz-pw+\phi(z)-\psi(w),
$$
where $\phi(z)=(1+z)^q-1-qz$ and $\psi(w)=(1+w)^p-1-pw$, therefore $0\leq \phi(z)\leq cz^2$ and $0\leq \psi(w)\leq cw^2$ for $t\leq t^*$. Next, from 
$(\ref{O5})$, 
$$qz+\phi(z)=pw+\psi(w)+a_{_M}e^\frac{\gs t}{p-q}(z_t-\gth(1+z)),
$$
where $a_{_M}=M^{-1}(\gg X_{_M})^{1-q}$, and $z(t)\to 0$ and $z_t(t)\to 0$ when $t\to-\infty$. Therefore the previous identity becomes
 \bel{O6}qz+\phi(z)=pw+\psi(w)-\gth a_{_M}e^\frac{\gs t}{p-q}(1+h(t)),
\ee
where $h(t)\to 0$ when $t\to-\infty$. Next
$$\BA {lll}\,qz\leq qz+\phi(z)\leq\left\{\BA {lll}
qz(1+\ge)&\text{ if } z>0\\
qz(1-\ge)&\text{ if } z\leq 0
\EA\right.:=qz(1+\ge_z)\\[4mm]
pw\leq pw+\psi(w)\leq pw(1+\ge_w),
\EA$$
where $\ge_z=\ge sign\, (z)$ and $\ge_w=w sign\, (w)$. It follows from $(\ref{O6})$ that
 \bel{O7}\BA {lll}
 qz-pw(1+\ge_w)+\gth a_{_M}e^\frac{\gs t}{p-q}(1+h(t))\leq 0\\
 \;qz(1+\ge_z)-pw+\gth a_{_M}e^\frac{\gs t}{p-q}(1+h(t))\geq 0.
\EA\ee
This leads to the following two inequalities verified by $w_t$
$$\BA {lll}
w_t=\gg (w-z)=\myfrac{1}{p-q}(qw-qz)\geq \myfrac{1}{p-q}\left(qw-pw(1+\ge_w)+\gth a_{_M}e^\frac{\gs t}{p-q}(1+h)\right)\\[4mm]
\phantom{w_t=\gg (w-z)=\myfrac{1}{p-q}(qw-qz)}
\geq-w\left(1+\myfrac{p}{p-q}\ge_w\right)+\myfrac{\gth a_{_M}}{p-q}e^\frac{\gs t}{p-q}(1+h),
\EA$$
and
$$\BA {lll}
w_t\leq \myfrac{q}{p-q}\left(w+\myfrac{1}{q(1+\ge_z)}(-pw+ \gth a_{_M}e^\frac{\gs t}{p-q}(1+h)\right)\\[4mm]
\phantom{w_t}
\leq w\left(-1+\myfrac{p}{p-q}\ge_z\right)+\myfrac{\gth a_{_M}}{(p-q)(1+\ge_z)}e^\frac{\gs t}{p-q}(1+h),
\EA$$
and we know from (i) that $w(t)$ keeps a constant sign when $t\to-\infty$. We deduce from the above inequalities that if $\gth<0$ the function 
$t\mapsto e^{(1-\ge)t}w(t)$ is decreasing for some $\ge>0$ and tends to $0$, hence it is negative, while, if $\gth>0$, $t\mapsto e^{(1+\ge)t}w(t)$ is increasing for another $\ge>0$ and tends to $0$, hence it is positive. Then, we can summarize as follows, with a new $\ge>0$
 \bel{O8}\BA {lll}
(1-\ge)\left(\gth ^{-1}w-\myfrac{a_{_M}}{p-q}e^\frac{\gs t}{p-q}\right)\leq-(\gth ^{-1}w)_t\leq  (1+\ge)\left(\gth ^{-1}w-\myfrac{a_{_M}}{p-q}e^\frac{\gs t}{p-q}\right).
\EA\ee
As we have $1\pm\ge+\frac{\gs}{p-q}=\frac{p(q-1)}{p-q}\pm\ge$, the function $t\mapsto e^{(1+\ge)t}\left(\gth^{-1}w-\frac{a_{_M}}{p(q-1)+\ge}e^{\frac{\gs t}{p-q}}\right)$ is increasing and tend to $0$ as $t\to-\infty$. Hence it is positive. In the same way, the function $t\mapsto e^{(1-\ge)t}\left(\gth^{-1}w-\frac{a_{_M}}{p(q-1)-\ge}e^{\frac{\gs t}{p-q}}\right)$ is decreasing, tends to $0$ hence it is negative. Therefore we infer that
 \bel{O9}\BA {lll}
w(t)=\myfrac{\gth a_{_M}}{p(q-1)}e^{\frac{\gs t}{p-q}}(1+o(1)).
\EA\ee
This implies $\psi(w)=O\left(e^{\frac{2\gs t}{p-q}}\right)$. From $(\ref{O6})$, $z=O\left(e^{\frac{\gs t}{p-q}}\right)$, then $\psi(z)=O\left(e^{\frac{2\gs t}{p-q}}\right)$ and $\ge_z=O\left(e^{\frac{\gs t}{p-q}}\right)$. Since $X-X_{_M}=X_{_M}w=\gg Y_{_M}w$, we deduce $(\ref{O1})$. Notice also that from $(\ref{O6})$ there holds
$$z=\myfrac{(2-q)\gth a_{_M}}{q(q-1)}e^{\frac{\gs t}{p-q}}+o\left(e^{\frac{\gs t}{p-q}}\right),
$$
hence
 \bel{O10}\BA {lll}
w_t=\gg(w(t)-z(t))=\myfrac{\gs\gth a_{_M}}{p(p-q)(q-1)}e^{\frac{\gs t}{p-q}}(1+o(1)).
\EA\ee
In particular $X_t$ has the sign of $\gth$, and therefore $X$ is monotone.\qeda\medskip

\nind\Remark If $q=2$ we obtain
 \bel{O11}\BA {lll}
X(t)=\left(M\gg^2\right)^{\frac{1}{p-2}}+\myfrac{2(N-1)-(N-2)p}{2pM}e^{\frac{2t}{p-2}}+O\left(e^{\frac{4t}{p-2}}\right)\quad\text{as }t\to-\infty,
\EA\ee
so we recover the result of \cite{RiVe}.
\subsection{Local or global existence results}
\subsubsection{The systems of order 3}
Since $q\neq \frac{2p}{p+1}$, we can perform the transformation $T_\ell$ and assume that $M=1$. For proving the existence of solutions to $(\ref{Z1})$ there are essentially three methods: the methods of sub and super solutions which has already been developed in Section 2.3, the method of fixed points, and the use of a specific autonomous system of order 3. This last method appears to be entirely new and we explain it below. This system uses the variables $(X,\xi,S)\in \BBR_+\ti\BBR_+\ti\BBR_+$,
 \bel{N1}\BA {lll}
X_t=X(\gg-S)\\
\;\xi_t=\xi(\gb-S)\\
\;S_t=S(S+2-N)+\xi^{q-1}(MS^q-X^{p-q}).
\EA\ee
\blemma{syst3} Let $1<q<p$ with $q\neq\frac{2p}{p+1}$ and $M>0$. If $u$ is a decreasing positive solution of  $(\ref{W1})$, then 
  \bel{N2}\BA {lll}
(X(t),\xi(t),S(t))=\left(r^\gg u(r),r^\gb u(r),r\myfrac{|u'(r)|}{u(r)}\right)\quad\text{with } t=\ln r,
\EA\ee
satisfies $(\ref{N1})$. Conversely, to each trajectory of $(\ref{N1})$ in $\BBR_+\ti\BBR_+\ti\BBR_+$ corresponds a unique solution of $(\ref{W1})$.
\es
\Proof Let $u$ be a decreasing solution of $(\ref{W1})$. We recall that $(X,Y)$ are solutions of $(\ref{U9})$, $S=\frac YX$ and $\xi(t)=r^\gb u(r)$ with $t=\ln r$. Then $(X,S)$ satisfies the following system which is equivalent to $(\ref{U9})$,
 \bel{N3}\BA {lll}
\!X_t=X(\gg-S)\\
S_t=S(S+2-N)+e^{-\frac{\gs t}{p-q}}X^{q-1}(MS^q-X^{p-q}).
\EA\ee
Using $(\ref{U13})$ we have that $\xi^{q-1}=e^{-\frac{\gs t}{p-q}}X^{q-1}$. Since by computation $\xi_t=\xi(\gb-S)$, we deduce that $(X,\xi,S)$ satisfies $(\ref{N1})$.\\
Conversely, let $(X,\xi,S)\in \BBR_+\ti\BBR_+\ti\BBR_+$ be a solution of $(\ref{N1})$, then
$$\myfrac{\xi_t}{\xi}-\myfrac{X_t}{X}=\gb-\gg=-\frac{\gs}{(p-q)(q-1)}.
$$
Hence $\xi(t)=be^{-\frac{\gs t}{(p-q)(q-1)}}X(t)$ for some $b>0$. If we set $a=-\frac{(p-q)(q-1)}{\gs}\ln b$, we see that 
$$\xi(t)=e^{-\frac{\gs (t+a)}{(p-q)(q-1)}}X(t).
$$
Hence 
\bel{N4}\BA {lll}
\!X_t=X(\gg-S)\\
S_t=S(S+2-N)+e^{-\frac{\gs (t+a)}{p-q}}X^{q-1}(MS^q-X^{p-q}).
\EA\ee
Setting $\gt=t+a$, $X^{(a)}(\gt),S^{(a)}(\gt)=(X(t),S(t))=(X(\gt-a),S(\gt-a))$, then 
\bel{N5}\BA {lll}
\!X^{(a)}_\gt=X^{(a)}(\gg-S^{(a)})\\
S^{(a)}_\gt=S^{(a)}(S^{(a)}+2-N)+e^{-\frac{\gs \gt}{p-q}}(X^{(a)})^{q-1}(M(S^{(a)})^q-(X^{(a)})^{p-q}).
\EA\ee
Then the function $\gr\mapsto u^{(a)}(\gr)=\gr^{-\gg}X^{(a)}(\ln\gr)$ satisfies $(\ref{W1})$. Let $(X,\xi,S)$ and $(\tilde X,\tilde\xi,\tilde S)$ be two solutions of 
$(\ref{N1})$. Then there exist $a,\tilde a$ such that 
$$\xi(t)=e^{-\frac{\gs (t+a)}{(p-q)(q-1)}}X(t)\quad\text{and }\,\tilde\xi(t)=e^{-\frac{\gs (t+\tilde a)}{(p-q)(q-1)}}\tilde X(t),
$$
and 
$$\BA {lll}u^{(a)}(\gr)=\gr^{-\gg}X^{(a)}(\ln\gr)=\gr^{-\gg}X(\ln\gr-a)\\
u^{(\tilde a)}(\gr)=\gr^{-\gg}X(\ln\gr-\tilde a).
\EA
$$
If $(X,\xi,S)$ and $(\tilde X,\tilde\xi,\tilde S)$ correspond to the same trajectory, there exists $h\in\BBR$ such that 
$(\tilde X,\tilde\xi,\tilde S)(t)=(X,\xi,S)(t+h)$ for all $t$, thus
$$\xi(t+h)=e^{-\frac{\gs (t+\tilde a)}{(p-q)(q-1)}}X(t+h)=e^{-\frac{\gs (t+ a+h)}{(p-q)(q-1)}}X(t+h).
$$
Therefore $\tilde a=a+h$. Hence
$$u^{(\tilde a)}(\gr)=\gr^{-\gg}\tilde X(\ln\gr-\tilde a)=\gr^{-\gg}\tilde X(\ln\gr- a-h)=\gr^{-\gg}X(\ln\gr-a)=u^{(a)}(\gr).
$$
In conclusion, there is a one to one correspondence between the trajectories of $(\ref{N1})$ and the solutions of $(\ref{W1})$.\qeda\medskip

\nind\Remark Using the relation $(\ref{U14r})$ one can see that $(\ref{N1})$ is equivalent to the  following system in the variables $(x,\xi,S)$, 
 \bel{N6}\BA {lll}
x_t=x(\ga-S)\\
\;\xi_t=\xi(\gb-S)\\
\;S_t=S(S+2-N)+\xi^{q-1}MS^q-x^{p-1}.
\EA\ee
This system is particularly suitable for constructing local solutions in $r^{-\ga}$ or $r^{2-N}$, in particular when $r\to\infty$, in the case $q>\frac{2p}{p+1}$.


\subsubsection{Singular solutions of eikonal type}
{\it Proof of \rth{T10}.} We recall that these solutions of eikonal type are the solutions which behave like $r^{-\gg}$ near $0$ or $\infty$. For $c>0$ and $A\geq 0$ we set $u_{_{c,A}}(r)=cr^{-\gg}+A$ and $u_{_{c}}=u_{_{c,0}}$.
Then there exist $a,b>0$ depending on $p$ such that  
$$ c^pr^{-\gg p}+aA^p+apc^{p-1}Ar^{-\gg(p-1)}\leq u^p_{_{c,A}}(r)\leq c^pr^{-\gg p}+bA^p+bpc^{p-1}Ar^{-\gg(p-1)}.
$$
\nind (i) {\it Subsolutions}. 
If $u_{_{c}}(r)=cr^{-\gg}$, then
 \bel{N7}\BA{lll}\myfrac 1c\CL^{M}_{p,q}u_{_{c}}=-\gg(\gg+2-N)r^{-\gg-2}+c^{p-1}r^{-\gg p}-Mc^{q-1}\gg^qr^{-(\gg+1)q}\\[0mm]\phantom{\myfrac 1c\CL^{M}_{p,q}u_c}
=r^{-\frac{pq}{p-q}}\left(c^{p-1}-c^{q-1}X_{_M}^{p-q}-\gg\gth r^{\frac{\gs}{p-q}}\right).
\EA\ee
Set $\Phi(c)=c^{p-1}-c^{q-1}X_{_M}^{p-q}$. Then $\Phi(X_{_M})=0$ and $\Phi$ achieves its minimum at $c_m=\left(\frac{q-1}{p-1}X_{_M}^{p-q}\right)^{\frac 1{p-q}}$ with minimal value $\Phi(c_m)=-\frac{p-q}{q-1}\left(\frac{q-1}{p-1}X_{_M}^{p-q}\right)^{\frac {p-1}{p-q}}$. Notice that if $\gth=0$ i.e. $q=\frac{(N-2)p}{N-1}$, then 
we find the explicit solution $u^*_{_{M,p}}$ defined in $(\ref{U11})$. \\
(i-a) If $N=1$ or $N=2$, or if $N\geq 3$ and $q>\frac{N-2}{N-1}p$, then $\gth>0$ and $u_{_{c}}$ is a subsolution in $\BBR^N\setminus\{0\}$ provided $c\leq X_{_M}$.\\
(i-b) If $q>\frac{2p}{p+1}$ and $\gth<0$ there exists $r_1>0$ such that $u_{_{c_m}}$ is a subsolution in $B_{r_1}\setminus\{0\}$ Hence 
$\tilde u_{_{c_m}}=c_m(r^{-\gg}-r_1^{-\gg})_+$ is a subsolution in $\BBR^N\setminus\{0\}$.\smallskip

\nind (ii) {\it Supersolutions}. We have 
 \bel{N8}\BA{lll}\myfrac 1c\CL^{M}_{p,q}u_{_{c,A}}\geq r^{-\frac{pq}{p-q}}\left(\Phi(c)-\gg\gth r^{\frac{\gs}{p-q}}\right)+\myfrac{aA^p}{c}+apc^{p-2}Ar^{-\gg(p-1)}
\EA\ee
(ii-a) If $\gth<0$, then for $c\geq X_{_M}$ and any $A\geq 0$, $u_{_{c,A}}$ is a supersolution in $\BBR^N\setminus\{0\}$.\\
(ii-b) If $q>\frac{2p}{p+1}$ and $\gth>0$, then for any $R>0$ we take $c>X_{_M}$ such that $\Phi(c)\geq \gg\gth R^\frac{\gs}{p-q}$, hence $\CL^{M}_{p,q}u_{_{c,A}}\geq 0$ in $B_R\setminus\{0\}$. Since  $-\frac{pq}{p-q}+\frac{\gs}{p-q}=\frac{q-2p}{p-q}<0$, we take $A>0$ such that $aA^p\geq\gg\gth cR^{-\gg-2}$, hence $\CL^{M}_{p,q}u_{_{c,A}}\geq 0$ in $B_R^c$.  Consequently  $u_{_{c,A}}$ is a supersolution in $\BBR^N\setminus\{0\}$.\\
(ii-c) If $q<\frac{2p}{p+1}$ and $\gth>0$, then we can take $c$ such that $\Phi(c)\geq\gg\gth R^\frac{\gs}{p-q}$ and obtain that $u_{_{c,A}}$ is a supersolution in $B_R^c$. \smallskip

\nind (iii) {\it Proof of statements 1 and 2}.\\
If $q>\frac{2p}{p+1}$ and whatever is the sign of $\gth$ there exist $c_m\leq c<X_{_M}<c'$ and $A>0$ such that $u_{_{c',A}}$ is a supersolution in $\BBR^N\setminus\{0\}$ larger than the subsolution $u_{_{c}}$. By \cite[Theorem 1.4.5]{Vebook} there exists a radial solution $u$ in $\BBR^N\setminus\{0\}$ satisfying $u_{_{c}}\leq u\leq u_{_{c',A}}$. Its behaviour at infinity is given by \rth{T8}. This solution is decreasing by the maximum principle and it is unique by \rth{uni}-(3). \\
The existence of a solution in a bounded domain $\Gw$ containing $0$ and vanishing on $\prt\Gw$ satisfying $(\ref{Z19})$ follows by \rth{Exist1} which is proved in Section 4. So we deduce statement 1.\\
\nind If $1<q<\frac{2p}{p+1}$ and $\gth>0$, one has a supersolution $u_{_{c',A}}$ in $B_R^c$ and a subsolution $u_{_{c}}$ in $\BBR^N\setminus\{0\}$.
Up to increasing the value of $A$ one has again a supersolution $u_{_{c',A}}$ larger than the subsolution $u_{_{c}}$. Hence there exists a solution $u$ in between satisfying $(\ref{Z29})$ which proves statement 2.
\qeda
\subsubsection{Riccati type singular solutions}
{\it Proof of \rth{T11}.} We recall that the Riccati equation $(\ref{Z16})$ admits the radial solution $\xi_{_M}|x|^{-\gb}$ if and only if $\frac N{N-1}<q<2$. This function is a supersolution of $(\ref{Z1})$ in $\BBR^N\setminus\{0\}$.\medskip

\nind 1- {\it Local existence in a neighborhood of $0$}. Since $q>\frac{N}{N-1}$ the system $(\ref{N1})$ in variables $(X,\xi,S)$ admits the equilibria $(0,0,0)$, $(0,0,N-2)$ and $(0,\xi_{_M},\gb)$.
Our aim is to construct local radial solutions of $(\ref{W1})$ satisfying $\displaystyle \lim_{r\to 0}r^\gb u(r)=\xi_{_M}$ and $\displaystyle \lim_{r\to 0}r^{\gb+1} |u'(r)|=\eta_{_M}=\gb\xi_{_M}$, equivalently 
 \bel{M1}\displaystyle \lim_{t\to -\infty}(X(t),\xi(t),S(t))=(0,\xi_{_M},\gb).\ee
Conversely, any solution $(X,\xi,S)$ satisfying $(\ref{M1})$ corresponds to a solution $u$ satisfying $\displaystyle \lim_{r\to 0}(r^\gb u(r),r^{\gb+1} |u'(r)|)=(\xi_{_M},\gb\xi_{_M})$. The system  $(\ref{N1})$ may be singular at $\xi_{_M}=0$ or at $X=0$; hence we desingularize it by setting $\hat X=X^{p-q}$ and $\hat \xi=\xi^{q-1}$. Then $(\hat X,\hat\xi,S)$ satisfies 
 \bel{M2}\BA {lll}
 \hat X_t=(p-q)\hat X(\gg-S)\\
 \;\hat \xi_t=(q-1)\hat\xi(\gb-S)\\
 \;S_t=S(S+2-N)+\hat\xi(MS^q-\hat X).
 \EA\ee
 So we are led to study solutions in a neighborhood of the equilibrium $(0,\hat\xi_{_M},\gb)$ where $\hat\xi_{_M}=\xi^{q-1}_{_M}=\frac{\gk}{M\gb^{q-1}}$. We set
 $\hat\xi=\hat\xi_{_M}+\bar \xi$, $\hat X=\bar X$ and $S=\gb+\bar S$ in order to reduce the study at $(0,0,0)$, and $(\bar \xi,\bar X,\bar S)$ satisfies the following linearized system
  \bel{M3}\BA {lll}
 \bar X_t=\frac{\gs}{q-1}\bar X\\
 \;\bar \xi_t=-(q-1)\hat\xi_{_M}\bar S\\
 \;\bar S_t=-\hat\xi_{_M}\bar X+M\gb^q\bar\xi+(\gb+\gk(q-1))\bar S.
 \EA\ee
 If we denote by $\CA$ the matrix of this system, then its charecteristic values are the roots of the polynomial
   \bel{M4}\BA {lll}
det (\CA-\gm I)=\left(\gm-\gm_1\right)\left(\gm-\gm_2\right)\left(\gm-\gm_1\right),
 \EA\ee
with $\gm_1=\myfrac{\gs}{q-1}$, $\gm_2=\gb$ and $\gm_3=(q-1)\gk=(N-1)q-N$. Since $q>\max\left\{\frac{2p}{p+1},\frac{N}{N-1}\right\}$ all the eigenvalues are positive. We find that 
 $${\bf u}_2=\left(0,1,-\myfrac{M\gb^q}{\gk (q-1)}\right)\;\,\text{and }\;{\bf u}_3=\left(0,1,-M\gb^{q-1}\right)
 $$
 are eigenvectors corresponding to $\gm_2$ and $\gm_3$ respectively. If $\gm_1\neq\gm_2$ and $\gm_1\neq\gm_3$, we can take for eigenvector corresponding to $\gm_1$ the vector ${\bf u}_1=\left(1,b,c\right)$ for some real numbers $b$ and $c$. Actually
 $$b=-\myfrac{\gm_3}{M\gm_2^{q-1}\gm_1}\quad\text{where }\; c=-\myfrac{\gk\gm_1}{M\gm_2^{q-1}(\gm_1-\gm_2)(\gm_1-\gm_3)}.
 $$
 
 Then there exists 
 one trajectory of $(\ref{M2})$  with $X(t)>0$ when $t\to-\infty$ such that $\xi(t)=\xi_{_M}+O(e^{\frac{\gs t}{q-1}})$ when $t\to-\infty$. Hence there exists at least one solution $u$ of $(\ref{W1})$ such that $u(r)=r^{-\gb}\xi_{_M}+Cr^{-\gb+(N-1)q-N}(1+o(1))$ when $r\to 0$. \medskip

\nind 2- {Local existence at infinity}. Here we assume $\frac{N}{N-1}<q<\frac{2p}{p+1}$. Then $\gm_1<0$, $\gm_2=\gb>0$ and $\gm_3=\gk(q-1)>0$. Then there exists a unique local trajectory which converges to $(0,\hat\xi,\gb)$ when $t\to\infty$, it corresponds to the stable manifold of this point. If there exists a positive solution in $\BBR^N\setminus\{0\}$, the  solution can be extended as a solution in $\BBR^N$ by \cite [Theorem 1.1]{BVGHV3} since in this range of values of $q$ one has $p>\frac{N}{N-2}$. By \rprop{Oss} such a solution is identically $0$. \qeda
\medskip

\nind\Remark Note that we have many types of trajectoriess converging to the origin and their geometry depends in their sign and their relative order. In this respect we denote
\bel{M7}f(q):=(N-1)q-N+\frac{N-(N-2)q}{2-q},\ee
and we have
   \bel{M6}\BA {lll}
(i)\qquad &\gm_3=\gk(q-1)>\gm_2=\gb\Longleftrightarrow q>1+\myfrac{1}{\sqrt{N-1}}\\
(ii) \qquad &\gm_1>\gm_2\Longleftrightarrow p<\myfrac{2(q-1)}{2-q}\Longleftrightarrow q>\frac{2(p+1)}{p+2}\\
(iii) \qquad &\gm_1>\gm_3\Longleftrightarrow p<f(q).
 \EA\ee
 We have that $\gm_1=\gm_2=\gm_3$ only if $p=\frac{2}{\sqrt{N-1}-1}$ and $q=1+\frac{1}{\sqrt{N-1}}$, a condition which is compatible with $p>1$ only if $2\leq N\leq 9$. \medskip


 Global (necessarily singular) solutions in $r^{-\gb}$ are difficult to construct. We give below a range of exponents in which there exists at least one. 
 \bth{lamglo} Let $M>0$, $p>1$ and $1<q<2$, $q\neq1+\frac{1}{\sqrt{N-1}}$. If there holds
    \bel{M5_*}\BA {lll}
    p<\max\left\{\myfrac{2(q-1)}{2-q}, f(q)\right\},
     \EA\ee
in particular if $p<\frac{N}{N-2}$ and $q>\frac{N}{N-1}$, then there exists a positive radial solution of $(\ref{Z1})$ defined in $\BBR^N\setminus\{0\}$ satisfying 
    \bel{M5}\BA {lll}\displaystyle
\lim_{x\to 0}|x|^\gb u(x)=\xi_{_M}.
 \EA\ee
 \es
\Proof The function $U(x)=\xi_{_M}|x|^{-\gb}$ is  a supersolution of  $(\ref{Z1})$ in $\BBR^N\setminus\{0\}$. We look for a subsolution under the form
$\tilde\xi(t)=\xi_{_M}(1-Ae^{dt})_+$ for some $d,A>0$. Set
$$H[\tilde\xi](t)=\tilde\xi_{tt}+D\tilde\xi_t-\gk\gb\tilde\xi-e^{\gm_1t}\tilde\xi^p+M|\gb\tilde\xi-\tilde\xi_t|^q,
$$
where 
$$D=\myfrac{Nq-N-2}{q-1}=\gk-\gb.
$$
Then on the interval $I_A:=(-\infty,-\frac{\ln A}{d})$ one has 
$$\gb\tilde\xi-\tilde\xi_t=\xi_{_M}\left(\gb-A(\gb-d)e^{dt}\right).$$
 In order $H[\tilde\xi]\geq 0$, one needs
 $$\BA {lll}
-A\left(d^2+Dd-\gk\gb\right)e^{dt}+M \xi^{q-1}_{_M}\gb^q\left(1-A\myfrac{\gb-d}{\gb}e^{dt}\right)^q-\gk\gb\\[3mm]
\phantom{----------------------}
- \xi^{p-1}_{_M}e^{\gm_1t}\left(1-Ae^{dt}\right)^p\geq 0.
 \EA$$
Since $M\xi_{_M}^{q-1}\gb^q=\gk\gb$, if we set $Z=Ae^{dt}$, then $0<Z\leq 1$ on  $I_A$ and  the previous inequality to be verified becomes
 $$\BA {lll}
\myfrac{\xi^{p-1}_{_M}}{A}e^{(\gm_1-d)t}Z\left(1-Z\right)^p\leq -\left(d^2+Dd-\gk\gb\right)Z+\gk\gb\left(1-\myfrac{\gb-d}{\gb}Z\right)^q-\gk\gb.
 \EA$$
We first impose $d\leq \gm_1$, then $\frac{e^{(\gm_1-d)t}}{A}\leq A^{-\frac{\gm_1}{d}}$ on $ I_A$. We set 
    \bel{M9}\BA {lll}
   Q(Z)=\gk\gb\left(1-\myfrac{\gb-d}{\gb}Z\right)^q-\left(d^2+Dd-\gk\gb+\xi^{p-1}_{_M}
    A^{-\frac{\gm_1}{d}}\right)Z-\gk\gb.
  \EA\ee
  Then
 $$Q'(Z)=-q\gk(\gb-d)\left(1-\myfrac{\gb-d}{\gb}Z\right)^{q-1}-\left(d^2+Dd-\gk\gb+\xi^{p-1}_{_M}
    A^{-\frac{\gm_1}{d}}\right),
 $$
 and 
  $$Q''(Z)=\myfrac{\gk(\gb-d)^2q(q-1)}{\gb}\left(1-\myfrac{\gb-d}{\gb}Z\right)^{q-2}.
 $$
 Since $\gk>0$, $Q$ is convex on $[0,1]$. Furthermore $Q(0)=0$. Hence $H(\tilde \xi)\geq 0$ if $Q'(0)\geq 0$. 
     \bel{M11}\BA {lll}
Q'(0)=-q\gk(\gb-d)-\left(d^2+Dd-\gk\gb+\xi^{p-1}_{_M}
    A^{-\frac{\gm_1}{d}}\right)
 \EA\ee
 Replacing $D$ by its value, $(\ref{M11})$ will be achieved, provided $A$ is large enough, if
      \bel{M12}\BA {lll}
-d^2+((q-1)\gk+\gb)d -(q-1)\gk\gb=-(d-\gm_2)(d-\gm_3) >0.
 \EA\ee
 The condition is that $\gm_1\geq d$ with $d$ satisfying  $(\ref{M12})$. It necessitates $\gm_2\neq\gm_3$, equivalently 
 $q\neq1+\frac{1}{\sqrt{N-1}}$, and\\
 \nind (i) either $\gm_2<\inf\{\gm_1,\gm_3\}$, then we can choose any $d\in (\gm_2,\inf\{\gm_1,\gm_3\})$,\\
 \nind (ii) or $\gm_3<\inf\{\gm_1,\gm_2\}$, then we can choose any $d\in (\gm_2,\inf\{\gm_1,\gm_2\})$.\\
 These conditions are satisfied if $\gm_1>\inf\{\gm_2,\gm_3\}$ which is equivalent to $(\ref{M5_*})$.
 If one of the above conditions is satisfied, it follows by \cite[Corollary 1.4.5]{Vebook}
that there exists a radial positive solution $u$ of $(\ref{Z1})$ in $\BBR^N\setminus\{0\}$ which satisfies
      \bel{M13}\BA {lll}
\xi_{_M}\left(1-A|x|^{d}\right)_+|x|^{-\gb}\leq u(x)\leq \xi_{_M}|x|^{-\gb}\quad\text{for all }\,x\in \BBR^N\setminus\{0\}.
 \EA\ee
 \qeda\medskip

\nind\Remark Condition $(\ref{M5_*})$ is equivalent to 
 \bel{M15}\BA{lll}
 (i)\qquad &1<p<f(q)\qquad &\text{if }\, 1<q<1+\myfrac{1}{\sqrt{N-1}}\qquad \qquad \qquad \\[3mm]
 (ii) \qquad &1<p<\myfrac{2(q-1)}{2-q}\qquad &\text{if }\, 1+\myfrac{1}{\sqrt{N-1}}<q<2.
\EA \ee
Condition (ii) is equivalent to 
\bel{M15'}
\myfrac{2(p+1)}{p+2}<q<2.
\ee
Note that the nature of the variations of the function $p=f(q)$ differs according to the value of $N$. \\
If $N=3$ or $4$, $f$ is increasing and onto from $(\frac 32,2)$ to $(3,\infty)$ when $N=3$ and from $(\frac43,2)$ to $(2,\infty)$ when $N=4$.\\
If $N\geq 5$, $f$ achieves a maximal value $\tilde p$ for $q=\tilde q$ with 
        \bel{M14}\BA {lll}
\tilde q=2-\sqrt{\myfrac{N-4}{N-1}}\, \text{ and }\;\tilde p=2\left(N-2-\sqrt{(N-4)(N-1)}\right).
 \EA\ee
In particular one has
$$f\left(\frac N{N-1}\right)=f\left(\frac N{N-2}\right)=\frac N{N-2}.$$

\subsubsection{Emden-Fowler type singular solutions}
\nind{\it  Proof of \rth{T12}-(1)}. Since $1<p<\frac{N}{N-2}$ the function $x\mapsto U_{x_0}(x)=x_0|x|^{-\ga}$ is a subsolution of $(\ref{Z1})$ in $\BBR^N\setminus\{0\}$. In order 
$x\mapsto U_{_C}(x):=C|x|^{-\ga}$ to be a supersolution, one needs 
      \bel{L1}C^{p-1}\geq x_0^{p-1}+\ga^qC^{q-1}M|x|^{-\frac{\gs}{p-1}}.
\ee
The function $C\mapsto C^{p-q}-x_0^{p-1}C^{1-q}$ is increasing and onto from $[x_0,\infty)$ to $[x_0,\infty)$. Hence there exists $C>x_0$ such that 
$\ga^qC^{q-1}M=C^{p-1}-x_0^{p-1}$. For such a value we have that 
$$ \CL^{M}_{p,q}U_{_C}=M\ga^q(1-r^{-\frac{\gs}{p-1}})r^{-\frac{2p}{p-1}}.
$$
Since $\gs<0$ the function $U_{_C}$ is a supersolution of $(\ref{Z1})$ in $B_1\setminus\{0\}$. For $A>0$ we set 
$U_{_{C,A}}=U_{_C}+A$. Then 
      \bel{L2}\BA {lll}
      \CL^{M}_{p,q}U_{_{C,A}}=\CL^{M}_{p,q}U_{_C}+U^p_{_{C,A}}-U^p_{_C} \geq A^p+M\ga^q(1-|x|^{-\frac{\gs}{p-1}})|x|^{-\frac{2p}{p-1}}.
\EA\ee
Clearly $\CL^{M}_{p,q}U_{_{C,A}}\geq 0$ in $B_1\setminus\{0\}$, and for $|x|>1$, one has
$$ \CL^{M}_{p,q}U_{_{C,A}}\geq A^p-M\ga^q|x|^{-\frac{(p+1)q}{p-1}}\geq A^p-M\ga^q.
$$
Therefore, if $A=M\ga^{\frac qp}$, the function $U_{_{C,A}}$ is a supersolution in $\BBR^N\setminus\{0\}$. Since $U_{x_0}\leq U_{_{C,A}}$, it follows by \cite[Theorem 1.4.5]{Vebook} that there exists a solution $u$ of $(\ref{Z1})$ in $\BBR^N\setminus\{0\}$ such that $U_{x_0}\leq u\leq U_{_{C,A}}$. Then by 
\rth{T7}-(1), u satisfies $(\ref{Z32})$-(i), and by \rth{T9}-(2), $(\ref{Z32})$-(ii) holds. Furthermore 
$r^\ga u(r)\geq x_0$ for any $r>0$. Uniqueness (not only for radial solutions) is a consequence of \rth{uni}-(2). Obviously $|x|\mapsto u(x)$ is decreasing. Existence of a positive solution in a bounded domain $\Gw$ containing $0$ is a consequence of \rth{Exist1}, see Section 4.{\hspace{15mm}\hfill $\phantom{--\square}$}\qeda\smallskip

\nind{\it Proof of \rth{T12}-(2)}. It is a consequence of \rth{T10} and \rth{T8}.{\hspace{10mm}\hfill $\phantom{\square}$}\qeda
    \subsubsection{Solutions behaving like the Newtonian potential}
There exist also solutions which behave like the Newtonian kernel at $0$. They are described in the next result.

 \bth{weaksing} Let $1<p<\frac{N}{N-2}$ and $1<q<\frac{N}{N-1}$. Then for any $M\geq 0$ and $k>0$ there exists a minimal positive solution $u_k$ of $(\ref{Z1})$ in 
 $\BBR^N\setminus\{0\}$ such that $(\ref{Z8})$ holds. Furthermore it is radial and nonincreasing. If we assume $1<q\leq \frac{2p}{p+1}$, this solution is unique among all the positive solutions. 
 \es
 \Proof {\it Proof of existence.} If $M=0$ the result is classical and for $k>0$ we denote by $v_k$ the solution of $\CL_p v=0$ in $\BBR^N\setminus\{0\}$ satisfying $(\ref{Z11})$. This is a natural subsolution of $(\ref{Z1})$. \\
 The construction of the supersolution is more involved. \\
 \nind{(i) We first assume that $N\geq 3$ and prove that for any $k>0$ there exists $M_k>0$ such that for any $0<M\leq M_k$ there exists a supersolution of $(\ref{Z1})$ satisfying $(\ref{Z8})$. Let $a>0$ set
 $$w_k(x)=k|x|^{2-N}+k^q|x|^{2-(N-1)q}+a
 $$
Then there exist $c_5, c_6>0$ depending on $N$ and $q$ such that.
\begin{equation}\label{X5}\BA {lll}\!
\CL_{p,q}^Mw_k=k^q((N\!-\!1)q-\!2)(N-(N-\!1)q))|x|^{(1-N)q}+(k|x|^{2-N}\!+k^q|x|^{2-(N-1)q}\!+a)^p\\[1mm]
\phantom{\CL_{p,q}^Mw}
-M\left((N\!-\!2)k|x|^{1-N}\!+\!((N\!-\!1)q\!-\!2)k^q|x|^{1-(N-1)q}\right)^q\\[1mm]
\phantom{\CL_{p,q}^Mw_k}
\geq c_5k^q|x|^{(1-N)q}+a^p-c_6M\left(k^q|x|^{(1-N)q}+k^{q^2}|x|^{q-(N-1)q^2}\right)\\[1mm]
\phantom{\CL_{p,q}^Mw_k}
\geq k^q\left(c_5-c_6M(1+k^{q^2-q})\right)|x|^{(1-N)q}+a^p-k^q\left(c_6M+k^{q^2-q}\right).
\EA\end{equation}
Note that we have only used inequalities $2\leq (N-1)q\!\leq N$.
Set $M_k=\frac{c_5}{c_6M(1+k^{q^2-q})}$. Then, for $M\leq M_k$ we take $a^p=k^q\left(c_6M+k^{q^2-q}\right)$ and we derive that $\CL_{p,q}^Mw_k\geq 0$.
The supersolution $w_k$ satisfies $v_k\leq k|x|^{2-N}\leq w_k$.\\
\nind (ii) If $N=2$ and for $b>0$ we denote by $\psi_k $ the solution of 
\begin{equation}\label{X6}\BA {lll}\!
-\Gd\psi+\psi^p=|x|^{-q}+2\gp k\gd_0,
\EA\end{equation}
and we set $w_k=\psi_k+b$. Since $1<q<2$, $w_k=-k\ln|x| (1+o(1))$ and $\nabla w_k=-k|x|^{-1}(1+o(1))$  as $x\to 0$. Hence
\begin{equation}\label{X7}\BA{lll}
0\leq \psi_k(x)\leq c_7(-k\ln|x|+1)\;{and }\;\,|\nabla\psi_k(x)|\leq c_7(k+1)|x|^{-1}\quad\text{for } 0<|x|\leq 1.
\EA\end{equation}
Furthermore, by Keller-Osserman technique combined with scaling method, there holds in $\BBR^2\setminus B_1$,
\begin{equation}\label{X8}\BA{lll}
(i)\qquad &0\leq \psi_k(x)\leq c_8\max\left\{|x|^{-\ga}, |x|^{-\frac{q}{p}}\right\},\qquad\qquad\qquad\qquad\qquad\\[3mm]
(ii)\qquad& |\nabla\psi_k(x)|\leq c_8\max\left\{|x|^{-\frac{p+1}{p-1}}, |x|^{-\frac{p+q}{p}}\right\}.
\EA\end{equation}
In the above inequalities, $c_7$ and $c_8$ are positive constants depending on $p$ and $q$. Hence
\begin{equation}\label{X9}\BA {lll}\!
\CL_{p,q}^Mw_k=|x|^{-q}+(\psi_k+b)^p-\psi_k^p-M|\nabla \psi_k|^q.
\EA\end{equation}
We infer
\begin{equation}\label{X10}\BA {lll}\!
\CL_{p,q}^Mw_k\geq |x|^{-q}+b^p-Mc_7^qk^q|x|^{-q}\quad\text{if }0<|x|\leq 1,
\EA\end{equation}
and
\begin{equation}\label{X11}\BA {lll}\!
\CL_{p,q}^Mw_k\geq b^p-Mc_8^qk^q\quad\text{if }x|\geq 1.
\EA\end{equation}
If $k$ is fixed,  $M\leq M_k:=k^{-q}c_7^{-q}$ and $b\geq Mc_8^qk^q$ we conclude that $w_k$ is a supersolution in $\BBR^2\setminus\{0\}$ larger than $v_k$. \smallskip

\nind We deduce from (i) and (ii) that for any $k>0$ there exists $M_k>0$ such that for any $0<M\leq M_k$ there exists a positive radial solution $u_k$ of $(\ref{Z1})$ satisfying $(\ref{Z8})$. Furthermore $u_k$ satisfies $(\ref{Y1})$. Therefore $u_k$ is necessarily decreasing. \smallskip

\nind{\it End of the proof of existence.} Let $q<\frac{N}{N-1}$ and $q_1$ such that $q<q_1<\frac{N}{N-1}$. For  $k>0$ let $\ge>0$ such that for any 
$0<M'\leq \ge$ there exists a positive radial solution $w_k$ to $\CL_{q_1,M'}w=0$ satisfying $(\ref{Z8})$. 
If $M>M'$ there holds 
$$M|X|^q\leq M'|X|^{q_1}+C\qquad\text{for all }X\in\BBR^N,
$$
where $C=\left(\frac{qM}{q_1M'}\right)^{\frac{q}{q_1-q}}\left(M-\frac{q}{q_1}M'\right)>0$. Then
$$\CL_{p,q}^Mw_k=\CL_{q_1,M'}w_k+(M'-M)|\nabla w_k^q|\geq -C
$$
which implies that $w_k+C^{\frac 1p}$ is s supersolution of $(\ref{Z1})$ and $v_k\leq w_k+C^{\frac 1p}$. We conclude as in the first step. \smallskip

\nind{\it Uniqueness.} It is proved in \rth{uni}, this ends the proof. \\
\nind When we do not assume $q\leq \frac{2p}{p+1}$ we have only the existence of a minimal positive solution. This is due to the fact that for two solutions 
$u$ and $u'$ as above, $\min\{u,u'\}$ is a supersolution larger that $v_k$. The conclusion follows easily.\medskip
\qeda

In the next statements we prove the existenc of radial solutions defined in the complement of a ball of $\BBR^N$, $N\geq 3$ which behaves like the Newtonian potential at infinity. We start with the following lemma dealing with the positive radial solutions of $\CL_pv=0$ in the complement of a ball.

\blemma{ext} Assume $N\geq 3$ and $p>\frac{N}{N-2}$. Then for any $c>0$ there exists $k_c>0$ such that the unique solution $v_c$ of $\CL_pv=0$ in 
$B_1^c$ verifying $v\lfloor_{\prt B_1}=c$ satisfies 
  \bel{L3}\displaystyle
  \lim_{|x|\to\infty}|x|^{N-2} v_c(x)=k_c.
  \ee
  Furthermore the mapping $c\mapsto k_c$ is continuous and increasing from $(0,\infty)$ onto $(0,k_\infty)$ for some $k_\infty<\infty$.
\es 
\Proof The existence and uniqueness of a solution $v_c$ in an exterior domain and the fact that $(\ref{L3})$ holds is classical (see e.g. \cite {Veasym}). However the fact that $k_c>0$ and the continuity of $c\mapsto k_c$ is not proved there.  By the maximum principle $c\mapsto k_c$ is nondecreasing. Next we set 
$s=\frac{r^{N-2}}{N-2}$ and $v_c(r)=r^{2-N}\gr(s)$. Then $\gr_c:=\gr$ satisfies 
  \bel{L4}\displaystyle
s^2\gr_{ss}=c_{N,p}s^{\frac{N}{N-2}-p}\gr^p\;\text{ on }\,((N-2)^{-1},\infty)\;\text{ and }\,\gr((N-2)^{-1})=c,
  \ee
  where $c_{N,p}=(N-2)^{\frac{4-N}{N-2}-p}$. By the maximum principle $v_c(r)\leq cr^{2-N}$ ($v_c$ is the positive harmonic function in $B_1^c$ with value $c$ on $\prt B_1$), hence $\gr (s)$ is bounded. Since $\gr$ is convex and bounded, it is decreasing and $\gr's)\to 0$ as $s\to \infty$. Hence 
    \bel{L5}-\gr'(s)=c_{N,p}\myint{s}{\infty}\gt^{\frac{N}{N-2}-p-2}\gr^p(\gt)d\gt\leq c'_{N,p}s^{\frac{N}{N-2}-p-1}\gr^p(s).
\ee
Hence by integration the function $\gr\mapsto\Gf(\gr)=\gr^{1-p}-c''\gr^{\frac{N}{N-2}-p}$ is increasing and bounded. Then it has a finite limit $\ell$ when $\gr\to\infty$ and $\Gf^{1-p}(\gr)$ has the same limit $\ell$. Thus $\ell\neq 0$ and 
consequently $k_c>0$. Let $\{c_n\}$ be a decreasing sequence in $\BBR_+$ converging to 
  $c^*>0$. Then the sequence of corresponding solutions $\{v_{c_n}\}$ is decreasing to $v_{c^*}$ the sequence $\{k_{c_n}\}$ is nonincreasing with limit 
  $k^*\geq k_{c^*}$.    From $(\ref{L5})$ one get
    \bel{L6}\displaystyle
c_n-k_{c_n}=c_{N,p}\myint{(N-2)^{-1}}{\infty}\myint{s}{\infty}\gt^{\frac{N}{N-2}-p-2}\gr_{c_n}^p(\gt)d\gt ds,
  \ee
  and the same identity holds in $c_n$ is replaced by $c^*$. By the dominated convergence theorem, one has that 
      \bel{L7}\displaystyle
c^*-k^*=c_{N,p}\myint{(N-2)^{-1}}{\infty}\myint{s}{\infty}\gt^{\frac{N}{N-2}-p-2}\gr_{c^*}^p(\gt)d\gt ds,
  \ee
  which implies that $k^*=k_{c^*}$. A similar result holds if $\{c_n\}$ is an increasing sequence in $\BBR_+$ converging to 
  $c_*>0$. Hence $c\mapsto k_c$ is increasing and continuous. When $c\uparrow\infty$ $v_c$ increases and converges to the unique positive solution 
  $v_\infty$ of $\CL_pv=0$ in $B_1^c$ such that $\displaystyle\lim_{r\to 1}v(r)=\infty$. Hence $k_c\uparrow k_\infty$ and $k_\infty<\infty$.\qeda\medskip
  
  \nind\Remark Since the equation $\CL_pv=0$ is invariant by the transformation $T_\ell$ defined in $(\ref{Z2})$, the ball $B_1$ can be replaced by $B_R$ for any 
  $R>0$. The range of $k_c$, that we call $k_{c,R}$ is modified accordingly and $\displaystyle\lim_{c\to\infty}k_{c,R}=k_{\infty,R}$. Then one has 
  \bel{L7bis}k_{\infty,R}=R^{N-2-\ga}k_{\infty,1}.\ee
  
\bth{dirac34} Let $N\geq 3$, $M>0$, $p>\frac{N}{N-2}$ and $\frac{N}{N-1}<q<p$. \smallskip

\nind 1- For any $k>0$ there exist $R:=R_k>0$ and a positive radial solution $u$ of $(\ref{Z1})$ in $B_R^c$ satisfying
  \bel{L8}\BA {lll}
      \displaystyle \lim_{|x|\to\infty}|x|^{N-2}u(x)=k.
            \EA\ee
            
\nind 2- If $\frac{2p}{p+1}<q<p$ there exist $\tilde k>0$ and a  positive radial solution, unique among all the positive solutions, $u$ of $(\ref{Z1})$ in $\BBR^N\setminus\{0\}$ satisfying $\displaystyle
\lim_{r\to 0}r^\gg u(r)=X_{_M}$  and $(\ref{L8})$ with $k=\tilde k$. In the particular case $q=\frac{(N-2)p}{N-1}$ we have $u=u^*_{_{M,p}}$ (see $(\ref{U11})$).\smallskip

\nind 3- If $\frac{2p}{p+1}<q<2$  and the assumption $(\ref{M5_*})$ of \rth{lamglo} is satisfied, 
            there exist $ k>0$ and a radial positive solution $u$ of $(\ref{Z1})$ in $\BBR^N\setminus\{0\}$ satisfying $(\ref{L8})$ and $\displaystyle \lim_{x\to 0}|x|^\gb u(x)=\xi_{_M}$. Furthermore $u$ is unique among all the positive solutions satisfying $(\ref{L8})$.
\es
\Proof 1- If $w$ is a positive radial and decreasing function such that $\CR^M_{q}w=0$ it satisfies (see $(\ref{T11})$)
$$-w'(r)=r^{1-N}\left(C+\myfrac{M}{\gk}r^{N-(N-1)q}\right)^{-\frac 1{q-1}},
$$
where, $\gk=\frac{(N-1)q-N}{q-1}$ and $C\in\BBR$. If $C>0$, $w$ is defined on $(0,\infty)$. Hence if $w(r)\to 0$ as $r\to\infty$, one has
   \bel{L10}\BA {lll}
w(r)=\myint{r}{\infty} s^{1-N}\left(C+\myfrac{M}{\gk}s^{N-(N-1)q}\right)^{-\frac 1{q-1}}ds.
            \EA\ee
            Then $w(r)=\frac{1}{(N-2)C^{\frac{1}{q-1}}}r^{2-N}(1+o(1))$ as $r\to\infty$. Hence, if $k>0$ is given, we choose $C>0$ such that 
            $\frac{1}{(N-2)C^{\frac{1}{q-1}}}=k$. In order that $k$ is in the range of the application $c\mapsto k_{c,R}$, one takes $R>0$ such that 
            $k<R^{N-2-\ga}k_{\infty,1}$. For such an $R$, there exists $c>0$ such that the solution $v_c$ of $\CL_pv=0$ in $B_R^c$ verifying 
            $v=c$ on $\prt B_R$ satisfies $(\ref{L3})$. We then set $C=\left(\frac{1}{(N-2)k}\right)^{q-1}$. The function $w:=w_C$ defined by $(\ref{L10})$ is a supersolution of $(\ref{Z1})$ in $B_R^c$, larger than the subsolution $v_c$ and both $v_c$ and $w_C$ satisfy $(\ref{L8})$. Then by \cite[Theorem 1.4.5]{Vebook} there exists a radial positive solution $u$ of $(\ref{Z1})$ in $B_R^c$ such that $v_c\leq u\leq w_C$, hence $(\ref{L8})$ follows.\\
            \nind 2- The existence of a unique positive and radial solution in $\BBR^N\setminus\{0\}$ satisfying $(\ref{Z19})$ follows from \rth{T10}. The asymptotic behaviour is a consequence of \rth{T8}-(2).\\
            \nind 3- Under the condition $(\ref{M5_*})$ of \rth{lamglo} there exists a unique positive solution in $\BBR^N\setminus\{0\}$ satisfying $(\ref{M5})$. From \rth{T8}, and since $p>\frac{N}{N-2}$, its behaviour at infinity is given by $(\ref{Z27})$ for some specific $k^*>0$. Uniqueness follows from $\rth{uni}$.
           \qeda
\mysection{Isolated singularities of non-radial solutions}
\subsection{Existence and uniqueness of singular solutions}
The results of this paragraph are independent of the description of the radial singular solutions performed in the previous sections and they provide a general tool for constructing singular solutions. The existence of singular solutions is based upon the next variant of  \cite[Theorem 2.1]{BMP} proved in \cite[Corollary 1.4.5]{Vebook}.
\bth{BMPLV} Let $G$ be a bounded domain in $\BBR^N$, $B\in C(G\ti\BBR\ti\BBR^N)$ a real valued function, $\Gg\in C(\BBR_+,\BBR_+)$ an increasing function 
such that
\begin{equation}\label{X1}
|B(x,r,\xi)|\leq \Gg(|r|)(1+|\xi|^2)\quad\text{for all }(x,r,\xi)\in G\ti\BBR\ti\BBR^N.
\end{equation}
Let $\CQ$ be the operator defined by 
\begin{equation}\label{X2}
\CQ(u)=-\Gd u+B(x,u,\nabla u).
\end{equation}
If there exist a supersolution $\phi\in W^{1,\infty}(G)$ and a subsolution $\psi\in W^{1,\infty}(G)$ such that $\psi\leq\phi$, then for any 
$\chi\in W^{1,\infty}(G)$
satisfying $\psi\leq\chi\leq\phi$ there exists a function $u\in W^{1,2}(G)$ verifying $\psi\leq u\leq\phi$, solution of $\CQ(u)=0$ and such that $u-\chi\in W^{1,2}_0(G)$.
\es

One of the main application of this result is \rth{Exist1} which is proved below\medskip

\nind{\it Proof of \rth{Exist1}}. 
 Let $\{\ge_n\}$ be a sequence decreasing to $0$ and such that $\ge_1< \dist (0,\prt\Gw)$ and set 
$\displaystyle m=\max_{z\in\prt\Gw}v(z)+\max_{z\in\prt\Gw}(\phi(z)-v(z))_+$. Then $m\geq\gf$ on $ \prt\Gw$ and the function 
$\overline v=v+m$ satisfies $\CL_{p,q}^M\overline v\geq 0$ in $\Gw\setminus \{0\}$. The function $\displaystyle\underline v=(v-\max_{z\in\prt\Gw}\phi(z))_+$
satisfies $\CL_{p,q}^M\underline v\leq 0$. Put $\chi=\sup\left\{\underline v,\inf\left\{\overline v,\phi\right\}\right\}$.  Then $\chi\in W_{loc}^{1,\infty}(\Gw\setminus \{0\})$, $\underline v\leq\chi\leq \overline v$ and $\chi=\phi$ on $\prt\Gw$. By \rth{BMPLV} for any $n\in\BBN^*$ there exists a function 
$u_{n}\in W^{1,2}(\Gw\setminus \overline B_{\ge_n})$ such that $(u_n-\chi)\lfloor_{\Gw\setminus \overline B_{\ge_n}}\in W^{1,2}_0(\Gw\setminus \overline B_{\ge_n})$ satisfying $\CL_{p,q}^Mu_n=0$ in $\Gw\setminus \overline B_{\ge_n}$. Furthermore $u_n$ is unique by the maximum principle. Since 
$u_n=\overline v$ on $\prt B_{\ge_n}$, $v$, and therefore $\overline v$, is radially decreasing  and $u_n=\chi$ on $\prt\Gw$ we infer that $u_n\leq u_{n'}$ in $\Gw\setminus B_{\ge_n}$ if $n'\geq n$. Hence the sequence $\{u_n\}$ is increasing and it satisfies 
\begin{equation}\label{X4}\displaystyle 
(v(x)-\max_{z\in\prt\Gw}\phi(z))_+\leq u_n(x)\leq v(x)+\max_{z\in\prt\Gw}(\phi(z)-v(z))_+\quad\text{for all }x\in \Gw\setminus B_{\ge_n}.
\end{equation}
 By standard regularity estimates, $u_n$ is relatively compact in $C_{loc}^1(\Gw\setminus \{0\})$. Hence it converges to a solution $u$ of $\CL_{p,q}^Mu=0$ in $\Gw\setminus \{0\}$ which coincides with $\phi$ on $\prt\Gw$ and satisfies $(\ref{X3})$.\qeda\medskip
 
 As a  first application we have the following:
\bcor {Bd} Let $\Gw$ be any  bounded smooth domain containing $0$ and $\gf\in W^{1,\infty}(\Gw)$ be nonnegative. There exists a positive solution $u$ of $\CL_{p,\frac{2p}{p+1}}^Mu=0$ in $\Gw\setminus \{0\}$ with value $\gf$ on $\prt\Gw$ such that $u(x)-a|x|^{-\ga}$ remains bounded in $\Gw$ where 
 $a$ is equal to $x_{_M}$ or $x_{_{j,M}}$ (j=1,2) or $x_{m^*}$ according to we are in the cases (1)-(2) or (3) or (4) of \rth{T1}.\es
 
 The existence of singular solutions is not restricted to the case $q=\frac{2p}{p+1}$ where they are explicit. The following easy to prove corollary shows that 
 existence, and sometimes uniqueness, holds when $1<q<\frac{N}{N-1}$. This range of exponents is analysed in \cite{BVGHV3} in connection with problems with Dirac measure data.
 
\bcor {Bd1} Let $\Gw\subset\BBR^N$, $N\geq 1$, be any  bounded smooth domain containing $0$. Assume $1<p<\frac N{N-2}$ if $N\geq 3$ or any $p>1$ if $N=1,2$, $1<q<\min\left\{p,\frac N{N-1}\right\}$ if $N\geq 2$ or any $q>1$ if $N=1$, $M>0$ and $k>0$. Then for any $\gf\in W^{1,\infty}(\Gw)$, $\gf\geq 0$, there exists a positive solution $u$ of $\CL_{p,q}^Mu=0$ in $\Gw\setminus \{0\}$ with value $\gf$ on $\prt\Gw$ satisfying $(\ref{Z8})$. \es

\nind{\it Proof of \rth{T2} and \rth{T3}}.  It is a direct consequence of the above results.\qeda\medskip

More general uniqueness results valid for any positive solution, not necessarily radial, are obtained below. Furthermore the problems involved are either considered in $\BBR^N\setminus\{0\}$ or in a punctured bounded domain. If $b$ is a positive parameter we define a continuous group of transformations acting on functions $u$ defined in an open set $G$, $u\mapsto u_{\ell}^{(b)}$, for $\ell>0$ by the formula
\bel{UN1}
\displaystyle u_{\ell}^{(b)}(x)=\ell^bu(\ell x)\quad\text{ for all $\,\ell>0$ and $\,x\,$}\in\ell^{-1}G.
\ee
If $u$ satisfies $(\ref{Z1})$ in $G$, then $u_{\ell}^{(b)}$ satisfies 
\bel{UN2}
-\Gd u_{\ell}^{(b)}+\ell^{2-b(p-1)}(u_{\ell}^{(b)})^p-M\ell^{2-q-b(q-1)}|\nabla{u_{\ell}^{(b)}}|^q=0\quad\text{in }\;\ell^{-1}G.
\ee
If $\ell>1$, $u_{\ell}^{(b)}$ is a supersolution of $(\ref{Z1})$ if and only if
\bel{UN3}\BA{lll}
(i)\qquad &2-b(p-1)\leq 0&\Longleftrightarrow \ga\leq b\phantom{--------------}\\[2mm]
(ii)&2-q-b(q-1)\geq 0&\Longleftrightarrow \gb\geq b.
\EA\ee
This conditions are compatible if and only if $1<q\leq \frac{2p}{p+1}$. Similarly, if $\ell<1$, $u_{\ell}^{(b)}$ is a supersolution of $(\ref{Z1})$ if and only if
\bel{UN4}\BA{lll}
(i)\qquad &2-b(p-1)\geq 0&\Longleftrightarrow \ga\geq b\phantom{--------------}\\[2mm]
(ii)&2-q-b(q-1)\leq 0&\Longleftrightarrow \gb\leq b.
\EA\ee
This conditions are compatible if and only if $\frac{2p}{p+1}\leq q<2$. \medskip

\nind{\it Proof of \rth{uni1}} First we note that two terms on the right hand-side of $(\ref{UNX})$ in the statement of the theorem coincide only if $q=\frac{2p}{p+1}$ since $\ga\leq\gb$ is equivalent to $q\leq\frac{2p}{p+1}$ . We first study the problem in 
$\BBR^N\setminus\{0\}$. We have to consider two cases:\\
\nind 1- Suppose $\ga\leq\gb$. We choose $b$ such that 
\bel{UNX2}b\in \left(a,\infty\right)\cap \left[\ga,\gb\right]. 
\ee
Let $u$ and $\tilde u$ be two positive solutions satisfying $(\ref{UNX1})$. For $\ell>1$, $u_{\ell}^{(b)}$ is a supersolution. Since 
$$u_{\ell}^{(b)}(x)=\Gl\ell^{b-a}|x|^{-a}(-\ln |x|)^{\tilde a} (1+o(1))\quad\text{as }x\to 0
$$
and $u(x)\to 0$ as $|x|\to\infty$, for any $\ge>0$ the function $x\mapsto u_{\ell}^{(b)}(x)+\ge$ which is a supersolution is larger than $\tilde u$ near $0$ and at infinity. Then $u_{\ell}^{(b)}+\ge\geq \tilde u$ in $\BBR^N\setminus\{0\}$. Letting $\ge \downarrow 0$ and $\ell\downarrow 1$, yields $u\geq\tilde u$. Similarly 
$\tilde u\geq u$. \\
\nind 2- Suppose $\ga>\gb$. We choose $b$ such that 
\bel{UNX3}b\in \left(0,a\right)\cap \left[\gb,\ga\right]. 
\ee
Then for $\ell<1$, $u_{\ell}^{(b)}+\ge$ is a supersolution in $\BBR^N\setminus\{0\}$ which is larger than $\tilde u$ at $0$ and at $\infty$. Hence 
$\tilde u\leq u_{\ell}^{(b)}+\ge$ and we conclude as in the first case.\\
Next we consider the problem in $\Gw$. Since the solutions are continuous in $\overline\Gw\setminus\{0\}$, for $\ge>0$ we have that for $\ell>1$ 
$u_{\ell}^{(b)}+\ge >\tilde u$ near $\prt(\ell^{-1}\Gw)$ provided $\ell-1$ is small enough. Hence $u_{\ell}^{(b)}(x)+\ge \geq\tilde u$ in $\ell^{-1}\Gw\setminus\{0\}$.
This implies that $u\geq \tilde u$ by letting $\ell\uparrow 1$ and then $\ge\to 0$. If $\ell<1$ then $\Gw\subset \ell^{-1}\Gw$, and we compare $u_{\ell}^{(b)}+\ge$
and $\tilde u$ in $\Gw$. The proof follows.\qeda\medskip

The previous result necessitates to find some $b$ satisfying either $(\ref{UNX2})$ or $(\ref{UNX3})$ which is not always possible in practice. We give below a variant of the result 
which necessitates a slightly sharper blow-up estimate. 
\bth{uni2} Assume $N\geq 1$, $p>1$, $1<q\leq \tfrac{2p}{p+1}$ and $M>0$. Let $a$ such that 
\bel{UNX'}
0\leq a\leq\gb
\ee
There exists at most one positive solution of $(\ref{Z1})$ in $\BBR^N\setminus\{0\}$ satisfying
\bel{UNX'1}
u(x)=\Gl|x|^{-a}+\Gl'|x|^{-a'}(1+o(1))\quad\text{as }x\to 0,
\ee
or
\bel{UNX''1}
u(x)=\Gl|x|^{-a}(-\ln |x|)^{-a''}(1+o(1))\quad\text{as }x\to 0,
\ee
where $\Gl,\Gl'$ are some positive constants and $a>a'$ and $a''>0$.
\es
\Proof The principle of the proof is to replace $(\ref{UNX2})$ by
\bel{UNX'2}b\in \left[a,\infty\right)\cap \left[\ga,\gb\right]. 
\ee
when $\ga\leq\gb$. Then, for $\ell>1$, $u_{\ell}^{(b)}$ is a supersolution. If $u$ satisfies $(\ref{UNX'1})$ then, as $x\to 0$,
$$
u_{\ell}^{(b)}(x)=\Gl\ell^{b-a}|x|^{-a}+\Gl'\ell^{b-a'}|x|^{-a}|x|^{-a'}(1+o(1))
$$
Since $b\geq a>a'$, $u_{\ell}^{(b)}$ is larger than another solution $\tilde u$ near $0$. Thus $u_{\ell}^{(b)}+\ge\geq\tilde u$ for any $\ge>0$, which implies the claim. \\
If $u$ satisfies $(\ref{UNX''1})$, then
$$\BA{lll}
u_{\ell}^{(b)}(x)=\Gl\ell^{b-a}|x|^{-a}(-\ln|x|-\ln\ell)^{-a''}(1+o(1))\\[2mm]
\phantom{u_{\ell}^{(b)}(x)}
=\Gl\ell^{b-a}|x|^{-a}(-\ln |x|)^{-a''}\left(1+a''\myfrac{\ln\ell}{-\ln|x|}\right)(1+o(1)).
\EA$$
Again $u_{\ell}^{(b)}$ is larger than another solution $\tilde u$ in a neighborhood of $0$ and we end the proof as in the first case. \qeda\medskip

\nind\Remark The method developed above allows to give uniqueness result for large solutions under some starshapedness assumption. Let $G\subset\BBR^N$ be a domain with compact boundary and $\gr_G(x)=\dist(x,\prt G)$, we consider the problem
\bel{UNX4}\BA{lll}
-\Gd u+u^p-M|\nabla u|^q=0\quad\text{in }G\\\phantom{---,;;}
\displaystyle\!\lim_{\gr_G(x)\to 0}u(x)=\infty.
\EA\ee
Such a solution, if it exists is called a {\it large solution}.
\bth{uni3} Assume $N\geq 1$, $M>0$ and $p,q>1$ and $\Gw$ is a bounded domain starshaped with respect to $0$.  There exists at most one positive function satisfying $(\ref {UNX4})$ in one of the following case:\smallskip

\nind 1- $\frac{2p}{p+1}\leq q <2$ and $G=\Gw$.\smallskip

\nind 2- $1<q\leq \frac{2p}{p+1}$ and $G=\overline\Gw^c$.
\es
\Proof Let $u$ and $\tilde u$ be two positive solutions of $(\ref {UNX4})$. In the first case with $G=\Gw$. Then for $\ga\leq b\leq\gb$ and $\ell>1$, $u^{b}_\ell$ is a supersolution of $(\ref {UNX4})$ in $\Gw_\ell:=\ell^{-1}\Gw$. Since 
$\overline \Gw_\ell\subset \Gw$, $u^{b}_\ell\geq \tilde u$, it follows that $u\geq\tilde u$.\\
In the second case with $G=\overline\Gw^c$, then for $0<\ell<1$ and $\gb\leq b\leq\ga$. Then $u^{b}_\ell$ is a supersolution in $\ell^{-1}\overline\Gw^c\subset \overline\Gw^c$. Then for $\ge>0$, $u^{b}_\ell+\ge\geq \tilde u$ in $\ell^{-1}\overline\Gw^c$. Letting $\ge\to 0$ and $\ell\uparrow 1$ yields $u\geq\tilde u$.
This ends the proof.\qeda\medskip

If we combine the results of existence of radial singular solutions in $\BBR^N\setminus\{0\}$ with the uniqueness results of \rth{uni1} and \rth{uni2} we have the following:
\bth{uni} Assume $N\geq 3$, $p,q>1$ and $M>0$. There exists one and only one positive solution $u$ of $(\ref{Z1})$ in $\BBR^N\setminus\{0\}$, if one of the following conditions holds:\smallskip

\nind 1- $1<p<\tfrac{N}{N-2}$, $q=\frac{2p}{p+1}$, $M>0$ and $u$ satisfies  $(\ref{Z8})$-(i) for some $k>0$.\\
\nind 2- $\frac{2p}{p+1}<q<p$ and $\displaystyle \lim_{x\to 0}|x|^\gg u(x)=X_{_M}$. \\
\nind 3- $1<p<\tfrac{N}{N-2}$, $1<q<\frac{2p}{p+1}$, $M>0$ and either $\displaystyle \lim_{x\to 0}|x|^\ga u(x)=x_0$, or $u$ satisfies  $(\ref{Z8})$-(i) for some $k>0$.\\
\nind 4- $p=\frac{N}{N-2}$, $q=\frac{2p}{p+1}$, $M>0$  and $u$ satisfies  $(\ref{Z11})$-(i). 
\smallskip

\nind Furthermore, existence and uniqueness of a solution holds if the equation $(\ref{Z1})$ is considered in $\Gw\setminus\{0\}$ where $\Gw$ is a bounded smooth domain starshaped with respect to $0$ and is the function $u$ is equal to some $\gf$ on $\prt\Gw$ where $\gf\in C^{1}(\prt\Gw)$ is nonnegative. 

\es
\Proof By applying \rth{Exist1} and \rth{uni1} the proof is reduced to use results of existence of radial positive singular solutions in $\BBR^N\setminus\{0\}$ and to check that the parameters fulfill the conditions of \rth{uni1}. \\
Case 1- If $q=\tfrac{2p}{p+1}$, $\ga=\gb$ and $N-2<\gb$. Existence of radial positive solutions satisfying $(\ref{Z8})$-(i) is proved in \rth{T3}. \\
Case 2- Then $\gg>\gb$. Existence of a radial positive solution satisfying $\displaystyle \lim_{x\to 0}|x|^\gg u(x)=X_{_M}$ is proved in \rth{T10}-1.\\
Case 3- If $1<q<\frac{2p}{p+1}$, then $\ga<\gb$, and since $p<\frac{N}{N-2}$, $N-2<\ga$. Existence of a a radial positive solution satisfying $\displaystyle \lim_{x\to 0}|x|^\ga u(x)=x_0$ is proved in \rth{T12}-1. Furthermore the assumptions on $p$ and $q$ imply that $q<\frac{N}{N-1}$. The existence of a positive solution satisfying $(\ref{Z8})$-(i) for any $k>0$ is proved in \rth{weaksing}.
\\
Case 4- When $p=\frac{N}{N-2}$, $q=\frac{2p}{p+1}$ there exists a radial global solution satisfying $(\ref{Z11})$-(i) by \rth{T4}. We apply estimate $(\ref{UNX''1})$ in \rth{uni2} with $a=N-2=\ga=\gb$ and $a''=N-1$. The result follows.\qeda\medskip

\nind\Remark In the case $p=\frac{N}{N-2}$, $q=\frac{N}{N-1}$ we conjecture that the function $u_{x_M}$ is the only positive solution of $(\ref{Z1})$ defined in $\BBR^N\setminus\{0\}$ satisfying $\displaystyle\lim_{x\to 0}|x|^{N-2}u(x)= x_{_M}$.

\subsection{Characterization of singular solutions}
In this section we give some results showing how the characterization of singularities of radial solutions can be extended to nonradial solutions. An important tool for 
studying positive isolated singularities is Harnack inequality.
\bprop{harnak} Assume $M>0$, $p>1$ and $1<q\leq \frac{2p}{p+1}$. If $u$ is a positive solution of $(\ref{Z1})$ in $B_{R_0}\setminus\{0\}$, there exists $c_9=c_9(N,p,q,R_0,M)>0$ such that for any $R\in (0,\frac{R_0}{2}]$ there holds
   \bel{K1}\BA {lll}
\displaystyle\max_{|x|=R}u(x)\leq c_9\displaystyle\min_{|x|=R}u(x).
            \EA\ee
\es
\Proof We write $(\ref{Z1})$ under the form
   \bel{K2}\BA {lll}
-\Gd u+C(x)u+V(x)|\nabla u|=0,
            \EA\ee
where $C(x)=|u(x)|^{p-1}$ and $V(x)=M|\nabla u(x)|^{q-1}$. By \rprop{Oss}, 
$$C(x)\leq c_1^{p-1}\max\left\{M^{\frac{p-1}{p-q}}|x|^{-\frac{q(p-1)}{p-q}},|x|^{-2}\right\}
$$
and 
$$V(x)\leq c_4^{q-1}\max\left\{M^{\frac{q-1}{p-q}}|x|^{-\frac{p(q-1)}{p-q}},|x|^{-\frac{(q-1)(p+1)}{(p-1)}}\right\}
$$
Under the assumptions $1<q\leq \frac{2p}{p+1}$, the terms  $|x|^2 C(x)$ and $|x| V(x)$ are uniformly bounded in $B_{R_0}\setminus\{0\}$. The result follows by
\cite[Chapter 8]{GT}.\qeda 

\subsection{The case $1<p<\frac{N}{N-2}$ and $1<q<\frac{N}{N-1}$}

In this section, the results are obtained by a combination of \rth{weaksing} for existence of solutions and \rth{uni} for their uniqueness.
 
\bth{18M} Let $N\geq 3$, $M>0$, $1<p<\frac{N}{N-2}$, $1<q<\frac{N}{N-1}$ and $\Gw$ be a bounded domain containing $0$. For any $k>0$ there exists 
a unique positive function $u:=u_k$ solution of $(\ref{Z1})$ in $\Gw\setminus\{0\}$, vanishing on $\prt\Gw$ and satisfying 
\bel{K5}\displaystyle\lim_{x\to 0}|x|^{N-2}u(x)=k.
\ee
 Furthermore 
$k\mapsto u_k$ is increasing by the maximum principle and converges to a solution $u_\infty$ of $(\ref{Z1})$ in $\Gw\setminus\{0\}$, vanishing on $\prt\Gw$ and satisfying
\bel{K3}\BA {lll}
\displaystyle \lim_{x\to 0}|x|^\gg u(x)=X_M\qquad\text{if }\; \frac{2p}{p+1}<q<\frac{N}{N-1},
\EA\ee
where $X_{_M}$ is defined at $(\ref{Z19})$,  or
\bel{K4}
\displaystyle \lim_{x\to 0}|x|^\ga u(x)=\left\{\BA{lll}x_0\qquad&\text{if }\; 1<q<\frac{2p}{p+1}\\
x_{_M}\qquad&\text{if }\; q=\frac{2p}{p+1},
\EA\right.
\ee
where $x_{0}$, $x_{_M}$ are the unique positive root of equation $(\ref{Z3})$ with $M=0$ and $M>0$ respectively.
\es
\Proof Let $0<R_1<R_2$ be such that $\overline B_{R_1}\subset\Gw\subset\overline\Gw\subset B_{R_2}$. By \rth{T12}-1 and \rth{uni}, for $k>0$ there exists a unique solution $u_{1,k}$ (resp. $ u_{2,k}$) of $(\ref{Z1})$ in $B_{R_1}\setminus\{0\}$ (resp. $B_{R_2}\setminus\{0\}$)   satisfying $(\ref{Z8})$-(i) and vanishing on $\prt B_{R_1}$ (resp. $\prt B_{R_2}$).
If we extend $u_{1,k}$ by $0$ in  $B^c_{R_1}$, we have,
$$u_{1,k}(|x|)\leq u_k(x)\leq u_{2,k}(|x|)\quad \text{in }\,\Gw\setminus\{0\}.
$$
All the above functions are locally bounded in $\overline \Gw\setminus\{0\}$ and $\overline B_{R_2}\setminus\{0\}$ by \rprop{Oss}. Since the mappings $k\mapsto u_k$ and $k\mapsto u_{j,k}$ are increasing, we have, by letting $k\to\infty$, 
$$u_{1,\infty}(|x|)\leq u_\infty(x)\leq u_{2,\infty}(|x|)\quad \text{in }\,\Gw\setminus\{0\}.
$$
Then we obtain $(\ref{K3})$ by \rth{T6}-(1) and $(\ref{K4})$ by \rth{T7}-(1) and \rth{T3}.\qeda\medskip

The main characterization of isolated singularities is the next result.
\bth{19M} Let $N\geq 3$, $\Gw$ be an open subset containing $0$, $M>0$, $1<p<\frac{N}{N-2}$ and $1<q\leq \frac{2p}{p+1}$. If $u$ is a positive solution of $(\ref{Z1})$ in $\Gw\setminus\{0\}$, then either its behaviour at $x=0$ is given by $(\ref{K3})$ or $(\ref{K4})$, or there exists $k\geq 0$ such that $(\ref{K5})$ holds.
If $k=0$ the singularity at $0$ is removable.
\es
The proof needs a few intermediate steps.
\blemma{12adapt} Let $N\geq 3$, $\Gw$ be an open subset containing $0$, $M>0$, $1<p<\frac{N}{N-2}$ and $1<q< \frac{N}{N-1}$. If $u$ is a positive solution of $(\ref{Z1})$ in $\Gw\setminus\{0\}$ vanishing 
on $\prt\Gw$ and such that 
$$\displaystyle\limsup_{x\to 0}|x|^{N-2}u(x)<\infty.
$$
Then there exists $k\geq 0$ such that $(\ref{K5})$ holds.
If $k=0$, then $u$ coincides in $\Gw\setminus\{0\}$ with a $C^2(\Gw)$ solution of $(\ref{Z1})$ in $\Gw$. 
\es
\Proof By assumption $u(x)\leq c|x|^{2-N}$ and by \cite[Lemma 3.10]{BVGHV} we have the following: if $u$ is a solution of $(\ref{Z1})$ in $\Gw\setminus\{0\}$ (not necessarily positive) such that $|x|^m|u(x)|$ is bounded near $x=0$ for some $m<\inf\{\ga, \frac{2-q}{q-1}\}$, then $|x|^{m+1}|\nabla u(x)|$ is also bounded near $x=0$. Actually, in the reference the result is proved for a more general operator, without the absorption $u^p$, but the adaptation is straightforward. The result applies there with $m=N-2$ and in particular $|\nabla u(x)|\leq c'|x|^{1-N}$. We write $(\ref{Z1})$ under the form $(\ref{K2})$. Since $p<\frac{N}{N-2}$ and $q<\frac{N}{N-1}$ we have
$$|C(x)|\leq c|x|^{(2-N)(p-1)}\leq c|x|^{-2+\ge_1}\,\text{ and }\;|V(x)|\leq c|x|^{(1-N)(q-1)}\leq c|x|^{-1+\ge_2},$$ 
for some $\ge_1,\ge_2>0$. It follows by Serrin's result that either the singularity at $0$ is removable, or there exist $c_1>c_2>0$ such that 
$$c_2|x|^{2-N}\leq u(x)\leq c_1|x|^{2-N}\quad\text{for all }0<|x|\leq 1.
$$
In order to make the convergence precise, we denote by $u_2$ the solution of 
$$\BA{lll}
-\Gd u_2=M|\nabla u|^q-u^p\quad&\text{in } B_1\\\phantom{-\Gd }
u_2=0&\text{on } \prt B_1. 
\EA$$
Then $-v_2\leq u_2\leq v'_2$ where 
$$\BA{lll}
-\Gd v'_2=Mc'^q|x|^{(1-N)q}\quad&\text{in } B_1\\\phantom{-\Gd }
v'_2=0&\text{on } \prt B_1, 
\EA$$
and 
$$\BA{lll}
-\Gd v_2=c^p|x|^{(2-N)p}\quad&\text{in } B_1\\\phantom{-\Gd }
v_2=0&\text{on } \prt B_1, 
\EA$$
Because $(N-1)q<N$ and $(N-2)p<N$, $v_2$ and $v'_2$ satisfy
$$0\leq v_2(r)\leq c_{10}r^{2-N+\gd}\,\text{ and }\;0\leq v'_2(r)\leq c_{10}r^{2-N+\gd'},
$$
for some $c_{10}>0$ and $\gd=N-(N-2)p>0$, $\gd'=N-(N-1)q>0$. Then $u_2$ satisfies 
$$\displaystyle 
\lim_{x\to 0}|x|^{N-2}u_2(x)=0. 
$$
The function $u_1=u-u_2$ is harmonic in $B_1\setminus\{0\}$ and is bounded from below by $-v'_2$ which satisfies $\displaystyle 
\lim_{x\to 0}|x|^{N-2}v'_2(x)=0$. Hence by standard result on singularities of harmonic functions, $|x|^{N-2}u_1(x)$ admits a limit 
$k\geq 0$ when $x\to 0$. Combined with Serrin's estimates it follows  that either  $k=0$ and the singularity is removable, 
or $k>0$. Note that if $k=0$, then $u$ is a $C^2$ solution in $\Gw$.\\
Another proof based on a perturbation is the following: let $u(x)=u(r,s)=|x|^{2-N}\gf(t,s)$ with $r=|x|$ and $t=\ln r$. Then 
$$\BA{lll}\gf_{tt}+(N-2)\gf_t-e^{(N-p(N-2))t}\gf^p+\Gd'\phi\\[2mm]
\phantom{---------}
+Me^{(N-q(N-1))t}\left(((N-2)\gf-\gf_t)^2+|\nabla'\gf|^2\right)^{\frac q2}=0
\EA$$
Since $u=O(|x|^{2-N})$, we can write
$$\gf_{tt}+(N-2)\gf_t+\Gd'\phi=-e^{at}\psi,
$$
where $a=\min\{(N-p(N-2),(N-q(N-1)\}>0$ and $\psi$ is bounded. Then the result follows by \cite[Proposition 4.1]{BiRa}. 
\qeda\medskip

We give below another application of the perturbation method and specific to the case $q<\frac{2p}{p+1}$.
\bprop{reduc} Assume $\Gw$ is an open subset containig $0$, $M>0$, $1<p<\frac{N}{N-2}$ and $1<q<\frac{2p}{p+1}$. If $u$ is a solution of
\bel{K6_}-\Gd u+|u|^{p-1}u-M|\nabla u|^{q}=0
\ee
not necessarily nonnegative in $\Gw\setminus\{0\}$, then $r^\ga u(r,s)$ converges uniformly with respect to $s\in S^{N-1}$   when $r\to 0$ to a non-empty compact and connected subset of the set of solutions $\gw$ of
\bel{K6}
-\Gd'\gw+\ga (N-2-\ga)\gw+|\gw|^{p-1}\gw=0\qquad\text{on }\, S^{N-1}.
\ee
If $u\geq 0$, $\gw$ is either $x_0$ or $0$.
\es
\Proof We can assume that $\overline B_1\subset\Gw$ and set $\phi(t,s)=r^{\ga}u(r,s)$ with $t=\ln r$. Then $\phi$ satisfies
\bel{K7}\BA {lll}
\phi_{tt}+L\phi_t+\ga K\phi+\Gd'\phi-|\phi|^{p-1}\phi
+Me^{\frac{\gs t}{p-1}}\left((\ga\gf-\phi_t)^2+|\nabla '\phi|^2\right)^{\frac q2}=0.
\EA\ee
in $\BBR_-\ti S^{N-1}$ where $K=N-2-\ga$ and $L=K-\ga$. By assumption $\gs<0$, hence $(\ref{K7})$ is an exponentially small perturbation of the autonomous equation associated to the Emden-Fowler equation by the same change of variables. The result follows from \cite[Theorem 4.1]{BiRa} but for the sake of comprehension, we recall its proof. By \rprop {Oss-} the function $\phi$ is uniformly bounded, and by $(\ref{Y3''})$ $\phi_t$ and $\nabla'\phi$ are also uniformly bounded. By standard local regularity theory for elliptic equations, there holds
$$\norm{\prt_{t^{i}}\nabla'^j\phi}_{C[T-1,T+1]\ti S^{N-1}}\leq c_{11}\quad\text{for all }(i,j)\in \BBN\ti\BBN, i+j\leq 3\quad\text{and } T\leq -2,
$$
where $\nabla'^j$ is the covariant derivative of order $j$ on $S^{N-1}$. Then the omega-limit at $-\infty$ of the trajectory $\{\phi(t,.)\}_{t\in\BBR_-}$ in $C^2(S^{N-1})$ is a non-empty compact connected denoted by $\Gg_\phi\subset C^2(S^{N-1})$. From $(\ref{K7})$ we have that
$$\myfrac{d}{dt}E=L\myint{S^{N-1}}{}\phi^2_tdS(s)+Me^{\frac{\gs t}{p-1}}\myint{S^{N-1}}{}\left((\ga\gf-\phi_t)^2+|\nabla '\phi|^2\right)^{\frac q2}\phi_tdS(s),
$$
where 
$$J(t)=\myfrac{1}{2}\myint{S^{N-1}}{}\left(|\nabla'\phi|^2+\myfrac{2}{p+1}|\phi|^{p+1}-\phi^2_t\right) dS(s).
$$
Because $L\neq 0$ and $J$ is uniformly bounded, there holds
$$\myint{-\infty}{1}\myint{S^{N-1}}{}\phi^2_tdS(s)dt<\infty.
$$
Multiplying $(\ref{K7})$ by $w_{tt}$ and using the previous estimate, we obtain 
$$\myint{-\infty}{1}\myint{S^{N-1}}{}\phi^2_{tt}dS(s)dt<\infty.
$$
As $\phi_t$ and $\phi_{tt}$ are uniformly continuous in $(-\infty,1]\ti S^{N-1}$, this implies that
$$\displaystyle 
\lim_{t\to-\infty}\myint{S^{N-1}}{}\left(\phi^2_t(t,.)+\phi^2_{tt}(t,.)\right)dS(s)=0.
$$
Therefore $\Gg_\phi\subset\left\{\gw\in C^2(S^{N-1}): \ga K\gw+\Gd'\gw-|\gw|^{p-1}\gw=0\text{ in }S^{N-1}\right\}$. If $u\geq 0$, then $\gw\geq 0$ for any $\gw\in \Gg_\phi$ and the result follows by the maximum principle.\qeda\medskip


\blemma{13} Let $\Gw\subset\BBR^N$, $N\geq 3$, be a bounded open subset containing $0$, $M>0$, $1<p<\frac{N}{N-2}$ and $1<q<\frac{2p}{p+1}$. If $u$ is a 
nonnegative solution of $(\ref{Z1})$ in $\Gw\setminus\{0\}$ such that 
\bel{K8}\displaystyle 
\limsup_{x\to 0}|x|^{N-2}u(x)=\infty,
\ee
then 
\bel{K9}\displaystyle 
\lim_{x\to 0}|x|^{\ga}u(x)=x_0.
\ee
\es
\Proof  Without loss of generality we can assume that $\overline B_1\subset\Gw$. It follows from \rprop{harnak} that $\displaystyle 
\liminf_{x\to 0}|x|^{N-2}u(x)=\infty,
$
Hence there exists a decreasing sequence $\{r_n\}$ converging to $0$ such that 
$$\displaystyle \ga_n:=\inf_{|x|=r_n}|x|^{N-2}u(x)\uparrow\infty\quad\text{as }n\to\infty.
$$
Let $v_n$ be the solution of 
\bel{K10}\BA{lll}
-\Gd v+v^p=c_N\ga_n\gd_0\quad&\text{in } \CD'(B_1)\\
\phantom{-\Gd +v^p}
v=0\quad&\text{on } \prt B_1,\\
\EA\ee
where $c_N$ is the explicit constant such that $\BBG[c_N\gd_0](x)=|x|^{2-N}$, where $\BBG$ is the Newtonian potential in $\BBR^N$ ; note that $v_n$ is radial because of uniqueness. Then the sequence $\{v_n\}$ is increasing and converges to the function $v_\infty$
which satisfies (see \cite{Vesingsol})
\bel{K11}\BA{lll}
\phantom{i}-\Gd v+v^p=0\quad&\text{in } B_1\setminus\{0\}\\
\phantom{-\Gd i+iv^pi}
v=0\quad&\text{on } \prt B_1\\
\displaystyle\lim_{x\to 0}|x|^\ga v(x)=x_0.
\EA\ee
Moreover, $v_n(x)\leq \ga_n|x|^{2-N}$ and since $v_n$ is a subsolution of $(\ref{Z1})$, we have that $v_n\leq u$ in $\overline B_1\cap  B^c_{r_n}$. Letting $n\to\infty$ implies that $v_\infty(x)\leq u(x)$ in $ B_1\setminus\{0\}$. Therefore 
$$\displaystyle\liminf_{x\to 0}|x|^\ga u(x)\geq x_0.
$$
Combined with \rprop{reduc}, this inequality implies $(\ref{K7})$.\qeda\medskip

This result admits an extension to the case $q=\frac{2p}{p+1}$. 
\blemma{14} Let $\Gw\subset\BBR^N$, $N\geq 3$, be a bounded open subset containing $0$, $M>0$, $1<p<\frac{N}{N-2}$ and $q=\frac{2p}{p+1}$. If $u$ is a 
nonnegative solution of $(\ref{Z1})$ in $\Gw\setminus\{0\}$ such that 
\bel{K12}\displaystyle 
\limsup_{x\to 0}|x|^{N-2}u(x)=\infty,
\ee
then 
\bel{K13}\displaystyle 
\lim_{x\to 0}|x|^{\ga}u(x)=x_{_M}.
\ee
\es
\Proof Assuming for simplicity that  $\overline B_1\subset\Gw$ and using \rprop{harnak} we obtain that for some
  decreasing sequence $\{r_n\}$ converging to $0$ we have 
$$\displaystyle \ga_n:=\inf_{|x|=r_n}|x|^{N-2}u(x)\uparrow\infty\quad\text{as }n\to\infty.
$$
Therefore, $u$ is bounded from below in $B_1\setminus B_{r_n}$ by the (radial) solution $u_n$ of 
\bel{K14}\BA{lll}
-\Gd u+u^p-M|\nabla u|^q=0\quad&\text{in } B_1\setminus\overline B_{r_n}\\
\phantom{-\Gd +u^p-M|\nabla u|^q}
v=0\quad&\text{on } \prt B_1\\
\phantom{-\Gd +u^p-M|\nabla u|^q}
u=\ga_nr_n^{2-N}\quad&\text{on } \prt B_{r_n}.
\EA\ee
The sequence $\{u_n\}$ may not be monotone, but $u_n\geq v_n$ where $v_n$ has been defined in $(\ref{K10})$. Since $\{u_n\}$ is eventually locally bounded in 
$\overline B_1\setminus\{0\}$ by \rprop{Oss} and standard regularity results (see e.g. \cite{GT}), up to a subsquence, it converges locally uniformly in $\overline B_1\setminus\{0\}$ to a radial function $u_\infty$ which satisfies 
\bel{K15bis}\BA{lll}
-\Gd u+u^p-M|\nabla u|^q=0\quad&\text{in } B_1\setminus\{0\}\\
\phantom{-\Gd +u^p-M|\nabla u|^q}
v=0\quad&\text{on } \prt B_1,
\EA\ee
and 
\bel{K16}\BA{lll}
\displaystyle\liminf_{x\to 0}|x|^\ga u_\infty(x)\geq \liminf_{x\to 0}|x|^\ga v_\infty(x)=x_0.
\EA\ee
By \rth{subcrit-th}-(4) we have that 
\bel{K17}\BA{lll}
\displaystyle\liminf_{x\to 0}|x|^\ga u(x)\geq \displaystyle\lim_{x\to 0}|x|^\ga u_\infty(x)=x_{_M}. 
\EA\ee
The upper estimate is obtained as follows. By \rprop{harnak} the function $u$ is bounded from above in $B_1\setminus B_{r_n}$ by 
$\displaystyle u'_n+\max_{|y|=1}u(y)$ where $u'_n$ is the solution of 
\bel{K14bis}\BA{lll}
-\Gd u+u^p-M|\nabla u|^q=0\quad&\text{in } B_1\setminus\overline B_{r_n}\\
\phantom{-\Gd +u^p-M|\nabla u|^q}
v=0\quad&\text{on } \prt B_1\\
\phantom{-\Gd +u^p-M|\nabla u|^q}
u=c_9\ga_nr_n^{2-N}\quad&\text{on } \prt B_{r_n}.
\EA\ee
For the same reason as above there exists a subsequence, $\{u'_{n_\ell}\}$ which converges locally uniformly in $\overline B_1\setminus\{0\}$ to a radial solution $u'_\infty$ 
of $(\ref{K15bis})$. By \rth{subcrit-th}-(4) we have that $u'_\infty=u_\infty$. Then $\displaystyle\limsup_{x\to 0}|x|^\ga u(x)\leq \displaystyle\lim_{x\to 0}|x|^\ga u'_\infty(x)=x_{_M,}$ which ends the proof.$\phantom{azf}$\qeda\medskip

\nind{\it Proof of \rth{19M}} It is a direct consequence of \rlemma{12adapt}, \rlemma{13} and \rlemma{14}.\qeda

\medskip

When $q>\frac{2p}{p+1}$ and $|x|^{N-2}u(x)$ is not bounded, Harnack inequality may not hold. However we still have a dichotomy for the possible behaviour which extends \rth{19M}. 

\bth{20M} Let $\Gw\subset\BBR^N$ ($N\geq 3$) be an open set containing $0$, $M>0$, $1<p<\frac{N}{N-2}$ and $\frac{2p}{p+1}<q<\frac{N}{N-1}$. If $u$ is a nonnegative solution of $(\ref{Z1})$ in $\Gw\setminus\{0\}$, then the following dichotomy holds:\\
\nind 1- either there exists $k\geq 0$ such that $r^{N-2}u(r,.)$ converges to $k$ in measure on $S^{N-1}$ as $r\to 0$,\\
\nind 2- or  $\displaystyle\lim_{x\to 0}|x|^\gg u(x)=X_{_M}$.
\es
\Proof We recall Richard-V\'eron's isotropy theorem  \cite[Theorem 1.1]{RiVe}:\\
 {\it Let $g:\BBR_+\mapsto\BBR_+$ be a continuous nondecreasing function satisfying
\bel{K15}\BA{lll}
\myint{0}{1}g(r^{2-N})r^{N-1}dr<\infty.
\EA\ee
If $u\in C^1(\Gw\setminus\{0\})$ is a nonnegative function satisfying
\bel{K16bis}\BA{lll}
\Gd u\leq g(u)+f\quad\text{in }\,\Gw\setminus\{0\},
\EA\ee
where $f\in L^1_{loc}(\Gw\setminus\{0\})$ is a nonnegative radial function, then we have the following:\smallskip

\nind (i) either $r^{N-2}u(r,.)$ converges in measure on $S^{N-1}$ to some $k\geq 0$ as $r\to 0$,\\
(ii) or}
\bel{K17bis}\BA{lll}
\displaystyle \lim_{x\to 0}|x|^{N-2}u(x)=\infty.
\EA\ee
Since $p<\frac{N}{N-2}$, assumption $(\ref{K15})$ is satisfied with $g(r)=r^p$ and equation $(\ref{K16bis})$ with $f=0$. Then either the statement 1 holds, or $(\ref{K17bis})$ holds. If it is the case, then for any $k>0$,  $u$ is bounded 
from below in $B_1\setminus\{0\}$ by the solution $u_k$ of $(\ref{Z1})$ in $B_1\setminus\{0\}$ vanishing on $\prt B_1$ and satisfying $(\ref{K5})$. Such a solution exists and is unique by \rth{18M}. Letting $k\to\infty$, this implies
\bel{K18}\BA{lll}
\displaystyle \liminf_{x\to 0}|x|^{\gg}u(x)\geq X_{_M}.
\EA\ee
Next, we denote by $\gu$ the solution of $(\ref{Z1})$ on $(0,\infty)$ (hence $\gu$ is  1-dimensional) satisfying 
$$\displaystyle\lim_{r\to 0}r^\gg \gu(r)=X_{_M}\quad\text{and }\,\lim_{r\to \infty}r^\ga \gu(r)=x_0.
$$
Its existence is proved in \rth{T9}, \rth{T10}. It is decreasing. For $\ge>0$ the function $r\mapsto \gu_{\ge}(x)=\gu(|x|-\ge)$ satisfies $\CL^M_{p,q}( \gu_{\ge})\geq 0$ in 
$\overline B_\ge^c$ and $\gu_{\ge}(x)\to \infty$ when $|x|\downarrow \ge$. If $c=\max u\lfloor_{\prt B_\gd}$ for some $\gd>\ge>0$ such that $\overline B_\gd\subset \Gw$, then $\gu_{\ge}+c$ is a supersolution of $(\ref{Z1})$ in $B_\gd\setminus\overline B_\ge$ which is larger than $u$ for $|x|=\ge$ and $|x|=\gd$. Hence 
$u\leq gu_{\ge}+c$ in $B_\gd\setminus\overline B_\ge$. Letting $\ge\to 0$ yields $u(x)\leq \gu(x)+c$ for $0<|x|\leq \gd$ and finally
\bel{K19}\BA{lll}
\displaystyle \limsup_{x\to 0}|x|^{\gg}u(x)\leq X_{_M}.
\EA\ee
Combining $(\ref{K18})$ and $(\ref{K19})$ we obtain $\displaystyle\lim_{x\to 0}|x|^\gg u(x)=X_{_M}$.\qeda \medskip

\nind\Remark We conjecture that the convergence in \rth{20M}-1 holds in the strong sense. 
\medskip

\nind\Remark Most of the results of this section can be extended to the case $N=2$. The subcritical case corresponds then to $p>1$ and $1<q<2$. The main change is that $|x|^{2-N}$ has to be replaced by 
$-\ln|x|$. 
\subsection{The case $p\geq \frac{N}{N-2}$ and $q=\frac{2p}{p+1}$}
The cases that we consider are $q=\frac{2p}{p+1}$, $p\geq\frac{N}{N-2}$ and $M>0$. We recall that the stationary equation $(\ref{Z3})$ admits two positive constant solutions $x_{_{1,M}}<x_{_{2,M}}$ if $p>\frac{N}{N-2}$ and $M>m^*$, and only one denoted by $x_{_{M}}$ if  $p=\frac{N}{N-2}$ and $M>0$ or if $p>\frac{N}{N-2}$ and $M=m^*$. The following result is an improvement of \rprop{Oss}
 \blemma {sharp}Let $\Gw\subset\BBR^N$, $N\geq 3$, be a bounded domain containing $0$ such that $\overline B_R\subset\Gw$ and $p\geq\frac{N}{N-2}$. \\
 1- If  $p>\frac{N}{N-2}$ and $M>m^*$, then any positive solution $u$ of  $(\ref{Z7})$ in $\Gw\setminus\{0\}$ satisfies 
 \bel{K28}\BA{lll}
u(x)\leq x_{_{2,M}}|x|^{-\ga}+\displaystyle\sup_{|z|=R}u(z),
\EA\ee
 2- If $p>\frac{N}{N-2}$ and $M=m*$, or $p=\frac{N}{N-2}$ and $M>0$, the same inequality holds with $x_{_{2,M}}$ replaced by $x_{m^*}$ and $x_{_{M}}$  respectively.
 \es
 \Proof We assume first that $M>m^*$.\\
  1- {\it Construction of the maximal solution}. We claim that $x\mapsto x_{_{2,M}}|x|^{-\ga}$ is the maximal solution of $(\ref{Z7})$ in $\BBR^N\setminus\{0\}$. For $a>0$ 
 we set $\phi_a(s)=as^{\ga}$. Then 
 $$\BA {lll}
 \tilde L\phi_a(s):=-\phi''_a+\phi^p_a-M|\phi'_a|^{\frac{2p}{p+1}}=as^{-\ga p}\left(a^{p-1}-\ga^qMa^{q-1}-\ga(\ga+1)\right).
 \EA$$
 Taking $a$ large enough we obtain that $\phi_a$ is a supersolution in $(0,\infty)$. We set $\Gf_{a,\ge}(x)=\gf_a(x_1-\ge)$ for $x_1>\ge>0$ and as in the proof of \rth{20M} we deduce that the function $\tilde \Gf_{a,\ge}$ defined by 
 $$\tilde \Gf_{a,\ge}=\inf\left\{\CR[\Gf_{a,\ge}]:\CR\in SO(N)\right\}\qquad \text{($SO(N)$ is the group of rotations in $\BBR^N$)},
 $$
 is a positive and radial supersolution of $(\ref{Z7})$ in $\BBR^N\setminus\overline B_\ge$ which tends to infinity on $\prt B_\ge$. It dominates any positive solution of $(\ref{Z7})$ in $\BBR^N\setminus\overline B_\ge$. Next we set 
 $$\tilde \Psi_{\ge}(x)=\sup\left\{x_{_{2,M}}|x-z|^{-\ga}:|z|<\ge\right\}.
 $$
 It is a subsolution of $(\ref{Z7})$ in $\BBR^N\setminus\overline B_\ge$ dominated by $\tilde \Gf_{a,\ge}$. Since the supremum is achieved for 
 $z=\ge \frac{x}{|x|}$, the function $\tilde \Psi_{\ge}$ is radial and positive. By \cite[Theorem 1.4.5]{Vebook} there exists a solution $U_\ge$ of $(\ref{Z7})$ in $\BBR^N\setminus\overline B_\ge$ such that 
 $$\tilde \Psi_{\ge}\leq U_\ge\leq \tilde \Gf_{a,\ge}\quad\text{in }\,\BBR^N\setminus\overline B_\ge.
 $$
The function $U_\ge$ is positive and radial. Since any positive solution $u$ in $\BBR^N\setminus\{0\}$ is dominated by $\tilde \Psi_{\ge}$, the function $U_\ge$ is larger than $u$ in $\BBR^N\setminus\overline B_\ge$. This implies the relation, valid for any $\ell>0$, 
 \bel{K29}\BA{lll}T_\ell[U_\ge](x):=\ell^{\ga}U_\ge(\ell x)=U_{\ell^{-1}\ge}.
\EA\ee
When $\ge\downarrow 0$ the sequence $\{U_\ge\}$ decreases and converges to a positive radial solution $U_0$ of $(\ref{Z7})$ in $\BBR^N\setminus\{0\}$ which 
dominates any other positive solution. Hence $U_0$ is the maximal solution in $\BBR^N\setminus\{0\}$. Letting $\ge\downarrow 0$ in $(\ref{K29})$ we infer that 
$T_\ell[U_0]=U_0$ for any $\ell>0$. Hence $U_0$ is self-similar. Since it is radial and larger than any other positive solution, we deduce that 
 \bel{K30}\BA{lll}U_0(x)=x_{_{2,M}}|x|^{-\ga}\quad\text{for all }|x|>0.
\EA\ee
\nind 2- {\it End of the proof}. If $u$ is any positive solution in $\Gw\setminus\{0\}$, then $U_\ge+\displaystyle\sup_{|z|=R}u(z)$ is a supersolution larger than $u$ in $B_R\setminus B_\ge$. Letting $\ge\downarrow 0$ yields the result.
The proof in the other case is similar.\qeda

\bth{supercrit-p} Let $\Gw\subset\BBR^N$, $N\geq 3$, be a bounded domain containing $0$, $p>\frac{N}{N-2}$ and $M>m^*$. If $u$ is a positive solution of 
$(\ref{Z7})$ in $\Gw\setminus\{0\}$, there holds\smallskip

\nind 1- If 
\bel{K23}\BA{lll}
\displaystyle \liminf_{x\to 0}|x|^{\ga}u(x)=0,
\EA\ee
then $u$ can be extended as a $C^2$ solution of $(\ref{Z7})$ in $\Gw$.\smallskip

\nind 2- If
\bel{K24}\BA{lll}
\displaystyle \limsup_{x\to 0}|x|^{\ga}u(x)=x_{_{2,M}},
\EA\ee
then 
\bel{K25}\BA{lll}
\displaystyle \lim_{x\to 0}|x|^{\ga}u(x)=x_{_{2,M}},
\EA\ee
\smallskip

\nind 3- If
\bel{K26}\BA{lll}
\displaystyle \liminf_{x\to 0}|x|^{\ga}u(x)=x_{_{1,M}}\quad\text{or }\;\limsup_{x\to 0}|x|^{\ga}u(x)=x_{_{1,M}},
\EA\ee
then there exists a sequence $\{r_n\}\subset\BBR^+_*$ converging to $0$ such that 
\bel{K27}\BA{lll}
\displaystyle \lim_{r_n\to 0}r_n^{\ga}u(r_n,s)=x_{_{1,M}}\quad\text{uniformly on }S^{N-1}.
\EA\ee

\nind 4- Any one of the two following situations never occur
\bel{K31}\BA{lll}
(i)\qquad\displaystyle &\displaystyle 0<\limsup_{x\to 0}|x|^{\ga}u(x)<x_{_{1,M}}\qquad\qquad\qquad\qquad\qquad\\[2mm]
(ii)\qquad\displaystyle &\displaystyle x_{_{1,M}}<\liminf_{x\to 0}|x|^{\ga}u(x)<x_{_{2,M}}.
\EA\ee
\es
\Proof 1- If relation $(\ref{K23})$ holds, it follows by Harnack inequality proved in \rprop{harnak}, that there exists a sequence $\{r_n\}$ converging to $0$ as $n\to\infty$ such that 
\bel{K32}\BA{lll}
\displaystyle \lim_{r_n\to 0}r_n^{\ga}u(r_n,s)=0\quad\text{uniformly on }S^{N-1}.
\EA\ee
For any $\ge>0$ and $m=\sup\{ u(z):|z|=R\}$, the function 
$x\mapsto \ge |x|^\ga+m$ is a supersolution of $(\ref{Z7})$ in $B_R\setminus\{0\}$ which is larger than $u$ near $x=0$ and on $\prt B_R$. Hence 
$u(x)\leq \ge |x|^\ga+m$. Letting $\ge\to 0$ implies $u\leq m$, and the result follows by standard regularity.\\
\nind 2- If $(\ref{K24})$ holds there exists a sequence $\{r_n\}$ converging to $0$ such that
$$r_n^\ga\max\left\{u(r_n,s):s\in S^{N-1}\right\}=r_n^\ga u(r_n,s_n)\to x_{_{2,M}}\quad\text{as }r_n\to 0.
$$
Furthermore, we can assume that $s_n\to s^*$ when $n\to\infty$. Using \rlemma{sharp} there exist a nondecreasing sequence $\{w_{1,n}\}$ converging to $x_{_{2,M}}$ and a bounded sequence $\{w_{2,n}\}$ such that $r_n^\ga u(r_n,s_n)=w_{1,n}+r_n^\ga w_{2,n}$. We set $w(t,s)=r^\ga u(r,s)$ with $t=\ln r$, then there holds,
\bel{K33}\BA{lll} w_{tt}+Lw_t-\ga Kw+\Gd' w-w^p+M\left(\left(w_t-\ga w\right)^2+|\nabla' w|^2\right)^{\frac{p}{p+1}}=0,
\EA\ee
 on $\BBR_-\ti S^{N-1}$. By standard regularity estimates and Ascoli-Arzela theorem  there exist a subsequence $\{t_{n_j}\}$ of $\{t_n\}=\{\ln r_n\}$ and a 
 nonnegative $C^2$ function $W$ such that $w(t_n+t,s)$ converges to $W$ in the $C^2$ topology of $[-a,a]\ti S^{N-1}$, for any $a>0$, and $W$ is a solution of $(\ref{K33})$ in $\BBR\ti S^{N-1}$.  Furthermore $W(0,s^*)=x_{_{2,M}}$. By \rlemma{sharp}, $x_{_{2,M}}$ is the maximal solution of $(\ref{K33})$ in $\BBR\ti S^{N-1}$, it then follows from the strong maximum principle that $W=x_{_{2,M}}$ and $w(t_n+t,s) \to x_{_{2,M}}$ uniformly in $[-a,a]\ti S^{N-1}$. Let $\ge>0$, then there exists $n_\ge\in\BBN$ such that for any $n\geq n_\ge$, we have that $u(r_n,s)\geq r_n^{-\ga} (x_{_{2,M}}-\ge)$ for any $s\in S^{N-1}$. Since 
 $r\mapsto r^{-\ga} (x_{_{2,M}}-\ge)$ is a subsolution of $(\ref{Z7})$, it follows that for $m>n\geq n_\ge$, one has
 $$u(r,s)\geq r^{-\ga} (x_{_{2,M}}-\ge)\quad\text{for all }(r,s)\in [r_m,r_n]\ti S^{N-1}.
 $$ 
 Letting $r_m\to 0$ yields 
 \bel{K34}\BA{lll} u(x)\geq |x|^{-\ga} (x_{_{2,M}}-\ge)\qquad\text{for all } x\in B_{r_{n_\ge}}\setminus\{0\}.
\EA\ee
Since $\ge$ is arbitrary we infer that 
 \bel{K35}\BA{lll} \displaystyle \liminf_{x\to 0} |x|^\ga u(x)\geq x_{_{2,M}}.
\EA\ee
 By assumption $\displaystyle \limsup_{x\to 0} |x|^\ga u(x)\geq x_{_{2,M}}$; then $(\ref{K25})$ holds.\\
\nind 3- Let us assume that the first condition in $(\ref{K26})$ holds. If the function 
\bel {K35-1}\gt\mapsto \displaystyle w(\gt,s_\gt):=\min_{s\in S^{N-1}}w(\gt,s)\}\ee
is asymptotically monotone, nonincreasing or nondecreasing, then either $w(\gt,s_\gt)\uparrow x_{_{1,M}}$ in the first case, or $w(\gt,s_\gt)\downarrow x_{_{1,M}}$
 in the second case. Using again the uniform $C^{2,\ga}$ estimate and Ascoli-Arzela theorem we have that there exists a sequence $\{\gt_n\}$ converging to 
 $-\infty$ such that $w(t+\gt_n,s)$ converges in the $C^2$-topology of $[-a,a]\ti S^{N-1}$, for any $a>0$, to a positive solution $W$ of $(\ref{K33})$ in $\BBR\ti S^{N-1}$ such that $W(t,s)\geq x_{_{1,M}}$ and $\displaystyle \min_{s\in S^{N-1}}W(0,s)=W(0,s^*)=x_{_{1,M}}$ for some $s^*\in S^{N-1}$ in the first case.  By the strong maximum principle to $w$ and to $x_{_{1,M}}$ which are ordered solutions of $(\ref{K33})$ in $\BBR\ti S^{N-1}$ we infer that 
 $W\equiv x_{_{1,M}}$, hence $w(t,s)$ converges to $x_{_{1,M}}$ uniformly on $S^{N-1}$ when $t\to-\infty$. In the second case we obtain that the limit function $W$ satisfies $W(t,s)\geq x_{_{1,M}}$ and $\displaystyle \min_{s\in S^{N-1}}W(0,s)=W(0,s^*)=x_{_{1,M}}$. This implies again $W\equiv x_{_{1,M}}$ by the strong maximum principle. Finally we do not suppose that the function $w(\gt,s_\gt)$ defined in 
 $(\ref{K35-1})$ is monotone. By the definition of the liminf, there exist  sequences $\{t_n\}$ tending to $-\infty$ and $\{s_n\}\subset S^{N-1}$ such that 
 $$w(t_n,s_n)=\inf\{w(t,s):t\leq t_n,s\in S^{N-1}\}\uparrow x_{_{1,M}}\,\text{ as }\,n\to\infty.
 $$
 Using again Ascoli-Arzela theorem we deduce that, up to a subsequence $\{t_{n_k}\}$, $w(t+t_{n_k},s)$ converges in the $C^2$-topology of $[-a,a]\ti S^{N-1}$, for any $a>0$, to a positive solution $W$ of $(\ref{K33})$ in $\BBR\ti S^{N-1}$ and $W\geq x_{_{1,M}}$ and $W(0,s^*)=x_{_{1,M}}$. Hence $W\equiv x_{_{1,M}}$. \\
 \nind The proof of $(\ref{K27})$ under the second condition of $(\ref{K26})$ is similar.
\smallskip

\nind 4- Let $(\ref{K31})$-(i) be satisfied and $w$ be defined as in the previous steps. Then as in Step 3, $w(t_n,s)$ converges locally uniformly to a positive solution $W$ of $(\ref{K33})$ defined on $\BBR_-\ti S^{N-1}$, $w(t_n,s_n)\to W(0,s^*)$,  $W(0,s^*)$ is a local maximum of $W$ and it is smaller than $x_{_{1,M}}$. Hence 
$W_t(0,s^*)=|\nabla 'W(0,s^*)|=0$ and $W_{tt}(0,s^*)\leq 0$ and $\Gd'W(0,s^*)\leq 0$. 
Then
$$W^p(0,s^*)-M\ga^q W^\frac{2p}{p+1}(0,s^*)+\ga KW(0,s^*)\leq 0.
$$
This contradicts the fact that $\CP_{_M}$ defined in $(\ref{Z5})$ is positive on the interval $(0,x_{_{1,M}})$. \\
Similarly, if $(\ref{K31})$-(ii), we obtain that the limit function $W$ and the limit point $s^*$ where $W(0,s^*)$ is a local minimum of $W$ satisfies
$x_{_{1,M}}<W(0,s^*)<x_{_{2,M}}$ and 
$$W^p(0,s^*)-M\ga^q W^\frac{2p}{p+1}(0,s^*)+\ga KW(0,s^*)\geq 0.
$$
which is not compatible with the fact that $\CP_{_M}$ is negative on $(x_{_{1,M}},x_{_{2,M}})$.\qeda\medskip

\nind\Remark We conjecture that the following stronger form of \rth{supercrit-p} holds:\smallskip

\nind{\it 1- Either $u$ can be extended as a $C^2$ solution of $(\ref{Z7})$ in $\Gw$,\smallskip

\nind 2- or $\displaystyle\lim_{x\to 0}|x|^\ga u(x)= x_{_{2,M}}$,\smallskip

\nind 3- or $\displaystyle\lim_{x\to 0}|x|^\ga u(x)= x_{_{1,M}}$.
} \medskip

In the case $M=m^*$, we prove the following.
\bth{supercrit-m^*} Let $\Gw\subset\BBR^N$, $N\geq 3$, be a bounded domain containing $0$, $p>\frac{N}{N-2}$ and $M=m^*$. If $u$ is a positive solution of 
$(\ref{Z7})$ in $\Gw\setminus\{0\}$, there we have the following:\smallskip

\nind 1- either $u$ can be extended as a $C^2$ solution in $\Gw$,\smallskip

\nind 2- or there exists a sequence $\{r_n\}$ converging to $0$ such that $r_n^\ga u(r_n,s)$ converges to $x_{m^*}$ uniformly on $S^{N-1}$.
\es
\Proof If $u$ satisfies $(\ref{K23})$, then the singularity of $u$ at zero is removable since the function $P_{_M}$ is positive on $(0,m^*)$ and on $(m^*,\infty)$, see the argument in the proof of \rth{supercrit-p}-(1). Thus we are left with the case 
\bel{K36}\BA{lll}
\displaystyle \liminf_{x\to 0}|x|^{\ga}u(x)>0.
\EA\ee
If 
\bel{K36bis}\BA{lll}
\displaystyle \limsup_{x\to 0}|x|^{\ga}u(x)=x_{m^*},
\EA\ee
then, as in the proof of \rth{supercrit-p}-(3) we deduce that there exists a sequence $\{r_n\}$ converging to $0$ such that $r_n^\ga u(r_n,s)$ converges to 
$x_{m^*}$ uniformly on $S^{N-1}$. If there exists $m\neq {m^*}$ such that 
\bel{K37}\BA{lll}
\displaystyle \limsup_{x\to 0}|x|^{\ga}u(x)=m,
\EA\ee
then there exists a sequence $\{t_n\}$ converging to $-\infty$ and $\{s_n\}\subset S^{N-1}$ such that $w(t+t_n,.) $ converges in the $C^2$ topology 
of $[-a,a]\ti S^{N-1}$ for any $a>$ to a solution $W$ of $(\ref{K33})$. Furthermore $w_t(t_n,s_n)\to 0$ and $\displaystyle \liminf_{t_n\to-\infty} w_{tt}(t_n,s_n)\leq 0$. Since $\nabla 'w(t_n,s_n)=0$ and $\Gd 'w(t_n,s_n)\leq 0$, one has that
\bel{K37-1}-\ga K m-m^p+M(\ga m)^\frac{2p}{p+1}=-m\tilde P_{m^*}(m^{\frac{p-1}{p+1}})=-\Gd 'W(0,s^*) -w_{tt}(0,s^*)\geq 0.
\ee
Since $\tilde\CP_{m^*}\geq 0$ and vanishes only at $m^*$, it implies $m=m^*$, contradiction. The proof of the 
uniform convergence of $w(t_n,.)$ to $m^*$ follows from the strong maximum principle since $W$ is a positive solution $(\ref{K33})$ as in\rth{supercrit-p}.
\qeda\medskip

\nind\Remark We conjecture that assertion (2) holds under the form
$$\lim_{x\to 0}|x|^\ga u(x)=x_{m^*}.
$$

Finally we have the following result dealing with the case $p=\frac{N}{N-2}$ and $M>0$ where there exists a unique and explicit positive constant solution $x_{_M}$ to $(\ref{Z3})$.
\bth{supercrit-weak} Let $\Gw\subset\BBR^N$, $N\geq 3$, be a bounded domain containing $0$, $p=\frac{N}{N-2}$ and $M>0$. If $u$ is a positive solution of 
$(\ref{Z7})$ in $\Gw\setminus\{0\}$ which satisfies
\bel{K38}\BA{lll}
\displaystyle \limsup_{x\to 0}|x|^{\ga}u(x)=x_{_M},
\EA\ee
then
\bel{K39}\BA{lll}
\displaystyle \lim_{x\to 0}|x|^{\ga}u(x)=x_{_M},
\EA\ee
\es
\Proof Since the function $\CP_{_M}$ is negative on $(0,x_{_M})$ and positive on $(x_{_M},\infty)$, for any $\ge>0$ the function $x\mapsto (x_{_M}-\ge)|x|^{-\ga}$ is a subsolution of $(\ref{Z7})$. The proof follows as in the proof of \rth{supercrit-p}-(2).\qeda

\medskip

\nind\Remark We conjecture that the following dichotomy occurs:  {\it  if $u$ is a positive solution of $(\ref{Z7})$ in $\Gw\setminus\{0\}$ unbounded near $0$, then,\smallskip

\nind 1- either $(\ref{K39})$ holds,\smallskip

\nind 2- or $(\ref{Z11})$-(i) holds.}
\mysection{Behaviour at infinity of non-radial solutions}
In this section we present some results dealing with the asymptotic behaviour of solutions which extend to the non-radial case what has already been proved in the radial one. The results are more complete if there exists only one possible behaviour for radial positive solutions; they have to be compared with what was obtained in \cite{Veasym} when $M=0$. 
\subsection{The case $q=\frac{2p}{p+1}$}
\bth{ploplus} Let $p>1$, $M>0$ and $u$ be a positive solution of $(\ref{Z7})$ in $\BBR^N\setminus B_R$ with $N\geq 1$. \smallskip

\nind 1- If $N=1,2$ and $p>1$, or $N\geq 3$ and $1<p<\frac{N}{N-2}$, then
\bel{K40}\BA{lll}
\displaystyle \lim_{|x|\to \infty}|x|^{\ga}u(x)=x_{_M}.
\EA\ee

\nind 2- If $N\geq 3$, $p=\frac{N}{N-2}$, then
\bel{K41}\BA{lll}
\displaystyle \lim_{|x|\to \infty}|x|^{N-2}\left(\ln |x|\right)^{\frac{N-2}{2}}u(x)=\left(\myfrac{N-2}{\sqrt 2}\right)^{N-2}.
\EA\ee

\nind 3- If $N\geq 3$, $p>\frac{N}{N-2}$ and $M<m^*$, then there exists $k>0$ such that
\bel{K41+}\BA{lll}
\displaystyle \lim_{|x|\to \infty}|x|^{N-2}u(x)=k.
\EA\ee
\es
\Proof The method of the proof is firstly to construct two positive radial solutions  $u_j$, $j=1,2$ of $(\ref{Z7})$ in $B_R^c$ such that $u_1\leq u\leq u_2$, and to use 
\rprop{hardc}. The solution $v$ of $\CL_pv=0$ in $B_R^c$ with value $\min u\lfloor_{\prt B_R}$ for $|x|=R$ is a subsolution smaller than $u$. For cases 1 and 2, we can take for supersolution the function $u_{X_{_M}}+\max u\lfloor_{\prt B_R}$. Therefore there exist two positive radial solutions $u_1$ and $u_2$ of $(\ref{Z7})$ in $B_R^c$ with respective value $\min u\lfloor_{\prt B_R}$ and $\max u\lfloor_{\prt B_R}$ on $\prt B_R$ and such that
\bel{K42}\BA{lll}
v(x)\leq u_1(x)\leq u(x)\leq u_2(x)\leq u_{X_{_M}}(x)+\max u\lfloor_{\prt B_R}\quad\text{for }|x|\geq R,
\EA\ee
by \rth{subcrit-th} and \rth{crit-th}. Since $u_1$ and $u_2$ satisfy either $(\ref{K40})$ or $(\ref{K41})$ in cases 1 and 2 respectively, $u$ shares this behaviour.\smallskip

\nind In case 3 with $p>\frac{N}{N-2}$ the function $v$ satisfies the same behaviour $(\ref{K41+})$ up to the constant $c>0$ which is not fixed. By \rth{supcrit}-(3), $u_1$ and $u_2$ satisfy  $(\ref{K41+})$ with two different constants $0<c_1\leq c_2$. In order to prove that $(\ref{K41+})$ holds for some $c\in [c_1,c_2]$ we use the method introduced in 
\cite{Veasym}. We set $u(r,s)=r^{-\gn}w(t,s)$ with $\gn=N-2$ and $t=\ln r$, then $w$ satisfies 
\bel{K43}\BA{lll}
w_{tt}-\gn w_t+\Gd' w+e^{(N-p(N-2))t}\left(\gn^q M\left((w_t-\gn w)^2+|\nabla 'w|^2\right)^\frac q2-w^p\right)=0
\EA\ee
in $[0,\infty)\ti S^{N-1}$. Since $w$ and $\left(\gn^q M\left((w_t-\gn w)^2+|\nabla 'w|^2\right)^\frac q2-w^p\right)$ are bounded, it follows from \cite[Proposition 4.1]{BiRa} that there exists $c\geq 0$ such that $w(t,.)\to c$ uniformly on $S^{N-1}$ when $t\to\infty$. This ends the proof.\qeda

\bth{ploplusM=m*} Let $N\geq 3$, $p>\frac N{N-2}$  and $M=m^*$. If $u$ is a positive solution of $(\ref{Z7})$ in $\BBR^N\setminus B_R$, we have the following alternative,\smallskip

\nind 1- either there exists a sequence $\{r_n\}$ converging to infinity such that 
\bel{K40-1}\BA{lll}
\displaystyle \lim_{r_n\to \infty}r_n^{\ga}u(r_n,s)=x_{m^*},
\EA\ee
uniformly on $S^{N-1}$.\smallskip

\nind 2- or there exists $k>0$ such that
\bel{K41-1}\BA{lll}
\displaystyle \lim_{|x|\to \infty}|x|^{N-2}u(x)=k.
\EA\ee
\es
\Proof We can assume that $u$ is continuous in $B_R^c$.  By \rth{supcrit}-(2), $u$ is bounded from above by the function $u_R$ where 
$u_R$ is a radial soluion of $(\ref{Z7})$ in $\overline B_R^c$ which tends to $\infty$ when $r\downarrow R$ and satisfies $r^\ga u_R(r)\downarrow x_{m^*}$ when $r\to\infty$. Hence 
\bel{K41-11}m:=\limsup_{|x|\to \infty}|x|^{\ga}u(x)\leq \sup_{|x|\geq R}|x|^{\ga}u(x)\leq x_{m^*}.
\ee
{\it We claim that either $m=0$ or $m=x_{m^*}$} As in the proof of \rth{supercrit-m^*}, there exists a sequence  $\{t_n\}$ tending to $\infty$ and $\{s_n\}\subset S^{N-1}$ converging to $s^*$ such that $w_n(t,.):=w(t+t_n,.)$ converges in the $C^2$-topology of $[-a,a]\ti S^{N-1}$ for any $a>0$ to a solution $W$ of $(\ref{K33})$ in $\BBR\ti S^{N-1}$. Furthermore  $W$ achieves its maximal value $m$ at $(0,s^*)$, hence $W_t(0,s^*)= 0$, $\nabla'W(0,s^*)=0$, $\Gd'W(0,s^*)\leq 0$ and $W_{tt}(0,s^*)\leq 0$. Therefore 
$$-m\tilde \CP_{m^*}(m^{\frac{p-1}{p+1}})=-\ga Km-m^p+M(\ga m)^\frac{2p}{p+1}\geq 0.
$$
Since $(\ref{K37-1})$ holds this implies that either $m=0$ or $m=x_{m^*}$. If $m=x_{m^*}$ it follows by the strong maximum principle, as in the proof of \rth{supercrit-p}, that $w(t_n,s)\to x_{m^*}$ as $n\to\infty$, uniformly on $S^{N-1}$.

\nind If $m=0$, then
$$\displaystyle \lim_{|x|\to\infty}|x|^\ga u(x)=0.
$$
For any $a<x_{m^*}$ and $\gr>R$ such that $u\lfloor_{\prt B_\gr}\leq a$, we consider the problem
\bel{K41-12}\BA{lll}\displaystyle 
-v_{rr}-\frac{N-1}{r}v_r+v^p-M|v_r|^{\frac{2p}{p+1}}=0\quad\text{in }\,(\gr,\infty)\\
\phantom{--;,\frac{N-1}{r}v_r+v^p-M|v_r|^{\frac{2p}{p+1}}}
v(\gr)=\gr^{-\ga}a
\EA\ee
Since the solution of 
\bel{K41-121}\BA{lll}\displaystyle 
-\gu_{rr}-\frac{N-1}{r}\gu_r+\gu^p=0\quad\text{in }\,(\gr,\infty)\\
\phantom{--;,\frac{N-1}{r}\gu_r+\gu^p}
\gu(\gr)=\gr^{-\ga}a
\EA\ee
is a subsolution and $x_{m^*}|x|^{-\ga}$ a supersolution, the solution $v$ exists and it is unique. By the phase plane analysis of Figure 4, the function 
$\tilde\gu(t)= e^{\ga t}\tilde u(r)$ which initial value belongs to the region (F) converges to $0$ when $t\to\infty$. Since $(0,0)$ is a saddle point for the system
 $(\ref{W4})$ it follows that the corresponding trajectory is is the unstable one of this point.  The initial slope of this curve is $N-2$. By \rlemma{eigen} it follows that there exists $\ell>0$ such that 
 \bel{K41-122}\BA{lll}\displaystyle 
\lim_{r\to\infty}r^{N-2}v(r)=\ell. 
\EA\ee
Consequently $|x|^{N-2}u(x)$ is bounded, the proof of \rth{ploplus}-3 applies and deduce from $(\ref{K41-11})$ there exists $c>0$ such that $|x|^{N-2}u(x)\to c$ when $|x|\to\infty$.\qeda\medskip

In the case $M>m^*$ the situation is even more complicated and the results are still incomplete.
\bth{ploplusM>m*} Let $N\geq 3$, $p>\frac N{N-2}$  and $M>m^*$ If $u$ is a positive solution of $(\ref{Z7})$ in $\BBR^N\setminus B_R$, we have the following,\smallskip

\nind 1- 
\bel{K40-2}\BA{lll}
\displaystyle \limsup_{|x|\to \infty}|x|^{\ga}u(x)=x_{_{2,M}}\Longrightarrow\lim_{|x|\to \infty}|x|^{\ga}u(x)=x_{_{2,M}},
\EA\ee

\nind 2- If $\displaystyle \liminf_{|x|\to \infty}|x|^{\ga}u(x)=x_{_{1,M}}$, there exists a sequence $\{r_n\}$ tending to $\infty$ such that
\bel{K40-2-}\BA{lll}
\displaystyle \lim_{r_n\to \infty}r_n^{\ga}u(r_n,s)=x_{_{1,M}},
\EA\ee
uniformly on $S^{N-1}$.\smallskip

\nind 3-  If $\displaystyle \liminf_{|x|\to \infty}|x|^{\ga}u(x)=0$, there exists $k>0$ such that
\bel{K41-2}\BA{lll}
\displaystyle \lim_{|x|\to \infty}|x|^{N-2}u(x)=k.
\EA\ee
\es
\Proof By \rth{supcrit}-(1) $u$ is bounded from above by the solution $u_R$ of $(\ref{Z7})$ in $\BBR^N\setminus B_R$ which tends to infinity as 
$|x|\downarrow R$ and satisfies
$$
\displaystyle \lim_{|x|\to \infty}|x|^{\ga}u_R(x)=x_{_{2,M}}.
$$
Hence 
$$\displaystyle \tilde m:=\liminf_{|x|\to \infty}|x|^{\ga}u(x)\leq  m:=\limsup_{|x|\to \infty}|x|^{\ga}u(x)\leq x_{_{2,M}}
$$
1- If $m=x_{_{2,M}}$ there exists a sequence  $\{t_n\}$ tending to $\infty$ and $\{s_n\}\subset S^{N-1}$ converging to $s^*$ such that $w_n(t,.):=w(t+t_n,.)$ converges in the $C^2$-topology of $[-a,a]\ti S^{N-1}$ for any $a>0$ to a solution $W$ of $(\ref{K33})$ in $\BBR\ti S^{N-1}$. 
Furthermore  $W$ achieves its maximal value $m$ at $(0,s^*)$, hence $W_t(0,s^*)= 0$, $\nabla'W(0,s^*)=0$, $\Gd'W(0,s^*)\leq 0$ and $W_{tt}(0,s^*)\leq 0$. Therefore 
$$- m\tilde \CP_{_M}( m^{\frac{p-1}{p+1}})=-\ga K m- m^p+M(\ga m)^\frac{2p}{p+1}\geq 0.
$$
This implies that either $x_{_{1,M}}\leq  m\leq x_{_{2,M}}$ or $m=0$. For the liminf the same analysis yields that 
$$-\tilde m \tilde \CP_{_M}(\tilde m^{\frac{p-1}{p+1}})=-\ga K\tilde m-\tilde m^p+M(\ga \tilde m)^\frac{2p}{p+1}\leq 0,
$$
hence either $0\leq  \tilde m\leq x_{_{1,M}}$ or $\tilde m\geq x_{_{2,M}}$. Note that in the latter case $(\ref{K40-2})$ holds. \\
2- If $m=x_{_{2,M}}$, then using the function $W$ as  in the proof of \rth{ploplusM=m*}, we infer by the strong maximum principle that there exists a sequence $\{r_n\}$ tending to infinity such that $r_n^\ga u(r_n,s)$ converges to $x_{_{2,M}}$. For any $\ge>0$ there exists $n_\ge>0$ such that for $n\geq n_\ge$ we have $r_n^\ga u(r_n,s)\geq x_{_{2,M}}-\ge$ for all 
$s\in S^{N-1}$. Since $\CP_{_M}(x_{_{2,M}}-\ge) \leq 0$ the function $x\mapsto (x_{_{2,M}}-\ge)|x|^{-\ga} $ is a subsolution of $(\ref{Z7})$. Then, for any 
$r_n>r_{n_\ge}$, $u(x)\geq (x_{_{2,M}}-\ge)|x|^{-\ga} $ in $\{x:r_{n_\ge}\leq |x|\leq r_n\}$. This implies
$$\displaystyle \liminf_{|x|\to \infty}|x|^{\ga}u(x)\geq x_{_{2,M}}-\ge.
$$
Since $\ge$ is arbitrary, this yields $(\ref{K40-2})$.\\
If $\tilde m=x_{1,M}$ then we proceed as in the case above and deduce that there exist a function $W\geq 0$ satisfying $(\ref{K33})$ in $\BBR\ti S^{N-1}$ and a sequence $\{t_n\}$ tending to infinity such that 
$w(t_n+t,s)$ converges in the $C^2$-topology of $[-a,a]\ti S^{N-1}$ to $W$ for any $a>0$. The function $W$ is larger or equal to $x_{_{1,M}}$ and coincides with $x_{_{1,M}}$ at $(0,s^*)$ for some $s^*\in S^{N-1}$. By the strong maximum principle we have that $W\equiv x_{_{1,M}}$. This implies assertion 2.\smallskip

\nind3- If $\displaystyle \liminf_{|x|\to \infty}|x|^{\ga}u(x)=0$, then we deduce by Harnack inequality that there exists a sequence $\{r_n\}$ tending to $\infty$ such that 
$$\lim_{r_n\to\infty}r_n^\ga u(r_n,s)=0\quad\text{uniformly on }S^{N-1}. 
$$
Then for any $\ge>0$ the function there exists $n_\ge\in\BBN$ such that for any $n\geq n_\ge$, $u(r_n,s)\leq \ge |r_n|^{-\ga}$ for all $s\in S^{N-1}$. 
The function $x\mapsto\ge|x|^{-\ga}$ is a supersolution of $(\ref{Z7})$. Since there exists a sequence $\{r_n\}$ tending to infinity such that 
$u(r_n,s)\leq \ge r_n^{-\ga}$ for all $s\in S^{N-1}$ for 
$n\geq n_\ge$ it follows by the comparison principle applied to the sequence of annuli $\{x:r_{n_\ge}\leq |x|\leq r_n\}$, that $x(x)\leq \ge|x|^{-\ga}$. Since the function $x\mapsto \ge |x|^{-\ga}$ is a supersolution of $(\ref{Z7})$, it follows by the comparison principle that $u(x)\leq\ge |x|^{-\ga} $ in the annuli $\{x:r_{n_\ge}\leq |x|\leq r_n\}$. Letting $n\to\infty$ yields 
$$u(x)\leq\ge |x|^{-\ga}\quad\forall x\in B_{r_{n_\ge}}\Longrightarrow \limsup_{|x|\to \infty}|x|^{\ga}u(x)\leq \ge.
$$
Since $\ge$ is arbitrary we infer that $|x|^{\ga}u(x)$ converges to $0$ when $|x| \to\infty$. By the phase plane analysis of section 2.4 (see\rth{supcrit}), as in the proof of \rth{ploplusM=m*}-2, we have that $|x|^{N-2}u(x)$  is bounded. Hence $(\ref{K41-2})$ follows as in the previous proof.{\hspace{10mm}\hfill $\phantom{\square}$}\qeda\medskip

\nind\Remark We conjecture that the results of \rth{ploplusM=m*} and \rth{ploplusM>m*} hold under the following forms: \smallskip

\nind For \rth{ploplusM=m*}\smallskip

{\it \nind 1- either  $\displaystyle \lim_{x\to \infty}|x|^{\ga}u(x)=x_{m^*},$
\smallskip

\nind 2- or there exists $k>0$ such that $\displaystyle \lim_{|x|\to \infty}|x|^{N-2}u(x)=k$
}\smallskip

\nind For \rth{ploplusM>m*}\smallskip

{\it 
\nind 1- either  $\displaystyle \lim_{x\to \infty}|x|^{\ga}u(x)=x_{_{2,M}}$,
\smallskip

\nind 2- or  $\displaystyle \lim_{x\to \infty}|x|^{\ga}u(x)=x_{_{1,M}}$,
\smallskip

\nind 3- or there exists $k>0$ such that $\displaystyle \lim_{|x|\to \infty}|x|^{N-2}u(x)=k$.
}

\subsection{The case $q\neq\frac{2p}{p+1}$}
The next results extend the asymptotic behaviour described in \rth{T8} and \rth{T9} to non radial solutions. The following statement shows that equation $(\ref{Z1})$ inherits the properties of the Emden-Fowler equation $\CL_pu=0$ if $\frac{2p}{p+1}<q<p$.

\bth{T8bis}Let $N\geq 1$, $M>0$, $\frac{2p}{p+1}<q<p$ and $u$ be a positive solution of $(\ref{Z1})$ in $B_R^c$. Then\smallskip

\nind 1- If $N=1,2$ and $p>1$, or $N\geq 3$ and $1<p<\frac{N}{N-2}$, then 
\bel{K44}\BA{lll}\displaystyle
\lim_{|x|\to\infty}|x|^\ga u(x)=x_0.
\EA\ee
\smallskip

\nind 2- If $N\geq 3$ and $p>\frac{N}{N-2}$, then 
\bel{K45}\BA{lll}\displaystyle
\lim_{|x|\to\infty}|x|^{N-2}u(x)=k>0.
\EA\ee
\smallskip

\nind 3- If $N\geq 3$ and $p=\frac{N}{N-2}$, then 
\bel{K46}\BA{lll}\displaystyle
\lim_{|x|\to\infty}\left(\ln |x|\right)^{\frac{N-2}{2}}|x|^{N-2}u(x)=\left(\myfrac{N-2}{\sqrt 2}\right)^{N-2}.
\EA\ee
\es
\Proof In the first case, the solution $v$ of $\CL_pv=0$ in $B_R^c$ with value $\min u\lfloor_{\prt B_R}$ on $\prt B_R$ is a subsolution smaller than $u$ (it is obtained by minimization), and it has the behaviour expressed by $(\ref{K44})$. By \rth{T12} there exists a global positive solution $\tilde u$ of $(\ref{Z1})$ in $\BBR^N\setminus\{0\}$ satisfying 
$(\ref{K44})$. The difficulty is that this solution may not be larger than $u$ for $|x|=R$. In such a case, for $a>\max u\lfloor_{\prt B_R}$, the function $\tilde u_a:=\tilde u+a$ is a supersolution of $(\ref{Z1})$ in $B_R^c$. The solution $\tilde v$ which satisfies $\CL_p\tilde v=0$ in $B_R^c$ with value $a$ for $|x|=R$ is a subsolution smaller than $\tilde u_a$. Hence there exists a radial solution $u_a$ of $(\ref{Z1})$ in $B_R^c$ such that $u_a(R)=a+\tilde u(R)$ and it dominates $u$ in $B_R^c$.  By \rth{T8}-(1), the function $u_a$, and therefore $u$, satisfies $(\ref{K44})$.\smallskip

\nind In the second case, We proceed as in the proof of \rth{ploplus}-(2), with the help of \rth{T8}-(2). The function $u$ satisfies 
$u_1\leq u\leq u_2$ where $u_1$ and $u_2$ are radial solutions of  $(\ref{Z1})$ in $B_R^c$, hence $u(x)\leq C|x|^{2-N}$. If we set $w(t,s)=r^\gn u(r,s)$ with 
$t=\ln r$ and $\gn=N-2$, then $w$ satisfies 
\bel{K47}\BA{lll}
w_{tt}-\gn w_t+\Gd' w+\gn^q Me^{(N-q(N-1))t}\left((w_t-\gn w)^2+|\nabla 'w|^2\right)^\frac q2\\[2mm]\phantom{------------}
-e^{(N-p(N-2))t}w^p=0\qquad\text{in }[0,\infty)\ti S^{N-1}.
\EA\ee 
Since $w$, $w_t$ and $|\nabla' w|$ are bounded, it follows by the same argument \cite[Proposition 4.1]{BiRa} that $w(t,.)\to c\geq 0$ when $t\to\infty$ and $c>0$ since $u$ is bounded from below by the solution of $\CL^pv=0$ in $B_R^c$ which satisfies the same type of asymptotic behaviour with a positive limit of 
$r^{N-2}v(r)$ when $r\to\infty$.  \smallskip

\nind In the third case, it is proved in \cite[Th\'eor\`eme 3.1]{Veasym} that the solution $v$ of $\CL_pv=0$ in $B_R^c$ which coincide with $\min u\lfloor_{\prt B_R}$ for $|x|=R$ verifies the relation $(\ref{K46})$. In order to have the estimate from above, for $a>0$ the function $h_a(x)=a|x|^{2-N}$ satisfies 
$$\CL^M_{p,q}h_a(x)=a^p|x|^{-p(N-2)}-Ma^q(N-2)^q|x|^{-q(N-1)}=a^p|x|^{-N}-Ma^q(N-2)^q|x|^{-q(N-1)}.
$$
Since $q>\frac{2p}{p+1}=\frac{N}{N-1}$, we obtain for $|x|\geq R$
$$\BA {lll}\CL^M_{p,q}h_a(x)=a^q|x|^{(1-N)q}\left(a^{p-a}|x|^{(N-1)q-N}-M(N-2)^q\right)\\[2mm]
\phantom{\CL^M_{p,q}h_a(x)}\geq a^q|x|^{(1-N)q}\left(a^{p-a}R^{(N-1)q-N}-M(N-2)^q\right),
\EA$$
Therefore, for $a>0$ large enough $\CL^M_{p,q}h_a\geq 0$ in $B^c_R$. Since the solution $v_a$ of $\CL_pv=0$ in $B^c_R$ with value $aR^{2-N}$ for $|x|=R$ is a subsolution of $(\ref{Z1})$ smaller than $h_a$ it follows from \cite[Theorem 1.4.5]{Vebook} that there exists a  radial solution $u_a$ of $(\ref{Z1})$ in 
$B_R^c$ such that $v_a\leq u_a\leq h_a$. If we choose $a$ large enough so that $aR^{2-N}\geq \max u\lfloor_{\prt B_R}$, then $u_a$ is larger than $u$ in 
$B_R^c$. Finally, by \rth{T8}-(3) the function $u_a$ satisfies also $(\ref{K46})$, which ends the proof.\qeda\medskip

In the following result we extend  \rth{T9}-(2)  to the non-radial case.

\bth{T9bis} Let $N\geq 3$, $M>0$,  $p>1$ and $1<q<\frac{2p}{p+1}$ with $q\leq\frac{N}{N-1}$. If $u$ is a positive solution of $(\ref{Z1})$ in $\overline B_R^c$ 
then there holds
\bel{K48}\BA{lll}
\displaystyle \lim_{|x|\to\infty}|x|^\gg u(x)=X_{_M}.
\EA\ee 
\es
\Proof We can assume that $u$ is continuous in $B_\gr^c$ for any $\gr>R$.\\
For constructing a supersolution we proceed as in \rth{20M} using the solution $\gu$ of $(\ref{Z1})$ on $(0,\infty)$. The function 
$x\mapsto \gu(|x|-\gr)+\gd$ is a supersolution $\ref{Z1}$ in $B_\gr^c$ which is larger than $u$ for $|x|=\gr$ and at infinity. Hence it is larger than $u$ in $\overline B_\gr^c$. Letting $\gd\to 0$ yields $u(x)\leq \gu(|x|-\gr)$ for all $x>\gr$. \\

Next, we construct a subsolution: we set $c=\min u\lfloor_{\prt B_\gr}$. 
For $n>\gr$ we denote by $w_n$ the solution of 
$$\BA{lll}-\Gd w+w^p=0\quad\text{in }\, \Gg_{\gr,n}:=B_n\setminus\overline B_\gr\\
w_n=c\quad\text{on }\prt B_\gr\,,\; w_n=0\quad\text{on }\prt B_n.
\EA$$
The function $w_n$ which is unique is a subsolution of $\ref{Z1}$ in $ \Gg_{\gr,n}$ where it satisfies $w_n\leq \gu(|x|-\gr)$. By \rth{BMPLV}  there exists a solution $u_n$ of $(\ref{Z1})$ in $ \Gg_{\gr,n}$ which coincides with $w_n$ on $\prt\Gg_{\gr,n}$ and is radial as $w_n$ and $\gu(|.|-\gr)$ are (or by uniqueness), and $u_n\leq u$ in $ \Gg_{\gr,n}$. When $n\to\infty$, $u_n\uparrow u_\infty$. The function $u_\infty$ is a radial positive solution of $(\ref{Z1})$ in 
$\overline B_\gr^c$ and it satisfies 
$$u_\infty(x)\leq u(x)\leq \gu(|x|-\gr)\quad\text{for all }\, |x|>\gr.
$$
It follows from  \rth{T9}-(2) that $u_\infty$ satisfies $(\ref{K48})$, as $\gu$ does it too. This ends the proof.\qeda

\appendix

\mysection{Appendix}
\subsection{The a priori estimates} 
The following a priori estimates proved in \cite{BVGHV3}, \cite{BVGHV4} are fundamental throughout the paper. They do not depend on the sign of $q-\frac{2p}{p+1}$.
\bprop{Oss} Let $\Gw\subset\BBR^N$ be a domain containing $0$, $1<q<p$ and $M\in\BBR$. If $u\in C^1(\Gw\setminus\{0\})$ is a nonnegative solution of $(\ref{Z1})$ in $\Gw\setminus\{0\}$, then there holds for any $0<R\leq \frac 12\dist(x,\prt\Gw)$:\\
\nind 1-  If $M> 0$, 
\begin{equation}\label{Y1}
u(x)\leq c_1\max\left\{M^{\frac {1}{p-q}}|x|^{-\gg}, |x|^{-\ga}\right\}\quad\text{for all }x\in \overline B_{R}\setminus\{0\},
\end{equation}
where $c_1=c_1(N,p,q)>0$.\\
\nind 2-  If $M\leq 0$ and $q<2$ if $M<0$,
\begin{equation}\label{Y2}\displaystyle
u(x)\leq \min\left\{c_2|x|^{-\ga},c_3|M|^{-\frac {1}{q-1}}|x|^{-\gb}+\max_{|y|=R}u(y) \right\}\quad\text{for all }x\in \overline B_{R}\setminus\{0\},
\end{equation}
where $c_2=c_2(N,p)>0$  and  $c_3=c_3(N,q)>0$.
\es

When $u$ is a signed solution the following estimate holds \cite[Corollary 2.2]{BVGHV4}.
\bprop{Oss-} Under the assumptions on $\Gw$, $p$ and $q$ of \rprop{Oss} and assuming that $M>0$, any signed solution $u$ of $(\ref{K6_})$ in $\Gw\setminus\{0\}$ satisfies for any $0<R\leq \frac 12\dist(x,\prt\Gw)$,
\begin{equation}\label{Y1'}\BA{lll}
-\min\left\{c_4M^{\frac {1}{q-1}}|x|^{-\gb}, c_2|x|^{-\ga}\right\}\leq-u_-(x)\leq 0\\[3mm]
\phantom{--------------}\leq u_+(x)\leq c_1\max\left\{M^{\frac {1}{p-q}}|x|^{-\frac{q}{p-q}}, |x|^{-\ga}\right\},
\EA\end{equation}
for all $x\in \overline B_{R}\setminus\{0\}$, where $c_1=c_1(N,p,q)>0$, $c_2=c_2(N,p)>0$ and  $c_4=c_2(N,p,q)>0$.
\es

Using scaling method when $1<q\leq 2$ and the Bernstein method when $1<q<p$, it is proved in \cite[Proposition 2.3, Corollary 2.5]{BVGHV3} a gradient estimate that we recall. 
\bprop{Osgrad} Let $\Gw\subset\BBR^N$ be a domain containing $0$, $1<q<p$ and $M>0$. If $u\in C^1(\Gw\setminus\{0\})$ is a nonnegative solution of $(\ref{Z1})$ in $\Gw\setminus\{0\}$, there holds for any $0<R\leq R_0:=\frac 12\dist(x,\prt\Gw)$ and some constants $c_4=c_4(N,p,q,R_0)>0$:\\
1- When $1<q\leq \frac{2p}{p+1}$
\begin{equation}\label{Y3}
|\nabla u(x)|\leq c_4\max\left\{M^{\frac {1}{p-q}}|x|^{-\frac{p}{p-q}}, |x|^{-\frac{p+1}{p-1}}\right\},
\end{equation}
for all $x\in \overline B_{R}\setminus\{0\}$. If $q=\frac{2p}{p+1}$, $c_4$ is independent of $R_0$. \\
2- When $1<q<p$
\begin{equation}\label{Y3'}
|\nabla u(x)|\leq c_5\left(|x|^{-\frac{1}{q-1}}+\max\left\{M^{\frac {p}{q(p-q)}}|x|^{-\frac{p}{p-q}}, |x|^{-\frac{2p}{q(p-1)}}\right\}\right),
\end{equation}
for all $x\in \overline B_{R}\setminus\{0\}$. 
\es\medskip

\nind \Remark When $\frac{2p}{p+1}\leq q\leq 2$ any nonnegative solution of $(\ref{Z1})$ in $\BBR^N\setminus B_{\frac R2}$ satisfies $(\ref{Y3})$ in $\BBR^N\setminus B_{R}$.

\nind \Remark If $u$ is a signed solution, the Bernstein method that we developed in \cite{BVGHV3}, \cite{BVGHV4} cannot be applied, however the scaling method can be used if $1<q<p$. In the particular case $1<q\leq \frac{2p}{p+1}$ there holds in a neighborhood of $x=0$, 
\begin{equation}\label{Y3''}
|\nabla u(x)|\leq c_5\max\left\{M^{\frac {1}{p-q}}|x|^{-\gg-1}, |x|^{-\ga-1}\right\}.
\end{equation}

\subsection{Equilibrium with a simple eigenvalue}
Consider the system in $\BBR^2$
\bel{Ap1}\BA {lll}
x_t=ax+by+f(x,y)\\
y_t=cx+dy+g(x,y)
\EA\ee
where $a,b,c,d$ are real numbers and $f$ and $g$ two $C^1$ real functions satisfying 
\bel{Ap2}\BA {lll}
|f(x,y)+g(x,y)|\leq c\left(|x|^s+|y|^s\right)\qquad\text{for all }(x,y)\in B_1.
\EA\ee
Suppose that the matrix $A=\begin{pmatrix}a &b \\ c & d\end{pmatrix}$ of the system admits two eigenvalues $\gm_1\neq\gm_2$ and $\gm_2>0$. By reduction to the diagonal form 
$\begin{pmatrix}x \\ y\end{pmatrix}=P\begin{pmatrix}\tilde x \\ \tilde y\end{pmatrix}$ the system becomes
\bel{Ap3}\BA {lll}
\tilde x_t=\gm_1 \tilde x+\tilde f(\tilde x,\tilde y)\\
\tilde y_t=\gm_2 \tilde y+\tilde g(\tilde x,\tilde y).
\EA\ee
\blemma{eigen}
Under the above assumptions there exist at least two trajectories $\tilde \CT_1=(\tilde x_1,\tilde y_1)$ and $\tilde \CT_2=(\tilde x_2,\tilde y_2)$ tangent to the axis $0\tilde y$ converging to $(0,0)$ when $t\to-\infty$, one with $\tilde y_1(t)>0$, the other with $\tilde y_2(t)<0$ for $t\leq -T$. Any trajectory  $\tilde \CT=\{(\tilde x,\tilde y)\}_{t\leq T}$ converging to 
$(0,0)$ when $t\to-\infty$ and tangent at $(0,0)$ to the axis $0\tilde y$ satisfies for some $\ell\neq 0$
\bel{Ap4}\BA {lll}\displaystyle
\lim_{t\to-\infty}e^{-\gm_2t}\tilde y(t)=\ell.
\EA\ee
\es
\Proof The existence of the solutions tangent to the axis $0\tilde y$ is classical. Consider a $(\tilde x,\tilde y)$ converging to $(0,0)$ tangentialy to 
$0\tilde y$ and such that $\tilde y(t)>0$ for $t\leq -T$. Then
$$|\tilde g(\tilde x(t),\tilde y(t))|\leq c\left(|\tilde x(t)|^s+|\tilde y(t)|^s\right)\leq c'|\tilde y(t)|^s
$$ 
with $c'>c$, since $\frac{\tilde x(t)}{\tilde y(t)}\to 0$ when $t\to-\infty$. Put $\gu(t)=e^{-\gm_2t}\tilde y(t)$. Then 
$$|\gu_t(t)|=|\tilde g(\tilde x(t),\tilde y(t))|e^{-\gm_2t}\leq ce^{-\gm_2t+s\gm_2t}\gu^s(t).
$$
Therefore 
\bel{Ap5}\left|\left(\frac{\gu^{1-s}}{s-1}\right)_t\right|\leq ce^{\gm_2(s-1)t}\Longrightarrow \left(\gu^{1-s}\right)_t\in L^1(-\infty,T)
\ee
Then $\gu^{1-s}(t)$ admits a nonnegative limit $\ell$ when $t\to-\infty$. If $\ell=0$, it would follow from $(\ref{Ap5})$ that
$$\left(e^{\gm_2t}\gu(t)\right)^{s-1}\geq \tilde c\Longrightarrow \tilde y(t)\geq C>0,
$$
which contradict the fact that $\tilde y(t)\to 0$.\qeda\medskip

\nind\Remark This result is easily extendable to higher dimension, where $A$ is a $N\ti N$ matrix with a simple eigenvalue $\gm>0$ and such that 
$\BBR^N=ker(A-\gm I)\oplus E$, where $E$ is A-invariant. Consider the system 
\bel{Ap6}
X'=AX+F(X)
\ee 
where $|F(X)|\leq c|X|^s$ in $B_1$ for some $s>1$. If $X=X_1+X'$ where $X_1\in ker(A-\gm I)$ and $X'\in E$, then there exist two tajectories $X_j(t)$ of 
$(\ref{Ap6})$ admitting a limit direction $\gt\in ker(A-\gm I)\setminus\{0\}$  and , j=1,2, and they satisfy for some $a\neq 0$,
\bel{Ap7}
\displaystyle\lim_{t\to-\infty}e^{-\gm t}X(t)=a\gt.
\ee 



\end{document}